\documentclass[11pt]{amsart} 

\usepackage{amsthm} 
\usepackage{amsmath} 
\usepackage{amssymb} 
\usepackage{latexsym} 
\usepackage{enumerate} 
 
\def\limproj{\mathop{\oalign{{\rm lim}\cr 
\hidewidth$\longleftarrow$\hidewidth\cr}}} 
 
\def\Z{{\mathbb Z}_p} 
\def\Q{{\mathbb Q}_p} 
\def\F{{\mathbb F}_p} 
\def\g{{\rm G}_{\Q}} 
\def\i{{\rm I}_{\Q}} 
\def\G{{\rm GL}_2(\Q)} 
\def\K{{\rm GL}_2(\Z)} 
\def\B{{\rm B}(\Q)} 
\def\BK{{\rm B}(\Z)} 
\def\O{{\mathcal O}_E} 
 
\def\ra{\rightarrow} 
 
\def\RR{{\mathbb R}} 
\def\Qp{{\mathbb Q}_p} 
\def\Cp{\mathbb{C}_p} 
\def\eps{\varepsilon} 
\def\OO{\mathcal{O}}

\def\aa{\mathbb{A}} 
\def\bb{\mathbb{B}} 
\def\ee{\mathbb{E}} 
 
\def\ge{{\rm G}_{\ee}} 
 
\def\at{\widetilde{\mathbb{A}}} 
\def\bt{\widetilde{\mathbb{B}}} 
\def\et{\widetilde{\mathbb{E}}} 
 
\def\aplus{\mathbb{A}^+} 
\def\bplus{\mathbb{B}^+}

\def\atplus{\widetilde{\mathbb{A}}^+} 
\def\btplus{\widetilde{\mathbb{B}}^+} 
\def\etplus{\widetilde{\mathbb{E}}^+} 
 
\def\aminus{\mathbb{A}^-} 
\def\bminus{\mathbb{B}^-} 
 
\newcommand{\btdag}[1]{\widetilde{\mathbb{B}}^{\dagger #1}} 
\newcommand{\bdag}[1]{\mathbb{B}^{\dagger #1}} 
\newcommand{\bhol}[1]{\mathbb{B}^+_{\mathrm{rig} #1}} 
\newcommand{\btrig}[2]{\widetilde{\mathbb{B}}^{\dagger #1}_{\mathrm{rig} #2}} 
\newcommand{\bnrig}[2]{\mathbb{B}^{\dagger #1}_{\mathrm{rig} #2}}

\def\bcris{\mathbb{B}_{\rm cris}} 
\def\bcontp{\widetilde{\mathbb{B}}^+_{\rm rig}} 
\def\bst{\mathbb{B}_{\rm st}} 
\def\bdR{\mathbb{B}_{\rm dR}} 
 
\def\dcris{D_{\rm cris}} 
\def\ddR{D_{\rm dR}} 
\renewcommand{\ddag}[1]{D^{\dagger #1}} 
 
\def\on{\operatorname} 
 
\newtheorem{definit}{\bf{Definition}}[subsection] 
\newtheorem{prop}[definit]{\bf{Proposition}} 
\newtheorem{lem}[definit]{\bf{Lemma}} 
\newtheorem{thm}[definit]{\bf{Theorem}} 
\newtheorem{cor}[definit]{\bf{Corollary}} 
\newtheorem{conj}[definit]{\bf{Conjecture}} 
\newtheorem{rem}[definit]{\bf{Remark}} 
\newtheorem{ex}[definit]{\bf{Example}} 
 
\title{Towards a $p$-adic Langlands programme} 
\author{Laurent Berger \& Christophe Breuil}  
\date{Hangzhou, 16-27 August 2004} 
 
\begin{document} 
 
\maketitle 
 
{\bf WARNING: These notes are informal and we apologize for the inaccuracies, flaws and English
  mistakes that they surely contain.} 
 
\setcounter{tocdepth}{2} 
\tableofcontents 
 
\section{Aim and contents of the course (C.B.)} 
 
\subsection{Aim of the course} 
In 1998, M. Harris and R. Taylor (\cite{HT}), and slightly later
G. Henniart (\cite{He}), proved the celebrated local Langlands
correspondence for ${\rm GL}_n$ (see the course of Bushnell-Henniart
for the case $n=2$, which was known much earlier).
 
Let $p$ be a prime number and $K$, $E$ two finite extensions of $\Q$
with $E$ ``big'' (in particular containing $K$). Recently, Schneider
and Teitelbaum introduced a theory of locally analytic and continuous
representations of $p$-adic groups such as ${\rm GL}_n(K)$ over
$E$-vector spaces (see their course and also \cite{Em2}).
 
The aim of this course is to present some ongoing research on the
(quite ambitious) project of looking for a ``Langlands type local
correspondence'' between, on the one hand, some continuous
representations of ${\rm W}_K:={\rm Weil}(\overline\Q/K)$ on
$n$-dimensional $E$-vector spaces and, on the other hand, some
continuous representations of ${\rm GL}_n(K)$ on infinite-dimensional
$E$-vector spaces (although what we are going to study here is rather
microscopic with respect to such a programme!). Here and throughout,
$E$ will always denote the coefficient field of the vector spaces on
which the groups act. For some technical reasons, we will keep in this
course $E$ finite over $\Q$ but will always tacitly enlarge it if
necessary (which explains the above ``big''). 
 
\begin{ex} 
The case $n=1$ is just local class field theory. Let ${\rm
  W}_K\rightarrow E^{\times}$ be a continuous $1$-dimensional
representation, then, using the (continuous) reciprocity isomorphism
$\tau_K:{\rm W}_K^{\rm ab}\buildrel\sim\over\rightarrow K^{\times}$,
we deduce a continuous representation $K^{\times}\rightarrow
E^{\times}$.  
\end{ex} 
 
We would like to extend the case $n=1$ to $n>1$. One problem is that
the category of $p$-adic ( = continuous) finite dimensional
representations of ${\rm W}_K$ or ${\rm G}_K:={\rm
  Gal}(\overline\Q/K)$ is BIG, even for $n=2$. Fontaine has defined
important subcategories of $p$-adic representations of ${\rm G}_K$
called crystalline, semi-stable and de Rham with crystalline
$\subsetneq$ semi-stable $\subsetneq$ de Rham. They are important
because they are related to algebraic geometry. Moreover, these
categories are now well understood. Therefore it seems natural to
first look for their counterpart, if any, on the ${\rm GL}_n$-side.
  
In this course, we will do so for $n=2$, $K=\Q$ and crystalline
representations. As in the $\ell$-adic case (i.e. the usual local
Langlands correspondence), we will have an ``$F$-semi-simplicity''
assumption. We will also need to assume that the {\it Hodge-Tate
  weights} of the $2$-dimensional crystalline representation are
distinct. Finally, mainly for simplicity, we will also add the small
assumption that the representation is {\it generic} (which means that
we do not want to be bothered by the trivial representation of $\G$;
the non-generic cases behave a bit differently because of this
representation). If $f$ is an eigenform of weight $k\geq 2$ on
$\Gamma_1(N)$ for some $N$ prime to $p$, then the $p$-adic
representation $\rho_f\!\mid_{\g}$ of $\g$ associated to $f$ is known
to be crystalline generic with Hodge-Tate weights $(0,k-1)$ and is
conjectured to be always ``$F$-semi-simple'' (and really is if $k=2$),
so these assumptions should not be too restrictive. The case of equal
Hodge-Tate weights correspond to weight $1$ modular forms, that we
thus disregard.

Our basic aim in this course is to define and start studying some {\it
$p$-adic Banach spaces} endowed with a continuous action of $\G$,
together with their locally analytic vectors. These Banach will be
associated to $2$-dimensional ``$F$-semi-simple'' generic crystalline
representations of $\g$ with distinct Hodge-Tate
weights. Independently of the quest for a ``$p$-adic Langlands
programme'', we hope that these Banach spaces will also be useful in
the future for a ``representation theoretic'' study of the Iwasawa
theory of modular forms (see e.g. the end of \S\ref{finiso}).

\subsection{Why Banach spaces?} 
We will focus on Banach spaces rather than on locally analytic
representations, although there {\it are} some very interesting
locally analytic vectors inside the Banach (and we will study
them). In order to explain why, let us first recall, without proof,
some results in the classical $\ell$-adic case, due to
M.-F. Vign\'eras.

Recall that, if $\ell$ is any prime number and $E$ is a finite
extension of ${\mathbb Q}_{\ell}$, an $E$-Banach space is a
topological $E$-vector space which is complete with respect to the
topology given by a norm (we don't specify the norm in the data of a
Banach space). 
 
\begin{definit} 
Let $G$ be a group acting $E$-linearly on an $E$-vector space $V$. We
say a norm $|\!|\cdot|\!|$ on $V$ is invariant (with respect to $G$)
if $|\!|g\cdot v|\!|=|\!|v|\!|$ for all $g\in G$ and $v\in V$. 
\end{definit} 

\begin{definit} 
Let $G$ be a topological group. We call a unitary $G$-Banach space on
$E$ any $E$-Banach space $\Pi$ equipped with an $E$-linear action of
$G$ such that the map $G\times \Pi\rightarrow \Pi$, $(g,v)\mapsto g\cdot v$
($v\in\Pi$, $g\in G$) is continuous and such that the topology on $\Pi$ is
given by an invariant norm. 
\end{definit} 
 
Unitary $G$-Banach spaces on $E$ form an obvious category where the
morphisms are the continuous $E$-linear $G$-equivariant maps.
 
Now, let us assume $\ell\ne p$ so that $K$ is $p$-adic and $E$ is
$\ell$-adic. 

\begin{definit} 
A smooth representation of ${\rm GL}_n(K)$ on an $E$-vector space
$\pi$ is said to be integral if there exists an invariant norm on
$\pi$ (with respect to ${\rm GL}_n(K)$). 
\end{definit} 
 
It is easy to see that the absolutely irreducible smooth integral
representations of ${\rm GL}_n(K)$ are exactly those absolutely
irreducible smooth representations of ${\rm GL}_n(K)$ that correspond,
under the classical local Langlands correspondence, to those
$F$-semi-simple $\ell$-adic representations of ${\rm W}_K$ that {\it
  extend} to ${\rm G}_K$ (recall that an $\ell$-adic representation of
$W_K$ on a finite dimensional $E$-vector space is the same thing by
Grothendieck's $\ell$-adic monodromy theorem as a smooth
representation of the Weil-Deligne group on the same vector space).

Let $\pi$ be a smooth integral representation of ${\rm GL}_n(K)$ on an
$E$-vector space and let $|\!|\cdot|\!|$ be an invariant norm on
$\pi$. Then the completion of $\pi$ with respect to $|\!|\cdot|\!|$ is
obviously a unitary ${\rm GL}_n(K)$-Banach space on $E$. Conversely,
let $\Pi$ be a unitary ${\rm GL}_n(K)$-Banach space and denote by
$\pi\subset \Pi$ its subspace of smooth vectors (i.e. vectors fixed by
a compact open subgroup of ${\rm GL}_n(K)$). Then $\pi$ is obviously a
smooth integral representation of ${\rm GL}_n(K)$.

\begin{thm}[M.-F. Vign\'eras]\label{marie} 
The functor $\Pi\mapsto \pi$ that associates to a unitary ${\rm
  GL}_n(K)$-Banach space its subspace of smooth vectors induces an
equivalence of categories between the category of unitary ${\rm
  GL}_n(K)$-Banach spaces on $E$ which are topologically of finite
length and the category of finite length smooth integral
representations of ${\rm GL}_n(K)$ on $E$. An inverse functor
$\pi\mapsto \Pi$ is provided by the completion with respect to any
invariant norm on $\pi$. Moreover, length$(\Pi)$= length$(\pi)$. 
\end{thm} 
 
In particular, all the invariant norms on an integral representation
$\pi$ of finite length induce the same topology (in fact, Vign\'eras
has proved that the lattices $\{v\in\pi\mid |\!|v|\!|\leq 1\}$ for an
invariant norm $|\!|\!\cdot\!|\!|$ are always finitely generated over
$\O[\G]$, and thus are all commensurable). Moreover $\pi$ is
irreducible if and only $\Pi$ is topologically irreducible. Thus, we
see that the classical local Langlands correspondence can be
reformulated as a correspondence between ($F$-semi-simple) $\ell$-adic
representations of ${\rm G}_K$ and unitary ${\rm GL}_n(K)$-Banach
spaces on finite extensions of ${\mathbb Q}_{\ell}$ that are
topologically absolutely irreducible. Hence, we {\it could} work with
Banach spaces in the $\ell$-adic case, although everybody agrees that
smooth representations are much more convenient to handle than Banach
spaces. But we know, thanks to Th.\ref{marie}, that it wouldn't make
any difference.

We look for an analogous picture in the $p$-adic case, i.e. where $E$
is a finite extension of $\Q$ and where the smooth vectors inside the
unitary Banach spaces are replaced by the locally analytic
vectors. However, there is no straightforward generalization of
Th.\ref{marie} because:\\  
1) a strongly admissible locally analytic representation of ${\rm
  GL}_n(K)$ on $E$ (see the course of Schneider and Teitelbaum) can
have infinitely many non-equivalent invariant norms (ex: one can prove
this is the case for the locally algebraic representation ${\rm
  Sym}^{k-2}E^2\otimes_E{\rm Steinberg}$ of $\G$ if $k>2$)\\  
2) thanks to a theorem of Schneider and Teitelbaum (\cite{ST4}), we
know that if $\Pi$ is a unitary ${\rm GL}_n(K)$-Banach on $E$ that is
admissible (as a representation of some compact open subgroup) with
locally analytic vectors $\pi$, then length$(\Pi)\leq$ length$(\pi)$
(in fact, unitarity is useless here). However, in general, this is
only a strict inequality (ex: the continuous Steinberg is
topologically irreducible whereas the locally analytic Steinberg has
length $2$).
 
So it seems natural in the $p$-adic case to rather focus on Banach
spaces since the number of Jordan-H\"older factors is smaller. Also,
in the crystalline $\G$-case, we will see that the unitary $\G$-Banach
spaces are somewhat closer to the Galois representations than their
locally analytic vectors. 

\subsection{Some notations and contents} 
Here are some notations we will use: $p$ is a prime number,
$\overline\Q$ an algebraic closure of the field $\Q$ of $p$-adic
rationals, $\overline\Z$ the ring of integers in $\overline\Q$,
$\overline\F$ its residue field and $\g$ the Galois group ${\rm
  Gal}(\overline\Q/\Q)$. We denote by $\varepsilon : \g^{\rm
  ab}\rightarrow \Z^{\times}$ the $p$-adic cyclotomic character, ${\rm
  unr}(x)$ the unramified character of $\Q^{\times}$ sending $p$ to
$x$ (whatever $x$ is), val the $p$-adic valuation normalized by ${\rm
  val}(p)=1$ and $|\cdot|=\frac{1}{p^{{\rm val}(\cdot)}}$ the $p$-adic
norm. The inverse of the reciprocity map $\Q^{\times}\hookrightarrow
\g^{\rm ab}$ is normalized so that $p$ is sent to a geometric
Frobenius. In particular, for $x\in \Q^{\times}$, we have
$\varepsilon(x)=x|x|$ via this map. We let $\B\subset \G$ be the
subgroup of upper triangular matrices. We will loosely use the same
notation for a representation and its underlying $E$-vector space. We
will use without comment basics of $p$-adic functional analysis
(\cite{Schn}) as well as Schneider and Teitelbaum's theory of locally
analytic representations and the structure theorem of locally analytic
principal series for $\G$ (\cite{ST1},\cite{ST2}).

Here is a rough description of the contents of the course.

In lectures $2$ to $5$, the first author will introduce Fontaine's
rings $\bcris$ and $\bdR$ and will define crystalline and de Rham
representations of $\g$ using the theory of {\it weakly admissible
  filtered modules} (this restricted setting will be enough for our
purposes). He will give the explicit classification of $F$-semi-simple
crystalline representations of $\g$ in dimension $1$ and $2$ (in
higher dimension, it becomes more involved). Then he will give a {\it
  second} construction of crystalline representations of $\g$ using
the theory of $(\varphi,\Gamma)$-modules and he will describe the non
trivial links between $p$-adic Hodge theory and
$(\varphi,\Gamma)$-modules. He will then introduce and study the {\it
  Wach module} of a crystalline representation of $\g$ (a ``smaller''
module than the $(\varphi,\Gamma)$-module). These Wach modules will be
useful in the sequel for the applications to $\G$.

In lectures $6$ to $9$, the second author will define and start
studying the unitary $\G$-Banach spaces associated to $2$-dimensional
crystalline representations of $\g$ satisfying the above assumptions
(genericity, $F$-semi-simplicity and distinct Hodge-Tate weights). He
will first introduce the necessary material on distributions and
functions on $\Z$ (functions of class ${\mathcal C}^r$, tempered
distributions of order $r$, etc.) and will then turn to the definition
of the unitary $\G$-Banach space $\Pi(V)$ associated to the
($2$-dimensional) crystalline representation $V$. When $V$ is
reducible, $\Pi(V)$ is admissible, has topological length $2$ and is
split (as a representation of $\G$) if and only if $V$ is split (as a
representation of $\g$). When $V$ is absolutely irreducible, it seems
hard to prove anything on $\Pi(V)$ directly. For instance, it is not
even clear in general on its definition that $\Pi(V)\ne 0$. To study
$\Pi(V)$ for $V$ irreducible, one surprisingly needs to use the theory
of $(\varphi,\Gamma)$-modules for $V$. More precisely, in that case we
prove that there is a ``canonical'' (modulo natural choices)
$E$-linear isomorphism of vector spaces: 
$$\big(\limproj_{\psi}D(V)\big)^{\rm b}\simeq \Pi(V)^*$$ 
where $D(V)$ is the $(\varphi,\Gamma)$-module attached to $V$,
$\psi:D(V)\twoheadrightarrow D(V)$ is a certain canonical surjection
(that will be described), $\rm b$ means ``bounded sequences'' and
$\Pi(V)^*$ is the Banach dual of $\Pi(V)$. Such an isomorphism
(i.e. the link with $(\varphi,\Gamma)$-modules) was first 
found by Colmez (\cite{Co3}) in the case $V$ is irreducible
semi-stable non crystalline (and was inspired by work of the second
author \cite{Br3}, \cite{Br4}) but we won't speak of these cases
here. The left hand side of the above isomorphism being known to be
non zero, we deduce $\Pi(V)^*\ne 0$ and thus $\Pi(V)\ne 0$. The left
hand side being ``irreducible'' and ``admissible'' (what we mean here
is explained in the text), we also deduce that $\Pi(V)$ is
topologically irreducible and admissible as a $\G$-Banach space. It is
hoped that this isomorphism will be also useful in the future for a
more advanced study of the representations $\Pi(V)$, e.g. of their
locally analytic vectors and of their reduction modulo the maximal
ideal of $E$.

Finally, in the last lecture, the first author will state a conjecture
(contained in \cite{Br2}) relating the semi-simplifications of
$\Pi(V)$ and $V$ modulo the maximal ideal of $E$ and will prove in
some cases this conjecture using his Wach modules and previous
computations of the second author.

We didn't include all the proofs of all the results stated or used in
this text. But we have tried to include as many proofs, or sketches of
proofs, or examples, or references, as our energy enabled us to. We
apologize for the proofs that perhaps should be, and are not, in this
course.

\section{$p$-adic Hodge theory (L.B.)}\label{LB2} 
The aim of this lecture is to explain what crystalline representations 
are and to classify them (at least in dimensions $1$ and $2$). 
Let us start with 
a $p$-adic representation $V$, that is a finite dimensional $E$-vector 
space endowed with a continuous linear action of the group $\g$.  
For the time being, we will consider the underlying $\Q$-vector  
space of $V$ and $p$-adic representations will be $\Q$-vector spaces. We 
want to single out certain subcategories of the category of all 
$p$-adic representations. There is a general 
strategy for cutting out subcategories of the category of $p$-adic 
representations. Suppose that we are given a topological $\Q$-algebra 
$B$ which is equipped with a continuous linear action of $\g$. We say 
that $V$ is $B$-admissible if $B \otimes_{\Q} V = B^{\dim_{\Q}(V)}$ as 
$B[\g]$-modules. The set of all $B$-admissible representations is then 
a subset of the set of all $p$-adic representations. Furthermore, if 
$B$ is equipped with extra structures which commute with the action of 
$\g$ then $D_B(V):=(B \otimes_{\Q} V)^{\g}$ inherits those structures 
and provides us with invariants which we can use to classify 
$B$-admissible representations.  
In this lecture, we will define a 
ring of periods $\bcris$ which has two extra structures: a 
filtration and a Frobenius map $\varphi$. Hence, a 
$\bcris$-admissible representation $V$ gives rise to a filtered 
$\varphi$-module $\dcris(V)$ and we will see that one can  
actually recover the representation $V$ from $\dcris(V)$. 
 
\subsection{Some rings of periods} 
Let $\Cp$ be the completion of $\overline\Qp$  
for the $p$-adic topology and let 
\[ \et=\limproj_{x\mapsto x^p} \Cp  
=\{ (x^{(0)},x^{(1)},\cdots) \mid (x^{(i+1)})^p = x^{(i)} \}, \] 
and let $\etplus$ be the set of $x \in \et$  
such that $x^{(0)} \in  \OO_{\Cp}$.  
If $x=(x^{(i)})$ and $y=(y^{(i)})$ are two elements of $\et$, 
we define their sum $x+y$ and their product $xy$ by: 
\[ (x+y)^{(i)}= \lim_{j \ra + \infty} (x^{(i+j)}+y^{(i+j)})^{p^j}  
\quad\text{and}\quad 
(xy)^{(i)}=x^{(i)}y^{(i)}, \] which  
makes $\et$ into a field of 
characteristic $p$. If $x=(x^{(n)})_{n \geq 0}  
\in \et$, let $v_E(x)=v_p(x^{(0)})$. 
This is a valuation on $\et$ for which $\et$ is  
complete; the ring of integers of $\et$ is 
$\etplus$.  
 
Let $\epsilon=(\epsilon^{(i)}) \in \etplus$ where  
$\epsilon^{(0)}=1$ and $\epsilon^{(i)}$ is a primitive  
$p^i$th root of $1$.  
It is easy to see that  
$\F((\epsilon-1)) \subset \et = \etplus[(\epsilon-1)^{-1}]$ and one can show that  
$\et$ is a field which is the completion of 
the algebraic (non-separable!) closure of $\F((\epsilon-1))$, so it is 
really a familiar object.  
 
Let $\atplus$ be the ring $W(\etplus)$  
of Witt vectors with coefficients in 
$\etplus$ and let \[ \btplus=\atplus[1/p]=\{  
\sum_{k\gg -\infty} p^k [x_k],\ x_k \in \etplus 
\} \] where $[x] \in \atplus$ is the Teichm{\"u}ller  
lift of $x \in \etplus$.   
The topology of $\atplus$  
is defined by taking the collection of open sets  
$\{ ([\overline{\pi}]^k,p^n)\atplus \}_{k,n \geq 0}$  
as a family of neighborhoods of $0$, so that the natural map $\prod 
\etplus \ra \atplus$ which to $(x_k)_k$ assigns $\sum p^k [x_k]$ is a 
homeomorphism when $\etplus$ is given its valued field topology. 
This is not the $p$-adic topology of $\atplus$, which makes the map 
$(x_k)_k \mapsto \sum p^k [x_k]$ homeomorphism when $\etplus$ is given 
the discrete topology. 
  
The ring $\btplus$ is 
endowed with a map $\theta: \btplus \ra \Cp$ defined  
by the formula  
\[ \theta\left(\sum_{k\gg -\infty} p^k[x_k] \right) 
=\sum_{k\gg -\infty} p^k x_k^{(0)}. \] 
The absolute Frobenius $\varphi: \etplus \ra \etplus$ lifts by 
functoriality of Witt vectors 
to a map $\varphi: \btplus \ra \btplus$. It's easy to see that  
$\varphi(\sum p^k [x_k])=\sum p^k [x_k^p]$ and that $\varphi$ is 
bijective.  
 
Recall that $\epsilon=(\epsilon^{(i)})_{i \geq 0} 
\in\etplus$ where $\epsilon^{(i)}$  
is as above,  
and define $\pi=[\epsilon]-1$, $\pi_1= 
[\epsilon^{1/p}]-1$, $\omega=\pi/\pi_1$ 
and $q=\varphi(\omega)=\varphi(\pi)/\pi$.  
One can easily show that $\ker(\theta : \atplus \ra \OO_{\Cp})$ 
is the principal ideal generated by $\omega$. 
 
Here is a simple proof: obviously, the kernel of  
the induced map 
$\theta: \etplus \ra \OO_{\Cp}/p$ 
is the ideal of $x \in \etplus$ such that $v_E(x) \geq 1$. Let $y$ be any element of 
$\atplus$ killed by $\theta$ whose reduction modulo $p$ satisfies $v_E(\overline{y})=1$. 
The map $y \atplus \ra \ker(\theta)$ is then injective, and surjective modulo $p$; since 
both sides are complete for the $p$-adic topology, it is an isomorphism. Now, we just need 
to observe that the element $\omega$ is killed by $\theta$  
and that $v_E(\overline{\omega})=1$. 
 
\subsection{The rings $\bcris$ and $\bdR$}\label{2.2} 
Using this we can define $\bdR$;  
let $\bdR^+$ be the ring obtained by completing 
$\btplus$ for the $\ker(\theta)$-adic topology, so that 
\[ \bdR^+=\limproj_n \btplus / (\ker(\theta))^n. \]  
In particular, since $\ker(\theta)=(\omega)$, 
every element $x \in \bdR^+$ 
can be written (in many ways) as a sum $x=\sum_{n=0}^{+\infty} x_n \omega^n$ with $x_n \in 
\btplus$.  
The ring $\bdR^+$ is then naturally a $\Q$-algebra.  
Let us construct an interesting 
element of this ring; since $\theta(1-[\epsilon])=0$, the element 
$1-[\epsilon]$ is ``small''  
for the topology of $\bdR^+$ and the 
following series \[ -\sum_{n = 1}^{+\infty} \frac{(1-[\epsilon])^n}{n} \] will converge in 
$\bdR^+$, to an element which we call $t$. Of course, one should think of $t$ as 
$t=\log([\epsilon])$.  
For instance, if $g \in \g$, then  
\[ g(t)=g(\log([\epsilon]))=\log([g(\epsilon^{(0)},\epsilon^{(1)},\cdots)])= 
\log([\epsilon^{\eps(g)}])=\eps(g)t \] 
so that $t$ is a ``period'' for the cyclotomic character $\eps$.   
 
We now set $\bdR=\bdR^+[1/t]$, which is a field that  
we endow with the filtration defined by $\on{Fil}^i \bdR = t^i 
\bdR^+$. This is the natural filtration on $\bdR$ coming from the fact 
that $\bdR^+$ is a complete discrete valuation ring. 
By functoriality, all the rings we have defined are equipped with a continuous 
linear action of $\g$. One can show that $\bdR^{\g}=\Q$,  
so that if $V$ is a $p$-adic 
representation, then $\ddR(V)=(\bdR \otimes_{\Qp} V)^{\g}$ is naturally a  
filtered $\Q$-vector space.  
Note that $\dim_{\Q} (\ddR(V)) \leq d = \dim_{\Qp}(V)$ in general. 
\begin{definit}\label{defdr} 
We say that $V$ is  
\emph{de Rham} if $\dim_{\Q} \ddR(V) = d$. 
\end{definit} 
 
If $V$ is a de Rham representation, the \emph{Hodge-Tate weights} of 
$V$ are the integers $h$ such that $\on{Fil}^{-h} \ddR(V) \neq 
\on{Fil}^{-h+1} \ddR(V)$. The multiplicity of $h$  
is $\dim \on{Fil}^{-h} \ddR(V) / 
\on{Fil}^{-h+1} \ddR(V)$ so that $V$ has $d$ Hodge-Tate weights. 
 
One unfortunate feature of $\bdR^+$ is that it is too coarse a ring: there is 
no natural extension of the natural Frobenius 
$\varphi:\btplus \ra \btplus$ to a continuous map $\varphi: \bdR^+ \ra 
\bdR^+$. For example,  
if $\tilde{p} \in \etplus$ is an element such that 
$\tilde{p}^{(0)}=p$, then 
$\theta([\tilde{p}^{1/p}]-p) \neq 0$, so that  
$[\tilde{p}^{1/p}]-p$ is invertible in $\bdR^+$, and so $1/([\tilde{p}^{1/p}]-p) \in 
\bdR^+$. But if 
$\varphi$ is a natural extension of   
$\varphi:\btplus \ra \btplus$, then one should have 
$\varphi(1/([\tilde{p}^{1/p}]-p))=1/([\tilde{p}]-p)$, and since $\theta([\tilde{p}]-p)=0$,  
$1/([\tilde{p}]-p) \notin \bdR^+$.  
One would still like to have a Frobenius map, and there is a natural 
way to complete $\btplus$ (where one avoids inverting 
elements like $[\tilde{p}^{1/p}]-p$) such that the completion is still 
endowed with a Frobenius map. 
 
Recall that the topology of $\btplus$  
is defined by taking the collection of open sets  
$\{ ([\overline{\pi}]^k,p^n)\atplus \}_{k,n \geq 0}$  
as a family of neighborhoods of $0$. 
The ring $\bcris^+$ is defined as being 
\[ \bcris^+= \{ \sum_{n \geq 0} a_n \frac{\omega^n}{n!} 
\text{ where $a_n\in \btplus$ is sequence converging to $0$} \}, \] 
and $\bcris=\bcris^+[1/t]$. 
 
One could replace  
$\omega$ by any generator of $\ker(\theta)$ in $\atplus$. The 
ring $\bcris$ injects canonically into $\bdR$  
($\bcris^+$ is naturally defined as a subset of $\bdR^+$) and,  
in particular, it is endowed with the induced 
Galois action and filtration, as well as  
with a continuous Frobenius $\varphi$, extending the 
map $\varphi: \btplus \ra \btplus$. For example, $\varphi(t)=pt$.  
Let us point out once again that 
$\varphi$ does not extend continuously  
to $\bdR$. One also sets $\bcontp=\cap_{n=0}^{+\infty} 
\varphi^n(\bcris^+)$. 
 
\begin{definit}\label{defcris} 
We say that a representation $V$  
of $\g$ is \emph{crystalline} if  
it is $\bcris$-admissible.  
\end{definit} 
 
Note that $V$ is crystalline if and only if 
it is $\bcontp[1/t]$-admissible  
(the periods of crystalline representations 
live in finite dimensional $\Q$-vector  
subspaces of $\bcris$, stable by $\varphi$, and so in fact 
in $\cap_{n=0}^{+\infty} 
\varphi^n(\bcris^+)[1/t]$); this is equivalent  
to requiring that the $\Q$-vector space 
\[ \dcris(V)=(\bcris \otimes_{\Qp} V)^{\g}=(\bcontp[1/t] 
\otimes_{\Qp} V)^{\g} \]  
be of dimension $d=\dim_{\Qp}(V)$.  
Then $\dcris(V)$ is endowed with a Frobenius $\varphi$  
induced by that of $\bcris$ 
and $(\bdR\otimes_{\Q} V)^{\g}=\ddR(V)= 
\dcris(V)$ so that a crystalline  
representation is also de Rham and $\dcris(V)$ is a 
filtered $\Qp$-vector space.  
 
If $V$ is an $E$-linear representation, then we can do the above 
constructions (of $\dcris(V)$ and $\ddR(V)$) by considering $V$ as a 
$\Qp$-linear representation, and then one can check that $\dcris(V)$ 
and $\ddR(V)$ are themselves $E$-vector spaces and that: (1) $\varphi: 
\dcris(V) \ra \dcris(V)$ is $E$-linear (2) the filtration on $\ddR(V)$ 
is given by $E$-vector spaces. The Hodge-Tate weights of $V$ are then 
by convention taken to be with multiplicity $\dim_E 
\on{Fil}^{-h} \ddR(V) /\on{Fil}^{-h+1} \ddR(V)$ and so $V$ has $\dim_E(V)$ 
Hodge-Tate weights.  
 
Remark: one can also define a ring $\bst$ and study semistable 
representations, but we will not do so in this course. 
 
\subsection{Crystalline representations in dimensions $1$ and $2$} 
It is not true that every filtered $\varphi$-module $D$ over $E$ 
arises as the 
$\dcris$ of some crystalline representation $V$. If $D$ is a filtered 
$\varphi$-module of dimension $1$, and $e$ is some basis,  
define $t_N(D)$ as the valuation of the coefficient of $\varphi$ in that 
basis. Define $t_H(D)$ as the largest integer $h$ such that 
$\on{Fil}^h D \neq 0$. If $D$ is a filtered 
$\varphi$-module of dimension $\geq 1$, define $t_N(D):=t_N(\det D)$ and 
$t_H(D):=t_H(\det D)$. Note that if $D=\dcris(V)$, then $t_H(D)$ is 
minus the sum of the Hodge-Tate weights. 
We say that $D$ is an admissible filtered 
$\varphi$-module if $t_N(D)=t_H(D)$ and if for every subspace $D'$ of $D$ 
stable under $\varphi$, we have $t_N(D')-t_H(D') \geq 0$.  
 
One can prove that if $V$ is a crystalline representation, then 
the filtered $\varphi$-module  
$\dcris(V)$ is admissible and conversely, Colmez and Fontaine proved 
the following (cf \cite{CF00}): 
\begin{thm}\label{colfont} 
If $D$ is admissible, then there exists a crystalline 
representation $V$ of $\g$ such that $\dcris(V)=D$.  
The functor $V \mapsto \dcris(V)$ is then an equivalence of 
categories, from the category of $E$-linear  
crystalline representations to the 
category of filtered $\varphi$-modules over $E$.  
\end{thm} 
 
This allows us to make a list of all crystalline representations: we 
just have to make a list of all possible admissible filtered 
$\varphi$-modules. We will now do this in dimensions $1$ and $2$. Let 
$D = \dcris(V)$ be an admissible filtered $\varphi$-module. 
 
If $\dim_E(D) = 1$ we can write $D = E \cdot e$ 
where $\varphi(e)=p^n \lambda e$ with $\lambda \in \O^*$ and if $h$ is 
the Hodge-Tate weight of $V$, then $\on{Fil}^{-h}(D)=D$ and  
$\on{Fil}^{-h+1}(D)=0$. The admissibility condition in this case is 
nothing more than $n=-h$. We can then explicitly describe the 
representation $V$. Let $\mu_{\lambda}$ be the unramified character of 
$\g$ sending $\sigma$ to $\lambda$. One can then check that $V=\eps^n 
\mu_\lambda$. If $\lambda=1$ and $E=\Qp$,  
then $V=\eps^n$ is called $\Qp(n)$ and 
it has the following property: if $V$ is a de Rham representation with 
Hodge-Tate weights $h_1,\cdots,h_d$, then $V(n):=V \otimes_{\Qp} 
\Qp(n)$  has Hodge-Tate weights $h_1+n,\cdots,h_d+n$. The 
representation $V(n)$ is called a \emph{twist} of $V$. 
 
If $\dim_E(D)=2$  we can once again make a list of 
the possible $D$'s. However, we won't be able to make a list of the 
corresponding $V$'s because in general, it's very hard to give 
``explicitly'' a $p$-adic representation. This is why the theory of 
filtered $\varphi$-modules is so useful. 
 
Given a crystalline representation $V$, we can always twist it by a 
suitable power of the cyclotomic character, so that its Hodge-Tate 
weights are $0,k-1$ with $k \geq 1$, so now we assume this.  
Suppose furthermore that the Hodge-Tate weights of $V$ are distinct, 
so that $k \geq 2$, and that 
$\varphi:D \ra D$ is semi-simple  
(as we said in the introduction, we're 
only interested in those cases).  
 
We can then enlarge $E$ enough so that it contains  
the eigenvalues of $\varphi$, 
which we will call $\alpha^{-1}$ and 
$\beta^{-1}$. We then have $D=E e_\alpha \oplus E e_\beta$ with 
$\varphi(e_\alpha)=\alpha^{-1} e_\alpha$ and $\varphi(e_\beta) = \beta^{-1} 
e_{\beta}$. If $\alpha=\beta$, we choose any basis of $D$.

We have $t_H(D)=-(k-1)$ and $t_N(D)=-\on{val}(\alpha) - 
\on{val}(\beta)$ so if $D$ is admissible, then 
$\on{val}(\alpha)+\on{val}(\beta)=k-1$. But this is not the only 
condition. The filtration on $D$ is given as follows: $\on{Fil}^i D$ 
is $D$ if $i \leq -(k-1)$, it's a line $\Delta = E \delta$  
if $-(k-2) \leq i \leq 
0$ and it's $0$ if $i\geq 1$ (by definition of the Hodge-Tate 
weights). Note in particular that $t_H(e_\alpha)$ and $t_H(e_\beta)$ 
are both either $-(k-1)$ or $0$. 
 
The admissibility condition says that $-\on{val}(\alpha) = 
t_N(e_\alpha) \geq t_H(e_\alpha) \geq -(k-1)$ and similarly for 
$\on{val}(\beta)$. This implies 
that $(k-1) \geq \on{val}(\alpha) \geq 0$ and   
$(k-1) \geq \on{val}(\beta) \geq 0$. Suppose that $\on{val}(\beta) 
\leq \on{val}(\alpha)$ (switch $\alpha$ and $\beta$ if necessary).  
There are three cases to consider. 
\begin{enumerate} 
\item $0 < \on{val}(\beta) \leq \on{val}(\alpha)$.  
 If $\delta$ were $e_\alpha$ or $e_\beta$, then $\Delta$ would be a 
 sub-$\varphi$-module of $D$ and the admissibility condition says that 
 $-\on{val}(\alpha)$ or $-\on{val}(\beta)$ should be $\geq 
 t_H(\delta)=0$ which is a contradiction.  
 Therefore $\delta$ 
  cannot be $e_\alpha$ nor $e_\beta$. 
In particular, this implies that $\alpha \neq \beta$. 
 We can then rescale $e_\alpha$ and 
  $e_\beta$ so that $\delta=e_\alpha+e_\beta$.  
\item $\on{val}(\beta)=0$ and $\on{val}(\alpha)=(k-1)$. Then $\delta$ 
  can be $e_\beta$ or (up to rescaling as previously) 
  $e_\alpha+e_\beta$. Suppose it is $e_\alpha+e_\beta$. Then $E 
  e_\alpha$ is a subobject of $D$ which is itself admissible, so it 
  corresponds to a subrepresentation of $V$ so that $V$ is 
  reducible. It is however non-split, because $D$ has no other 
  admissible subobjects.  
\item $\on{val}(\beta)=0$ and $\on{val}(\alpha)=(k-1)$ and 
  $\delta=e_\beta$. In this case, $D$ is the direct sum of the two 
  admissible objects $E e_\alpha$ and $E e_\beta$ and so 
  $V$ is the direct sum of an unramified representation and another
  unramified representation twisted by
  $\Qp(k-1)$.   
\end{enumerate} 
 
As an exercise, one can classify the remaining cases (when $k=1$ or 
$\varphi$ is not semisimple). Finally, note that in the special case 
$\alpha+\beta=0$, it is possible to give an ``explicit'' description 
of $V$.  
 
To conclude this section, we note that if $V$ is a representation 
which comes from a modular form, whose level $N$ is not divisible by 
$p$, then $V$ is crystalline and we know exactly what $\dcris(V)$ is 
in terms of the modular form. This makes the approach via filtered 
$\varphi$-modules especially interesting for concrete applications.

\section{$(\varphi,\Gamma)$-modules (L.B.)}\label{LB3} 
 
In this lecture, we will explain what $(\varphi,\Gamma)$-modules are. The 
idea is similar to the one from the previous lecture: we construct a 
ring of periods and we use it to classify representations. In this 
case, we can classify all $p$-adic representations, as well as the 
$\g$-stable $\O$-lattices in them. The objects we use for this 
classification (the $(\varphi,\Gamma)$-modules) are unfortunately much 
more complicated than filtered $\varphi$-modules and not very well 
understood. In the next lecture, we will explain the links between the 
theory of filtered $\varphi$-modules and the theory of 
$(\varphi,\Gamma)$-modules for crystalline representations. 
 
Once more, we will deal with $\Qp$-linear representations and then at 
the end say what happens for $E$-linear ones (the advantage of doing 
this is that we can avoid an ``extra layer'' of notations). 
 
\subsection{Some more rings of periods}\label{d(v)} 
Let us go back to the field $\et$ which was defined in the previous 
lecture. As we pointed out then, $\et$ contains $\F((\epsilon-1))$ and is 
actually isomorphic to the completion of $\F((\epsilon-1))^{\rm alg}$. Let 
$\ee$ be the separable closure of $\F((\epsilon-1))$ (warning: \emph{not} 
the completion of the separable closure, which by a theorem of Ax is 
actually equal to $\et$) and write $\ee_{\Q}$ for $\F((\epsilon-1))$ (the 
reason for this terminology will be clearer later on). Write $\ge$ for 
the Galois group $\on{Gal}(\ee/\ee_{\Q})$ (a more correct notation for 
$\ge$ would be ${\rm G}_{\ee_{\Q}}$ but this is a bit awkward). 
 
Let $U$ be an $\F$-representation of $\ge$, that is a finite 
dimensional (=$d$)  
$\F$-vector space $U$ with a continuous linear action of 
$\ge$. Since $U$ is a finite set, this means that $\ge$ acts through 
an open subgroup. We can therefore apply ``Hilbert 90'' which tells us 
that the $\ee$-representation $\ee \otimes_{\F} U$ is isomorphic to 
$\ee^d$ and therefore that 
$D(U):=(\ee \otimes_{\F} U)^{\ge}$ is a $d$-dimensional  
$\ee_{\Q}$-vector space and that $\ee \otimes_{\F} U = \ee 
\otimes_{\ee_{\Q}} D(U)$.  
In addition, $D(U)$ is a \emph{$\varphi$-module over $\ee_{\Q}$}, meaning 
a vector space over $\ee_{\Q}$ with a semilinear Frobenius $\varphi$ such 
that $\varphi^* D(U) = D(U)$. 
  
It is therefore possible to recover $U$ from 
$D(U)$ by the formula $U=(\ee \otimes_{\ee_{\Q}} D(U))^{\varphi=1}$ 
and as a consequence, we get the following standard 
theorem going back at least to Katz: 
 
\begin{thm}\label{nickthecat} 
The functor $U \mapsto D(U)$ is an equivalence 
of categories from the category of $\F$-representations of $\ge$ to 
the category of $\varphi$-modules over $\ee_{\Q}$.  
\end{thm} 
 
Now let $\aa_{\Q}$ be the $p$-adic completion of $\Z[[X]][1/X]$ inside 
$\at$ where we have written $X$ for $\pi=[\epsilon]-1$. We can identify 
$\aa_{\Q}$ with the set of power series $\sum_{i=-\infty}^{\infty} a_i 
X^i$ where $a_i \in \Z$ and $a_{-i} \ra 0$ as $i \ra \infty$. The ring 
$\aa_{\Q}$ is a Cohen ring for $\ee_{\Qp}$. In particular, $\bb_{\Q} = 
\aa_{\Q}[1/p]$ is a field and we let $\bb$ denote the $p$-adic 
completion of the maximal unramified extension of $\bb_{\Qp}$ in $\bt$ 
so that $\aa = \bb \cap \at$ is a Cohen ring for $\ee$.  
 
Let $T$ be a $\Z$-representation of $\ge$, that is a finite type 
$\Z$-module with a continuous linear action of $\ge$. We can lift the 
equivalence of categories 
given by theorem \ref{nickthecat}
to get the following one:  
the functor $T \mapsto D(T) := (\aa \otimes_{\Z} T)^{\ge}$  
is an equivalence of categories from the category of  
$\Z$-representations of $\ge$ to the category of  
$\varphi$-modules over $\aa_{\Q}$.  
 
Finally, by inverting $p$ we get:  
\begin{thm}\label{notnick} 
The functor $V \mapsto D(V) := (\bb \otimes_{\Q} V)^{\ge}$  
is an equivalence of categories from the category of  
$\Q$-representations of $\ge$ to the category of  
\'etale $\varphi$-modules over $\bb_{\Q}$.   
\end{thm} 
 
Notice the appearance of 
``\'etale'' in the last statement. We say that a $\varphi$-module $D$ 
over $\bb_{\Q}$ is \emph{\'etale} if we can write $D = \bb_{\Q} 
\otimes_{\aa_{\Q}} D_0$ where $D_0$ is a $\varphi$-module over 
$\aa_{\Q}$ (which includes the requirement $\varphi^* D_0 = D_0$). 
 
This way, we have a convenient way of studying representations of 
$\ge$ using $\varphi$-modules over various rings. In the next 
paragraph, we will answer the question: what does this have to do with 
representations of $\g$? The answer is given by an amazing 
construction of Fontaine and Wintenberger. 
 
\subsection{The theory of the field of norms} 
Recall that we have chosen a primitive $p^n$th root of unity 
$\epsilon^{(n)}$ for every $n$. Set $F_n=\Q(\epsilon^{(n)})$ and $F_{\infty} = 
\cup_{n \geq 0} F_n$. The Galois group $\Gamma_{\Q}$ of 
$F_{\infty}/\Qp$ is isomorphic to $\Z^*$ (by the cyclotomic 
character). In this paragraph, we will explain a striking relation 
between $H_{\Q}$, the Galois group of $\overline\Qp/F_\infty$ and 
$\ge$ the Galois group of $\ee/\ee_{\Q}$. Indeed, Fontaine and 
Wintenberger proved that we have a natural isomorphism $H_{\Q} \simeq 
\ge$. We will now recall the main points of their construction. 
  
Let $K$ be a finite extension of $\Q$ and let $K_n=K(\epsilon^{(n)})$ as 
above. Let 
$\mathcal{N}_K$ be the set $\limproj_n K_n$ where the transition maps are given by 
$N_{K_n/K_{n-1}}$, so that $\mathcal{N}_K$ is the set of sequences 
$(x^{(0)},x^{(1)},\cdots)$ with $x^{(n)} \in K_n$ and $N_{K_n/K_{n-1}}(x^{(n)}) = 
x^{(n-1)}$. If we define a ring structure on $\mathcal{N}_K$ by 
\[ (xy)^{(n)}= x^{(n)}y^{(n)} \quad\text{and}\quad 
(x+y)^{(n)}=\lim_{m \ra +\infty} N_{K_{n+m}/K_n}(x^{(n+m)}+y^{(n+m)}), \] 
then $\mathcal{N}_K$ is actually a field,  
called the \emph{field of norms} of $K_{\infty}/K$. It is 
naturally endowed with an action of $H_K$.  
Furthermore, for every finite Galois extension 
$L/K$, $\mathcal{N}_L/\mathcal{N}_K$ is a  
finite Galois extension whose Galois group is 
$\on{Gal}(L_{\infty}/K_{\infty})$, and every  
finite Galois extension of $\mathcal{N}_K$ is 
of this kind so that the absolute 
Galois group of $\mathcal{N}_K$ is naturally  
isomorphic to $H_K$. 
 
We will now describe $\mathcal{N}_{\Q}$.  
There is a 
map $\mathcal{N}_{\Q} \ra \et$ which sends $(x^{(n)})$ to $(y^{(n)})$ 
where $y^{(n)}=\lim (x^{(n+m)})^{p^m}$ and it is possible to show that 
this map is a ring homomorphism and  
defines an isomorphism from $\mathcal{N}_{\Q}$ to 
$\ee_{\Qp}$ (the point is that by ramification theory, the 
$N_{K_n/K_{n-1}}$ maps are pretty close to the $x \mapsto x^p$ map).  
As a consequence of this and similar statements for the 
$\mathcal{N}_K$'s, we see that there is a natural bijection between 
separable extensions of $\mathcal{N}_{\Q}$ and separable extensions of 
$\ee_{\Q}$. Finally, we get out of this that $H_{\Q} \simeq \ge$. 
 
Let us now go back to the main results of the previous paragraph: we 
constructed an equivalence of categories from the category of  
$\Z$-representations of $\ge$ to the category of  
$\varphi$-modules over $\aa_{\Q}$ and by inverting $p$ 
an equivalence of categories from the category of  
$\Q$-representations of $\ge$ to the category of  
\'etale $\varphi$-modules over $\bb_{\Q}$. Since $H_{\Q} \simeq \ge$, 
we get equivalences of categories from the category of  
$\Z$-representations (resp. $\Q$-representations) of $H_{\Q}$ to the category of  
$\varphi$-modules over $\aa_{\Q}$ (resp over $\bb_{\Q}$).  
 
Finally, if we start with a $\Z$-representation or $\Q$-representation 
of $\g$, then what we get is a $\varphi$-module over $\aa_{\Q}$ (resp 
over $\bb_{\Q}$) which has a residual action of $\g/H_{\Q} \simeq 
\Gamma_{\Q}$. This is an \'etale $(\varphi,\Gamma)$-module: it is a  
module over $\aa_{\Q}$ (resp over $\bb_{\Q}$) with a Frobenius 
$\varphi$ (required to be \'etale) and an action of $\Gamma_{\Q}$ 
which commutes with $\varphi$. Note that on $\aa_{\Q}$, we have 
\[ \varphi(X)=\varphi([\epsilon]-1)= ([\epsilon^p]-1) 
= (1+X)^p-1 \] and also  
\[ \gamma(X)=\gamma([\epsilon]-1)=([\epsilon^{\eps(\gamma)}]-1) 
=(1+X)^{\eps(\gamma)}-1. \]  
 
\begin{thm}\label{pgmod} 
The functors $T \mapsto D(T)$ and $V \mapsto D(V)$ are  
equivalences of categories from the category of  
$\Z$-representations (resp. $\Q$-representations) of $\g$ to the 
category of \'etale $(\varphi,\Gamma)$-modules  over $\aa_{\Q}$ (resp 
over $\bb_{\Q}$).  
\end{thm} 
 
This theorem is proved in \cite{Fo91}, in pretty much the way 
which we recalled. 
 
\begin{ex} 
For example, let $V=\Qp(r)=\Qp \cdot e_r$. The restriction of $V$ to 
$H_{\Q}$ is trivial, so the underlying $\varphi$-module of $V$ is 
trivial: $D(V)=\bb_{\Q} \cdot e_r$ and $\varphi(e_r)=e_r$. The action of 
$\Gamma_{\Q}$ is then given by $\gamma(e_r) = \eps(\gamma)^r e_r$. In 
the lecture on Wach modules, we will see non-trivial examples of 
$(\varphi,\Gamma)$-modules.  
\end{ex} 
 
To conclude this paragraph, note that if $E$ is a finite extension of 
$\Qp$, and $V$ is an $E$-linear representation (or $T$ is an 
$\O$-representation) then $D(V)$ (resp $D(T)$) is a 
$(\varphi,\Gamma)$-module over $E \otimes_{\Q} \bb_{\Q}$ (resp $\O 
\otimes_{\Z} \aa_{\Q}$) and the resulting functors are again 
equivalences of categories.  
 
\subsection{$(\varphi,\Gamma)$-modules and the operator $\psi$}\label{psy} 
Let $V$ be a $p$-adic representation of $\g$ and let $T$ be a lattice 
of $V$ (meaning that $T$ is a free $\Z$-module of rank 
$d=\dim_{\Q}(V)$ such that $V = \Q \otimes_{\Z} T$ and such that $T$ 
is stable under $\g$). 
 
In this paragraph, we will introduce an operator $\psi: D(T) \ra D(T)$ 
which will play a fundamental role in the sequel. Note that $\varphi: 
\bb \ra \bb$ is not surjective, actually $\bb / 
\varphi(\bb)$ is an extension of degree $p$ of local fields, whose 
residual extension is purely inseparable. This makes it possible to 
define a left inverse $\psi$ of $\varphi$ by the formula: 
\[ \varphi(\psi(x)) = \frac{1}{p} \on{Tr}_{\bb / \varphi(\bb)} (x) \] 
and then some ramification theory shows that $\psi(\aa) \subset 
\aa$. Obviously, $\psi(\varphi(x))=x$ and more generally $\psi(\lambda 
\varphi(x)) = \psi(\lambda)x$,  
and also $\psi$ commutes with the 
action of $\g$ on $\bb$. Therefore, we get an induced map $\psi: D(T) 
\ra D(T)$ which is a left inverse for $\varphi$. One can give a slightly 
more explicit description of $\psi$. Given $y \in D(T)$, there exist 
uniquely determined $y_0,\cdots,y_{p-1}$ such that $y=\sum_{i=0}^{p-1} 
\varphi(y_i) (1+X)^i$ and then $\psi(y)=y_0$. Finally, 
if $y=y(X) \in \aa_{\Qp}$, we can give a slight variant of the definition 
of $\psi$:  
\[ \varphi(\psi(y)) = \frac{1}{p} \sum_{\eta^p=1} y((1+X)\eta-1). \] 
 
The $\Z$-module $D(T)^{\psi=1}$ will be of particular importance to 
us. It turns out to be of independent interest in Iwasawa 
theory. Indeed, it is possible to give explicit maps $h^1_n : 
D(T)^{\psi=1} \ra H^1(F_n,T)$ for all $n \geq 0$, and these maps are 
compatible with each other and corestriction (meaning that 
$\on{cor}_{F_{n+1}/F_n} \circ h^1_{n+1} = h^1_n$) so that we get a map 
$D(T)^{\psi=1} \ra \limproj_n  H^1(F_n,T)$ and this map is, by a 
theorem of Fontaine, an isomorphism. The $\Z[[\Gamma_{\Q}]]$-module  
$\limproj_n  H^1(F_n,T)$ is denoted by $H^1_{\rm Iw}(\Q,T)$ and is 
called the Iwasawa cohomology of $T$. Explaining this further is 
beyond the scope of the course, but it is an important application of 
the theory of $(\varphi,\Gamma)$-modules. One use we will make of 
the previous discussion is the following: for any $T$, $D(T)^{\psi=1} 
\neq \{0\}$ (see corollary \ref{brfoiw} for another application).
Indeed, in Iwasawa theory one proves that $H^1_{\rm 
Iw}(\Q,T)$ is a $\Z[[\Gamma_{\Q}]]$-module whose torsion free part is 
of rank $\on{rk}_{\Z}(T)$.  
 
Note that in a different direction, Cherbonnier and 
Colmez have proved that $D(T)^{\psi=1}$ contains 
a basis of $D(T)$ on $\aa_{\Q}$. We see that therefore, the kernel of 
$1-\psi$ on $D(T)$ is rather large. On the other hand, Cherbonnier and 
Colmez have proved that $D(T)/(1-\psi)$ is very small: indeed, if $V$ 
has no quotient of dimension $1$, then $D(T)/(1-\psi)$ is finite (and 
so $D(V)/(1-\psi)=0$). The reason is that we have an identification 
$D(T)/(1-\psi) = H^2_{\rm Iw}(\Q,T)$. 
We will use this to prove the following 
proposition due to Colmez which will be useful later on: 
 
\begin{prop}\label{psisurj} 
Let $V$ be an $E$-linear representation such that  
$\overline\Qp \otimes_{\Q} V$ has no quotient of dimension $1$. If $P 
\in E[X]$, then $P(\psi) : D(V) \ra D(V)$ is surjective. 
\end{prop} 
 
\begin{proof} 
By enlarging $E$ if necessary, we can assume that $P(X)$ splits 
completely in $E$ and so we only need to show the proposition when 
$P(X)=X-\alpha$ with $\alpha \in E$. There are several cases to 
consider:  
\begin{enumerate} 
\item $\alpha=0$. In this case for any $x \in D(V)$, we can write 
  $x=\psi(\varphi(x))$.  
\item $v_p(\alpha) < 0$. In this case, we use the fact that $D(T)$ is 
  preserved by $\psi$ so that the series 
  $(1-\alpha^{-1}\psi)^{-1}=1+\alpha^{-1}\psi + (\alpha^{-1}\psi)^2 + 
  \cdots$ converges and since $(\psi-\alpha)^{-1} = - \alpha^{-1} 
  (1-\alpha^{-1}\psi)^{-1}$ we're done. 
\item $v_p(\alpha) > 0$. By the same argument as the previous one, 
  $(1-\alpha\varphi)^{-1}$ converges on $D(V)$. If $x \in D(V)$, we can 
  therefore write $\varphi(x)=(1-\alpha \varphi)y$ and taking $\psi$, we get 
  $x=(\psi-\alpha)y$.   
\item $v_p(\alpha) = 0$. Let 
  $\mu_{\alpha}$ be the unramified character sending the frobenius to 
  $\alpha$. Then $D(V)/(\psi-\alpha) = D(V \otimes 
  \mu_{\alpha})/(1-\psi)$ and as we have recalled above, $D(V \otimes 
  \mu_{\alpha})/(1-\psi)=0$ if $\overline\Qp 
  \otimes_{\Q} V$ has no quotient of dimension $1$.  
\end{enumerate} 
\end{proof}

\section{$(\varphi,\Gamma)$-modules and $p$-adic Hodge theory (L.B.)}\label{LB4} 
 
In this chapter, we study the link between $p$-adic Hodge theory (the 
theory of $\bcris$ and $\bdR$) and the theory of 
$(\varphi,\Gamma)$-modules. The latter theory provides a very flexible 
way of studying $p$-adic representations, since to describe $p$-adic 
representations, one only needs to give two matrices over $\bb_{\Q}$ 
(one for $\varphi$ and one for a topological generator of $\Gamma_{\Q}$ 
(which is procyclic if $p \neq 2$)) satisfying some simple 
conditions (which say that $\varphi$ and $\Gamma_{\Q}$ commute and that 
$\varphi$ is \'etale).  
 
\subsection{Overconvergent representations} \label{LB41} 
In order to link $p$-adic Hodge theory  and the theory of 
$(\varphi,\Gamma)$-modules, we will apply the usual strategy: 
construct more rings of periods. Let $x \in \bt$ be given and write $x 
=\sum_{k \gg -\infty} p^k [x_k]$ with $x_k \in \et$. Note that if $y 
\in \et$, then $[y] \in \bdR^+$ but still, there is no reason for the 
series $x =\sum_{k \gg -\infty} p^k [x_k]$ to converge in $\bdR^+$ 
because the $x_k$'s could ``grow too fast'' (recall that $\et$ is
a valued field). We will  
therefore impose growth conditions on the $x_k$'s.  
Choose $r \in \RR_{> 0}$ and define 
\[ \btdag{,r} = \{ x \in \bt,\ x= \sum_{k \gg -\infty} p^k [x_k],\  
k+\frac{p-1}{pr} v_E(x_k) \ra +\infty \}. \] 
If $r_n = p^{n-1}(p-1)$ for some $n \geq 0$,  
then the definition of $\btdag{,r_n}$ boils down   
to requiring that  
$\sum_{k \gg -\infty} p^k x_k^{(n)}$ converge in $\Cp$,  
which in turn is equivalent to requiring that the series
$\varphi^{-n}(x) = \sum_{k \gg -\infty} p^k [x_k^{1/p^n}]$  
converge in $\bdR^+$. For example $X \in \bdag{,r}$ for all $r$'s and 
$\varphi^{-n}(X) = [\epsilon^{1/p^n}]-1 = \epsilon^{(n)} e^{t/p^n} -1 \in 
\bdR^+$.  
 
Let $\bdag{,r} = \btdag{,r} \cap \bb$ and let $\bdag{} = \cup_{r > 0} 
\bdag{,r}$. This is the subring of elements 
of $\bb$ which are ``related to $\bdR$'' in a way. We say that a 
$p$-adic representation $V$ of $\g$ is \emph{overconvergent} if $D(V)$ 
has a basis which is made up of elements of $\ddag{}(V)=(\bdag{} 
\otimes_{\Q} V)^{H_{\Q}}$. If there is such a basis, then by 
finiteness it will actually be in $\ddag{,r}(V)=(\bdag{,r} 
\otimes_{\Q} V)^{H_{\Q}}$ for $r$ large enough.  
 
We therefore need to know which $p$-adic representations are 
overconvergent, and this is given to us by a theorem of Cherbonnier 
and Colmez (cf \cite{CC98}):  
\begin{thm}\label{cherbcol} 
Any $p$-adic representation $V$ of $\g$ is 
overconvergent.  
\end{thm} 
This theorem is quite hard, and we will not comment 
its proof. For de Rham representations, one can give a simpler  
proof, based on the results of Kedlaya. 
 
If $V$ is overconvergent and if there is a basis of $D(V)$ in 
$\ddag{,r}(V)$, then $\ddag{,r}(V)$ is a $\bdag{,r}_{\Q}:= 
(\bdag{,r})^{H_{\Q}}$-module. Warning: $\varphi(\bdag{,r}) \subset 
\bdag{,pr}$ so one has to be careful about the fact that 
$\ddag{,r}(V)$ is not a $\varphi$-module. In any case, we need to know 
what $\bdag{,r}_{\Q}$ looks like. The answer turns out to be very 
nice. A power series $f(X) \in \bb_{\Q}$ belongs to $\bdag{,r}_{\Q}$ 
if and only if it is convergent on the annulus $0 < v_p(z) 
\leq 1/r$ (the fact that $f(X) \in \bb_{\Q}$ then implies that it is 
bounded on that annulus). 
In particular, one can think ``geometrically'' about the 
elements of $\bdag{,r}_{\Q}$ and this turns out to be very 
convenient.  
 
Given $r > 0$, let $n(r)$ be the smallest $n$ such that $p^{n-1}(p-1) 
\geq r$ so that $\varphi^{-n}(x)$ will converge in $\bdR^+$ for all $x 
\in \btdag{,r}$ and $n \geq n(r)$. We then have a map $\varphi^{-n} : 
\ddag{,r}(V) \ra (\bdR^+ \otimes V)^{H_{\Q}}$ and the image of this 
map is in $(\bdR^+ \otimes V)^{H_{\Q}}$ if 
$n \geq n(r)$. This is used by Cherbonnier and Colmez to prove a 
number of reciprocity laws, but it is not enough for our purposes, 
which is to reconstruct $\dcris(V)$ from $\ddag{}(V)$ if $V$ is 
crystalline.  
 
We do point out, however, that  
 
\begin{prop}\label{imgiotan} 
If $V$ is de Rham, then the image 
of the map $\varphi^{-n} : \ddag{,r}(V) \ra (\bdR^+ \otimes 
V)^{H_{\Q}}$ lies in $F_n((t)) \otimes_{\Q} \ddR(V)$ and so if $V$ is 
an $E$-linear de Rham representation, then the image 
of the map $\varphi^{-n} : \ddag{,r}(V) \ra (\bdR^+ \otimes 
V)^{H_{\Q}}$ lies in $(E \otimes_{\Q} F_n((t))) \otimes_{E} 
\ddR(V)$.  
\end{prop} 
 
This will be used later on. 
 
\subsection{A large ring of periods} 
In this paragraph, we take up the task of reconstructing  
$\dcris(V)$ from $\ddag{}(V)$ if $V$ is crystalline. For this purpose, 
we will need a few more rings of periods. Let $\bnrig{,r}{,\Q}$ be  
the ring of power series $f(X) = \sum_{i=-\infty}^{\infty} a_i X^i$ 
where $a_i \in \Q$ and $f(X)$ converges on the annulus $0 < v_p(z) 
\leq 1/r$ (but is not assumed to be bounded anymore).  
For example, $t=\log(1+X)$ belongs to $\bnrig{,r}{,\Q}$ but not to 
$\bdag{,r}_{\Q}$. Let $\bnrig{}{,\Q} =  \cup_{r > 0} 
\bnrig{,r}{,\Q}$. Note that this ring is often given a different name, 
namely $\mathcal{R}_{\Q}$. It is the ``Robba ring''. As Colmez says, 
$\mathcal{R}_{\Q}$ is its ``first name'' and $\bnrig{}{,\Q}$ is its ``last 
name''. Similarly, let $\bdag{}_{\Q}=  \cup_{r > 0} \bdag{,r}_{\Q} = 
(\bdag{})^{H_{\Q}}$. The first name of that ring is 
$\mathcal{E}^{\dagger}_{\Q}$.  
 
The main result which we have in sight is the following:  
 
\begin{thm}\label{rpaeqd} 
If $V$ is a 
$p$-adic representation of $\g$, then \[ \left(\bnrig{}{,\Q}[1/t] 
\otimes_{\bdag{}_{\Q}} \ddag{}(V) \right)^{\Gamma_{\Q}} \] is a 
$\varphi$-module and as a $\varphi$-module, it is isomorphic to 
$\dcris(V)$. If $V$ is crystalline, then in addition we have  
\[ \bnrig{}{,\Q}[1/t] 
\otimes_{\Q} \left(\bnrig{}{,\Q}[1/t] 
\otimes_{\bdag{}_{\Q}} \ddag{}(V) \right)^{\Gamma_{\Q}}  = \bnrig{}{,\Q}[1/t] 
\otimes_{\bdag{}_{\Q}} \ddag{}(V). \] 
\end{thm} 
 
The proof of this result is quite technical, and as one suspects, it 
involves introducing more rings of periods, so we will only sketch the 
proof. Let us mention in passing that there are analogous results for 
semistable representations.  
 
Recall that $\ddag{}(V) = (\bdag{} 
\otimes_{\Qp} V)^{H_{\Q}}$ and that $\dcris(V) = (\bcontp[1/t] \otimes_{\Qp} 
V)^{\g}$ as we have explained in paragraph \ref{2.2}.  
The main point is then to construct a big ring $\btrig{}{}$ which 
contains both $\btdag{}$ (and hence $\bdag{}$) and $\bcontp$. This 
way, we have inclusions  
\begin{equation}\label{inc1}  
\dcris(V) \subset \left( \btrig{}{}[1/t] \otimes_{\Qp} V 
\right)^{\g} 
\end{equation} 
and  
\begin{equation}\label{inc2} \left(\bnrig{}{,\Q}[1/t] 
\otimes_{\bdag{}_{\Q}} \ddag{}(V) \right)^{\Gamma_{\Q}} 
\subset \left( \btrig{}{}[1/t] \otimes_{\Qp} V \right)^{\g}. 
\end{equation} 
We then need to prove that these two inclusions are 
equalities. There are essentially two steps. the first is to use the 
fact that all the above are finite dimensional $\Q$-vector spaces 
stable under $\varphi$ and that $\varphi^{-1}$ tends to ``regularize'' 
functions (think of the fact that $\varphi(\btdag{,r}) = 
\btdag{,rp}$). The ``most regular elements'' of $\btrig{}{}$ are those 
which are in $\bcontp$ which explains why inclusion (\ref{inc1}) is an 
equality.  
 
For inclusion (\ref{inc2}), we first remark that \[ (\btrig{}{}[1/t] 
\otimes_{\Qp} V)^{H_{\Q}} = \btrig{}{,\Q}[1/t] \otimes_{\bdag{}_{\Q}} 
\ddag{}(V), \] where $\btrig{}{,\Q} = (\btrig{}{})^{H_{\Q}}$ and in a 
way,  the ring $\btrig{}{,\Q}$ is the completion of $\cup_{m \geq 0} 
\varphi^{-m}(\bnrig{}{,\Q})$ in $\btrig{}{}$. Inclusion (\ref{inc2}) then 
becomes the statement that: 
\[ \left(\bnrig{}{,\Q}[1/t] 
\otimes_{\bdag{}_{\Q}} \ddag{}(V) \right)^{\Gamma_{\Q}} 
 = \left(\btrig{}{,\Q}[1/t] 
\otimes_{\bdag{}_{\Q}} \ddag{}(V) \right)^{\Gamma_{\Q}},\] 
and the proof of this is a ``$\varphi$-decompletion'' process, similar to 
the one used by Cherbonnier and Colmez to prove that all $p$-adic 
representations are overconvergent, and about which we shall say no 
more.  
 
\subsection{Crystalline representations} 
In the previous paragraph, we explained the ideas of the  
(admittedly technical) proof of the 
facts that if $V$ is a crystalline representation, then  
\[ \dcris(V) = \left(\bnrig{}{,\Q}[1/t] 
\otimes_{\bdag{}_{\Q}} \ddag{}(V) \right)^{\Gamma_{\Q}} \] and  
\[ \bnrig{}{,\Q}[1/t] 
\otimes_{\Q} \dcris(V)  = \bnrig{}{,\Q}[1/t] 
\otimes_{\bdag{}_{\Q}} \ddag{}(V). \] 
In this paragraph, we will explain one important consequence of this 
theorem, namely that $p$-adic representations are of finite 
height. First, let us explain what this means. Let $\aplus = \atplus 
\cap \aa$ and $\bplus = \btplus \cap \bb$. since $\bplus \subset 
\btplus$, we have a natural map $\bplus \ra \bdR^+$ and the rings 
$\bplus$ are ``even better behaved'' than the $\bdag{,r}$'s. In 
addition, they are stable under $\varphi$. We say that a $p$-adic 
representation $V$ is of \emph{finite height} if $D(V)$ has a basis 
made of elements of $D^+(V) := (\bplus \otimes_{\Q} V)^{H_{\Q}}$. 
This last space is a module over $\bplus_{\Qp}$ and one can prove that 
$\bplus_{\Q} = \Q \otimes_{\Z} \Z[[X]]$ which is indeed a very nice 
ring. We will explain the proof of Colmez' theorem (cf \cite{Co99}):  
\begin{thm}\label{crishf} 
Every crystalline representation of $\g$ is of finite height.  
\end{thm} 
This means 
that we can study crystalline representations effectively using $\Q 
\otimes_{\Z} \Z[[X]]$-modules which we will do in the next lecture.  
 
Recall that Cherbonnier and Colmez have proved that $D(V)^{\psi=1}$ 
contains a basis of $D(V)$. The strategy for our proof that if $V$ is 
crystalline, then $V$ is of finite height is to prove that if $V$ is 
crystalline, and $S$ is the set \[ S = X,\ \varphi(X),\ 
\varphi^2(X),\cdots, \]  
then $D(V)^{\psi=1} \subset (\bplus[S^{-1}] \otimes_{\Q} 
V)^{H_{\Q}}$. This will show  
that $D(V)$ has a basis of elements which live in $(\bplus \otimes_{\Q} 
V)^{H_{\Q}}$. This proof is due to Colmez; there is a different proof, 
due to the first author, which relies on a result of Kedlaya. We will 
use the fact (due to Cherbonnier) that $D(V)^{\psi=1} \subset (\bdag{}  
\otimes_{\Q} V)^{H_{\Q}}$.  
 
For simplicity, we will assume that $\varphi$ is semi-simple on 
$\dcris(V)$, and we will consider $V$ as  
an $E$-linear crystalline representation, where $E$ contains the 
eigenvalues of $\varphi$ on $\dcris(V)$. Suppose then that $\dcris(V) = 
\oplus E \cdot e_i$ where $\varphi(e_i)=\alpha_i^{-1} e_i$.  
If $x \in \ddag{}(V)$, 
then we can write $x=\sum x_i \otimes e_i$ with $x_i \in 
E \otimes_{\Q} \bnrig{}{,\Q}[1/t]$ and if $\psi(x)=x$ then, since $\psi$ acts by 
$\varphi^{-1}$ on $\dcris(V)$, we have $x_i \in 
(E \otimes_{\Q} \bnrig{}{,\Q}[1/t])^{\psi=1/\alpha_i}$.  
 
Let $\bhol{,\Q}$ be the ring of power series $f(X)=\sum_{i \geq 0} a_i 
X^i$ which converge on the open unit disk. The first name of that ring 
is $\mathcal{R}_{\Q}^+$. By an argument of $p$-adic analysis,  
which we will give in paragraph \ref{psianal}, we can 
prove that for any $\alpha \in E$, we have 
$(E \otimes_{\Q} \bnrig{}{,\Q}[1/t])^{\psi=\alpha} = (E \otimes_{\Q} \bhol{,\Q} 
[1/t])^{\psi=\alpha}$. We conclude that $D(V)^{\psi=1} \subset \bhol{,\Q} 
[1/t] \otimes_{\Q} \dcris(V)$. Finally, $\bhol{,\Q} \subset \bcontp$ 
so that we have $D(V)^{\psi=1} \subset (\bcontp [1/t] \otimes_{\Q} 
V)^{H_{\Q}}$. Let us recall that $D(V)^{\psi=1}  \subset (\bdag{} \otimes_{\Q} 
V)^{H_{\Q}}$. This means that the periods of the elements of 
$D(V)^{\psi=1}$ live in $\bcontp[1/t] \cap \bdag{}$ and one can easily 
prove (using the definition of all these rings) that $\bcontp[1/t] 
\cap \bdag{} \subset \bplus[S^{-1}]$. Indeed, one should  
think of $\bcontp$ as ``holomorphic (algebraic) functions on the open 
unit disk'' and of $\bdag{}$ as ``bounded (algebraic) functions on 
some annulus''. The statement that $\bcontp[1/t] 
\cap \bdag{} \subset \bplus[S^{-1}]$ then says that a function which 
is both meromorphic with poles at the $\epsilon^{(n)}-1$ and bounded 
toward the external boundary is the quotient of a bounded holomorphic 
function on the disk (an element of $\bplus$) by $\varphi^n(X^k)$ for 
$k,n$ large enough. Therefore, $V$ is of finite height. 
 
In other words, the $(\varphi,\Gamma)$-module of a crystalline 
representation has a very special form: one can choose a basis such 
that the matrices of $\varphi$ and $\gamma \in \Gamma_{\Q}$ have 
coefficients in $\bplus_{\Q}$ (in $E \otimes_{\Q} \bplus_{\Q}$ if $V$ 
is $E$-linear).  
   
\subsection{Eigenvectors of $\psi$}\label{psianal} 
In this paragraph, we investigate further the action of the operator 
$\psi$ on certain rings of power series. 
The goal is to explain the proof of the following proposition:  
\begin{prop}\label{egvpsi} 
If $\alpha \in E^*$, then:  
\[ (E \otimes_{\Q}  
\bnrig{}{,\Q}[1/t])^{\psi=\alpha} = (E \otimes_{\Q} \bhol{,\Q} 
[1/t])^{\psi=\alpha}. \] 
\end{prop} 
  
\begin{proof} 
Since $\psi(t^{-h}f) = p^h t^{-h} \psi(f)$, 
it is enough to prove that \[ (E \otimes_{\Q}  
\bnrig{}{,\Q})^{\psi=\alpha} =  
(E \otimes_{\Q} \bhol{,\Q})^{\psi=\alpha} \] 
whenever $v_p(\alpha)<0$ (if we are given $f$ such that 
$\psi(f)=\alpha f$, then $\psi(t^h f) = p^{-h} \alpha \cdot t^h f$ so to 
make $v_p(\alpha)<0$ we just need to take $h \gg 0$) 
To avoid cumbersome formulas, we will do the case $E=\Q$ where the 
whole argument already appears. 
 
Recall that we have the following  
rings of power series: $\bhol{,\Q}$ and 
$\bnrig{}{,\Q}$ and $\bdag{}_{\Q}$ and $\bplus_{\Q} = \Q \otimes_{\Z} 
\Z[[X]]$ and $\bb_{\Q}$ as well as $\aa_{\Q}$ and $\aplus_{\Q} = 
\Z[[X]]$. Note that each of these rings is stable under $\psi$.  
Let $\aminus_{\Q}$ be the set of $f(X)=\sum a_i X^i \in \aa_{\Q}$ such 
that $a_i = 0$ if $i \geq 0$ (this is not a ring!) so that $\aa_{\Q} = 
\aminus_{\Q} \oplus \aplus_{\Q}$ and let 
$\bminus_{\Q}=\aminus_{\Q}[1/p]$.  
 
Let us first prove that $\psi(\aminus_{\Q}) \subset \aminus_{\Q}$.  
To see this, recall that $q=\varphi(X)/X=1+[\epsilon]+\cdots+[\epsilon]^{p-1}$ so 
that if $\ell \geq 1$, then
$\psi(q^\ell)$ is a polynomial in $[\epsilon]$ of degree $\leq 
\ell-1$ (since $\psi([\epsilon^i])=0$ if $p \nmid i$ and $[\epsilon^{i/p}]$ 
otherwise) and so $\psi(q^\ell)$ is a polynomial in  
$X=[\epsilon]-1$ of degree $\leq \ell-1$. We'll leave it as an exercise to 
show that its constant term is $p^{\ell-1}$ (this will be useful 
later).  
Now observe that 
\[ \psi(X^{-\ell})=\psi(q^\ell / \varphi(X)^\ell) = X^{-\ell} 
\psi(q^\ell), \] 
so that $\psi(X^{-\ell}) \in \aminus_{\Q}$ and consequently 
$\psi(\aminus_{\Q}) \subset \aminus_{\Q}$.  
 
If $f \in \bnrig{}{,\Q}$ then we can write $f=f^-+f^+$ with $f^+ \in 
\bhol{,\Q}$ and $f^- \in \bminus_{\Q}$. In order to prove our main 
statement, and since $\psi$ commutes with $f \mapsto f^\pm$, we 
therefore only need to study the action of $\psi$ on $\bminus_{\Q}$.  
But since $\psi(\aminus_{\Q}) \subset \aminus_{\Q}$, $\psi$ cannot 
have any eigenvalues $\alpha$ with $v_p(\alpha)<0$  
on $\bminus_{\Q}$, so we are done.  
\end{proof} 
 
Exercise: using this, determine all the eigenvalues of $\psi$ on 
$\bminus_{\Q} \cap \bdag{}_{\Q}$.  
 
\subsection{A review of the notation} 
The point of this paragraph is to review the notation for the various 
rings that have been introduced so far. There are some general 
guidelines for understanding ``who is who''. The rings $\bdR$ and 
$\bcris$ do not follow these patterns though. 
 
Most rings are $\ee^*_*$ or $\aa^*_*$ or $\bb^*_*$. An ``$\ee$'' 
indicates a ring of characteristic $p$, an ``$\aa$'' a ring of 
characteristic $0$ in which $p$ is not invertible, and a ``$\bb$'' a 
ring of characteristic $0$ in which $p$ is invertible. 
 
The superscripts are generally $+$ or $\dagger$ (or nothing). A 
``$+$'' indicates objects defined on the whole unit disk (in a sense) 
while a ``$\dagger$'' indicates objects defined on an annulus (of 
specified radius if we have a ``$\dagger,r$'') and nothing generally 
indicates objects defined ``on the boundary''. A ``$+$'' also means that 
$X=\pi$ is not invertible. So a ``$\dagger$'' indicates that $X$ is 
invertible but not too much (as on an annulus).  
 
The subscript ``${\rm rig}$'' denotes holomorphic-like growth conditions 
toward the boundary, while no subscript means that we ask for the 
much stronger condition ``bounded''. The subscript ``$\Q$'' means we 
take invariants under $H_{\Q}$.  
 
A tilde ``$\widetilde{\ \ }$'' means that $\varphi$ is invertible  
(eg in the algebraic closure $\et$ of $\ee_{\Q}$) while 
no tilde means that $\varphi$ is not invertible (eg the separable closure 
$\ee$ of $\ee_{\Q}$). 
 
Finally, as we said, $\bdR$ and $\bcris$ do not follow these patterns, 
in part because they are not defined ``near the boundary of the open 
disk'' like the other ones. It is better to use the ring 
$\bcontp[1/t]$ than the ring $\bcris$ which is linked to crystalline 
cohomology but has no other advantage. The notation $\bcontp$ is 
consistent with our above explanations. Note that this ring is also 
called $\bb_{\rm cont}^+$ by some authors (with inconsistent 
notation).

\section{Crystalline representations and Wach modules (L.B.)}\label{LB5} 
 
This lecture will be a little different from the previous three: 
indeed we have (finally!) defined all of the period rings which we needed 
and constructed most of the objects which describe, in various ways, 
crystalline representations. It remains to construct the Wach module
associated to a crystalline representation $V$, 
and then to prove some technical statements about the action of $\psi$ 
on this module, which will be crucial for the applications to 
representations of $\G$.  
 
\subsection{Wach modules}\label{LB51} 
Recall from the previous lecture that if $V$ is a crystalline 
representation, then it is of finite height, meaning that we can write 
$D(V) = \bb_{\Q} \otimes_{\bplus_{\Q}} D^+(V)$ where $D^+(V) = (\bplus 
\otimes_{\Q} V)^{H_{\Q}}$ is a free $\bplus_{\Q} = \Q \otimes_{\Z} 
\Z[[X]]$-module of rank $d$. If $V$ is an $E$-linear representation, 
then we have that $D^+(V)$ is a free $E \otimes_{\Q}  
\bplus_{\Q} = E \otimes_{\O} \O[[X]]$-module of rank 
$\dim_E(V)$. Once again, we will deal with $\Qp$-linear 
representations in this chapter and at the end explain how everything 
remains true for $E$-linear objects in an obvious way.  
 
The starting point for this paragraph is that all crystalline 
representations are of finite height, but there are finite height 
representations which are not crystalline (for example: a non-integer 
power of the cyclotomic character). Wach has characterized which 
representations $V$ of finite height are crystalline: they are the ones 
for which there exists a $\bplus_{\Q}$-submodule $N$ of $D^+(V)$  
which is stable under $\Gamma_{\Q}$ and such 
that: (1) $D(V) = \bb_{\Q} \otimes_{\bplus_{\Q}} N$ (2) there exists 
$r \in \mathbb{Z}$ such that $\Gamma_{\Q}$ acts trivially on $(N/ 
X\cdot N)(r)$. 
 
Suppose that the Hodge-Tate weights of $V$ are $\leq 0$. One can then 
take $r=0$ above and one can further ask that $N[1/X] = 
D^+(V)[1/X]$. In this case, $N$ is uniquely determined by the above 
requirements and we call it the Wach module $N(V)$ associated to 
$V$. It is then stable by $\varphi$. So to summarize, the Wach module 
$N(V)$ associated to $V$ satisfies the following properties: 
 
\begin{prop}\label{wachprp} 
The Wach module $N(V) \subset D(V)$ has the following properties if 
the weights of $V$ are $\leq 0$: 
\begin{enumerate} 
\item it is a $\bplus_{\Q}$-module free of rank $d$  
\item it is stable 
under the action of $\Gamma_{\Q}$ and $\Gamma_{\Q}$ acts trivially on 
$N(V)/X \cdot N(V)$  
\item $N(V)$ is stable under $\varphi$.  
\end{enumerate} 
\end{prop} 
 
Suppose that the Hodge-Tate weights of $V$ are in the interval 
$[-h;0]$ and 
let $\varphi^*(N(V))$ be the $\bplus_{\Q}$-module generated by 
$\varphi(N(V))$. One can show that $N(V) / \varphi^*(N(V))$ is killed by 
$(\varphi(X)/X)^h$. Recall that we write $q$ for the important element
$\varphi(X)/X$. 
 
Note that the Wach module of $V(-h)$ is simply $X^h N(V) e_{-h}$ 
(here $e_1$ is a basis of $\Qp(1)$ and $e_k=e_1^{\otimes k}$).  
If the Hodge-Tate weights of $V$ are no longer $\leq 0$ it is still 
possible to define the Wach module of $V$ by analogy with the above 
formula: we set $N(V) = X^{-r} N(V(-r)) e_r$  
where $r$ is large enough 
so that $V(-r)$ has negative Hodge-Tate weights. This obviously does 
not depend on the choice of $r$. The module $N(V) \subset D(V)$ is 
then stable under $\Gamma_{\Q}$ and $\Gamma_{\Q}$ acts trivially on 
$N(V)/X N(V)$ but $N(V)$ is no longer stable under $\varphi$. What is 
true is the following: 
\begin{prop}\label{wachcomplem} 
If in prop \ref{wachprp}, 
the Hodge-Tate weights of $V$ are in the interval 
$[a;b]$ with $b$ not necessarily $\leq 0$,  
then $\varphi(X^b N(V)) \subset X^b N(V)$ and $X^b 
N(V) / \varphi(X^b N(V))$ is killed by $q^{b-a}$.  
\end{prop} 
 
We finish this paragraph with a few examples. If $V=\Q(r)$ then 
$N(V)=X^{-r} \bplus_{\Q} e_r$ which means that $N(V)$ is a 
$\bplus_{\Q}$-module of rank $1$ with a basis $n_r$ such that $\varphi(n_r) 
= q^{-r} n_r$ and $\gamma(n_r) = (\gamma(X)/X)^{-r} 
\eps(\gamma)^r n_r$ if $\gamma \in \Gamma_{\Q}$.  
 
Suppose now that $V$ is the representation attached to a supersingular 
elliptic curve (with $a_p=0$), twisted by $\Q(-1)$ to make its weights 
$\leq 0$. This representation is crystalline and 
its Hodge-Tate weights are $-1$ and $0$ 
(and the eigenvalues of $\varphi$ on $\dcris(V)$ are 
$\pm p^{1/2}$).  
The Wach module $N(V)$ is of rank $2$ 
over $\bplus_{\Q}$ generated by $e_1$ and $e_2$ with 
$\varphi(e_1)=q e_2$ and $\varphi(e_2)=  - e_1$. 
The action of $\gamma \in \Gamma_{\Q}$ 
is given by: 
\[ \gamma(e_1)=\frac{\log^+(1+X)}{\gamma(\log^+(1+X))} e_1 \quad\text{and}\quad  
\gamma(e_2)=\frac{\log^-(1+X)}{\gamma(\log^-(1+X))} e_2 \] where 
\[ \log^+(1+X) =\prod_{n \geq 0} \frac{\varphi^{2n+1}(q)}{p}  
\quad\text{and}\quad 
\log^-(1+X) =\prod_{n \geq 0} \frac{\varphi^{2n}(q)}{p},  \] 
so that $t=\log(1+X)=X \log^+(1+X) \log^-(1+X)$.  
 
Finally, we point out that if $V$ is an $E$-linear representation then 
$N(V)$ is a free $E \otimes \O[[X]]$-module satisfying all the 
properties given previously. 
 
\subsection{The weak topology}\label{defweaktopo} 
Until now, we have largely ignored issues of topology, but for the 
sequel it will be important to distinguish between two topologies on 
$D(V)$. There is the \emph{strong topology} which is the $p$-adic 
topology defined by choosing a lattice $T$ of $V$ and decreeing the 
$p^{\ell} D(T)$ to be neighborhoods of $0$. This does not depend on 
the choice of $T$. Then, there is the \emph{weak topology} which 
Colmez has defined for all $p$-adic representations $V$, but which we 
will only define for crystalline ones, using the Wach module. Choose a 
lattice $T$ in $V$ and define the Wach module associated to $T$ to be 
$N(T) = N(V) \cap D(T)$. One can prove that $N(T)$ is a free 
$\Z[[X]]$-module of rank $d$ such that $N(V) = \Q \otimes_{\Z} N(T)$.  
 
For each $k \geq 0$, we define a semi-valuation 
$\nu_k$ on $D(T)$ as follows: if $x \in D(T)$ then $\nu_k(x)$ is 
the largest integer $j \in \mathbb{Z} \cup \{ + \infty\}$  
such that $x \in X^j N(T) + p^k D(T)$. The weak topology on $D(T)$ is 
then the topology defined by the set $\{\nu_k\}_{k \geq 0}$ of all 
those semi-valuations. 
Remark that the weak topology induces on $N(T)$ the $(p,X)$-adic 
topology. The weak  
topology on $D(V)$ is the inductive 
limit topology on $D(V) = \cup_{\ell \geq 0} D(p^{-\ell}T)$.  
This does not depend on the choice of $T$. Concretely, if we  
have a sequence $(v_n)_n$ of elements of $D(V)$, and 
that sequence is bounded for the weak topology, then there is a 
$\g$-stable lattice $T$ of $V$ such that $v_n \in D(T)$ for every $n \geq 
0$ and furthermore, for every $k \geq 0$, there exists $f(k) \in 
\mathbb{Z}$ such that $v_n \in X^{-f(k)} N(T) + p^k D(T)$.  
 
In lecture 8, we will use the weak topology on $D(T)$.  
 
\subsection{The operator $\psi$ and Wach modules}\label{wachpsi} 
In this lecture, we assume that $V$ is a crystalline representation 
with $\geq 0$ Hodge-Tate weights. The goal of this paragraph is to 
show that   
\begin{thm}\label{psionwach} 
If $V$ is a crystalline representation 
with $\geq 0$ Hodge-Tate weights 
and if $T$ is a $\g$-stable lattice in $V$, then 
\[ \psi(N(T)) \subset N(T) \] and more generally, if 
$\ell \geq 1$, then \[ \psi(X^{-\ell} N(T)) \subset p^{\ell-1}  
X^{-\ell} N(T) + X^{-\ell+1} N(T). \] 
\end{thm} 
 
\begin{proof} 
Say that the weights of $V$ are in $[0,h]$. We 
know that $\varphi(X^h N(T)) \subset X^h N(T)$ and  
that $X^h N(T) / \varphi^*(X^h N(T))$ is killed by $q^h$. If $y \in N(T)$, 
then $q^h X^h y \in \varphi^*(X^h N(T)) = \varphi(X^h) \varphi^*(N(T))$ and 
since by definition  $q^h X^h  =\varphi(X^h)$ this shows that $y \in 
\varphi^*(N(T))$ and therefore $N(T) \subset \varphi^*(N(T))$. If $y \in 
\varphi^*(N(T))$ then by definition we can write $y=\sum y_i \varphi(n_i)$ 
with $y_i \in \O[[X]]$ and $n_i \in N(T)$ and $\psi(y) = \sum 
\psi(y_i) n_i \in N(T)$. This proves that $\psi(N(T)) \subset N(T)$.  
 
Now choose $y \in N(T)$ and $\ell \geq 1$. Once again, since $N(T) 
\subset \varphi^*(N(T))$, we can write  $y=\sum y_i \varphi(n_i)$ 
with $y_i \in \O[[X]]$ and $n_i \in N(T)$ so that $\psi(X^{-\ell} y) = 
\sum \psi(X^{-\ell} y_i) n_i$. Therefore, in order 
to prove that \[ \psi(X^{-\ell} N(T)) \subset p^{\ell-1}  
X^{-\ell} N(T) + X^{-\ell+1} N(T), \] 
it is enough to prove that 
\[ \psi(X^{-\ell} \O[[X]]) \subset p^{\ell-1}  
X^{-\ell} \O[[X]] + X^{-\ell+1} \O[[X]], \] 
which we did in \S\ref{psianal}. 
\end{proof} 
 
\subsection{The module $D^0(T)$}\label{herrmodule} 
The goal of this paragraph is to explain the proof of the following 
proposition:  
 
\begin{prop}\label{bh} 
If $V$ is crystalline irreducible and  
$\overline\Qp \otimes_{\Q} V$ 
has no quotient of $\dim 1$ (e.g. $V$ is absolutely irreducible), 
then there exists a \emph{unique}  
non-zero $\Q$-vector subspace $D^0(V)$ of $D(V)$  
possessing an $\Z$-lattice $D^0(T)$ 
which is a compact  
(for the induced weak topology) 
$\Z[[X]]$-submodule  
of $D(V)$ preserved by $\psi$ and $\Gamma_{\Q}$ with $\psi$ surjective. 
\end{prop} 
 
If $V$ has a quotient of $\dim 1$, then such an $\Z[[X]]$-module exists 
but is not unique.  
For example, if $V=\Qp$, the we can take either $\Q \otimes_{\Z} 
\Z[[X]]$ or $X^{-1} \Q \otimes_{\Z} \Z[[X]]$. 
If $V$ is not crystalline, it still exists (this is 
a result of Herr) but in our case, the situation is greatly simplified 
if we assume $V$ to be crystalline. Indeed, in this case $D^0(V)$ is 
free of rank $d$ over $\Q \otimes_{\Z} \Z[[X]]$.  
 
\begin{proof} 
First, we'll prove the existence of such a module. 
Suppose that the Hodge-Tate weights of $V$ are in $[-h;0]$. We can 
always twist $V$ to assume this, because the underlying $\varphi$-module 
of $D(V)$ is unchanged by twists. Choose a $\g$-stable lattice $T$ in $V$. 
Our assumption implies that $N(T) \subset \varphi^*(N(T))$ and that 
$\varphi^*(X^h N(T)) \subset X^h N(T)$ as we've seen previously.  
As a consequence, $\psi(N(T)) \subset N(T)$ and $X^h N(T) \subset 
\psi(X^h N(T))$. Consider the sequence  \[ X^h N(T),\ \psi(X^h N(T)), 
\ \psi^2(X^h N(T)),\cdots \] 
It is increasing and all its terms are  
contained in $N(T)$. Since $\Z[[X]]$ is 
noetherian, this implies that the sequence is eventually constant, and 
so there exists $m_0$ such that if $m \geq m_0$, then $\psi^{m+1}(X^h 
N(T)) = \psi^m(X^h N(T))$. One can then set $D^0_1(T)=\psi^{m_0}(X^h 
N(T))$ and $D^0_1(V) = \Q \otimes_{\Z} 
D^0_1(T)$. 
This proves the existence of one module  
$D^0_1(T)$ (and the associated $D^0_1(V)$) satisfying the 
conditions of the proposition. Before we prove the  
uniqueness statement, 
we want to set $D^0(V)$ to be the largest such module and first we 
need to prove some properties of those modules satisfying the 
conditions of the proposition. 
 
Suppose therefore that some module  
$M(V)$ satisfies the condition of the proposition, that is we can 
write $M(V) = \Q \otimes_{\Z} M(T)$ with $\psi : M(T) \ra M(T)$ 
surjective and $M(T)$ compact in $D(T)$ for the induced weak 
topology. This last property means that for any $k \geq 0$, there 
exists $f(k)$ such that $M(T) \subset X^{-f(k)}N(T) + p^k 
D(T)$. Recall that in proposition \ref{psionwach}, we've proved that 
if $f(k) \geq 1$, then: 
\[ \psi(X^{-f(k)}N(T)) \subset p^{f(k)-1} X^{-f(k)}N(T) + 
X^{-f(k)+1}N(T). \] 
Since $\psi(M(T))=M(T)$, we have \[ M(T)=\psi^m(M(T)) \subset 
p^{m(f(k)-1)} X^{-f(k)}N(T) + X^{-f(k)+1}N(T) + p^k 
D(T) \] and if $m$ is large enough, we finally get that  
$M(T) \subset X^{-f(k)+1}N(T) + p^k D(T)$. What this shows is that we 
can take $f(k)=1$ for all $k \geq 0$ and so $M(T) \subset 
X^{-1} N(T) + p^k D(T)$ for all $k \geq 0$ which implies that $M(T) 
\subset X^{-1} N(T)$. Therefore, any such $M(T)$ is contained in 
$\cap_{m \geq 0} \psi^m(X^{-1}N(T))$.  
 
On the other hand, if we set $D^0(T) = \cap_{m \geq 0} 
\psi^m(X^{-1}N(T))$, then obviously $D^0(T)$ satisfies the conditions 
of the proposition, and it contains $D^0_1(T)$ constructed previously 
so it is certainly non-zero, and $D^0(V)$ is therefore a free 
$\bplus_{\Q}$-module of rank $d$ contained in $X^{-1}N(V)$.   
 
Any other $M(T)$ satisfying the conditions of the proposition is 
included in $D^0(T)$.  
If $M(V)$ was a free $\bplus_{\Q}$-module of rank $<d$, then 
we would have a submodule of $N(V)$ stable under $\psi$ and 
$\Gamma$. One can prove that it would then have to be stable under 
$\varphi$ (this is easy to see in dimension $d=2$, less so in higher 
dimensions) and this would correspond  
to a subrepresentation of $V$, which 
we assumed to be irreducible. Therefore, any such $M(V)$ is a 
$\bplus_{\Q}$-module of rank $d$ contained in $X^{-1}N(V)$.  
 
Recall that by proposition \ref{psisurj}, for any 
polynomial $P$, the map $P(\psi) : D(V) \ra D(V)$  
is surjective. By 
continuity, we can find a set $\Omega \subset D(T)$  
bounded for the weak topology and $\ell \geq 0$ 
such that $P(\psi)(\Omega) \supset p^\ell D^0(T)$.  
This means that for any $k \geq 
0$, there exists $f(k)$ such that \[ p^\ell D^0(T) =  
\psi^m(p^\ell D^0(T)) \subset \psi^m \circ 
P(\psi) (X^{-f(k)} N(T) + p^k D(T)) \]  
for all $m \geq 0$, and just as 
above this proves that  
$p^\ell D^0(T) \subset P(\psi) X^{-1}N(T)$ and then that  
\[ p^\ell D^0(T)\subset P(\psi) \cap_{m \geq 0} \psi^m(X^{-1}N(T)) = 
P(\psi)(D^0(T)), \]  
so that $\psi: D^0(V) \ra D^0(V)$ is surjective. 
 
We now prove the uniqueness of an $M(V)$ satisfying the conditions of 
the proposition. By the previous discussion, $D^0(V)/M(V)$ is a finite 
dimensional $\Q$-vector space stable under $\psi$ and such that the 
operator $P(\psi)$ is surjective for all $P \in \Q[X]$ which is 
impossible unless $D^0(V)=M(V)$. 
\end{proof} 
 
\subsection{Wach modules and $\dcris(V)$}\label{LB54} 
We end this lecture with two ways of recovering $\dcris(V)$ from the 
Wach module $N(V)$ of a crystalline representation $V$. 
 
Let $V$ be a crystalline representation. Recall that $N(V) \subset 
D(V)$ is not necessarily stable under $\varphi$, but that in any case 
$\varphi(N(V)) \subset N(V) [1/q]$. 
We can define a filtration on $N(V)$ in the following way: $\on{Fil}^i 
N(V) = \{ x \in N(V),\ \varphi(x) \in q^i N(V) \}$. The results of 
\ref{LB51} show that if the weights of $V$ are in $[a;b]$, then 
$\on{Fil}^{-b} N(V) = N(V)$.  
 
Note also that since $X$ and $q=\varphi(X)/X$ are coprime, we have a natural 
identification $N(V)/X \cdot N(V) = N(V) [1/q] / X \cdot N(V) [1/q]$ and we endow 
$N(V)/ X \cdot N(V)$ with the induced filtration from $N(V)$, and the 
Frobenius induced from $N(V)/X \cdot N(V) =  
N(V) [1/q] / X \cdot N(V) [1/q]$. This 
makes $N(V)/ X \cdot N(V)$ into a filtered $\varphi$-module. 
The main result is then (see \cite[\S III.4]{LB04}): 
\begin{thm}\label{sdlkfj} 
The two filtered $\varphi$-modules 
$N(V)/ X \cdot N(V)$ and $\dcris(V)$ are isomorphic. 
\end{thm} 
 
We shall now prove a technical result which will be used later 
on. Suppose now that the weights of $V$ are $\geq 0$. We've seen that 
both $N(V)$ and $X^{-1} N(V)$ are stable under $\psi$, so we get a map 
$\psi: X^{-1} N(V) / N(V) \ra X^{-1} N(V) / N(V)$.  
\begin{prop}\label{psidcris} 
If we identify 
$X^{-1} N(V) / N(V)$ with $\dcris(V)$ by  
the map $X^{-1} y \mapsto 
\overline{y}$, then the map  
$\psi: X^{-1} N(V) / N(V) \ra X^{-1} N(V) / N(V)$ 
coincides with the map $\varphi^{-1} : \dcris(V) \ra \dcris(V)$.  
\end{prop} 
 
\begin{proof} 
Indeed, given $y \in N(V)$, choose $z \in N(V)[1/q]$ such that $y-\varphi(z) \in 
X \cdot N(V)[1/q]$. We then have $\psi(y/X)=\psi(\varphi(z)/X+w)$ with $w 
\in N(V)$ and $\psi(\varphi(z)/X)=z/X$ so that $\psi(y/X)\equiv z 
\mod{N(V)[1/q]}$.  
\end{proof} 
 
Now, we will give a different relation between $\dcris(V)$ and 
$N(V)$. Recall that we have an isomorphism 
\[ \bnrig{}{,\Q}[1/t] \otimes_{\Q} \dcris(V) \simeq \bnrig{}{,\Q}[1/t] 
\otimes_{\bdag{}_{\Q}} \ddag{}(V) \]  and in particular an inclusion 
\[ N(V) \subset \ddag{}(V) \subset \bnrig{}{,\Q}[1/t] \otimes_{\Q} \dcris(V). \]  
One can then prove that: 

\begin{prop}\label{inccris}\label{groumpf} 
If $V$ is crystalline, then 
\[ N(V) \subset \bhol{,\Q}[1/t] \otimes_{\Q} \dcris(V), \] 
and \[ \bhol{,\Q}[1/t] \otimes_{\bplus_{\Q}} N(V) 
= \bhol{,\Q}[1/t] \otimes_{\Q} \dcris(V). \]
Furthermore, if the weights of $V$ are $\geq 0$, then  
\[ N(V) \subset \bhol{,\Q} \otimes_{\Q} \dcris(V). \] 
\end{prop} 
 
Exercise: work out explicitly all the above identifications for the 
Wach modules given at the end of \ref{LB51}.  
 
\section{Preliminaries of $p$-adic analysis (C.B.)} 
We give here some classical material of $p$-adic analysis which will
crucially be used in the coming lectures. The results can be found in
\cite{Schi}, \cite{Co1}, \cite{Schi2} or any treatise on $p$-adic
analysis. Thus, we do not provide most of the proofs. 
 
\subsection{Functions of class ${\mathcal C}^r$}\label{c^r} 
 
Fix $E$ a finite extension of $\Q$. We have extensively seen in
Schneider and Teitelbaum's course what a locally analytic function
$f:\Q\rightarrow E$ is. Let $r$ be an element in ${\mathbb
R}^+$. Here, we study the larger class of functions $f:\Z\rightarrow
E$ which are {\it of class ${\mathcal C}^r$}. For $f:\Z\rightarrow E$
any function and $n\in {\mathbb Z}_{\geq 0}$, we first set:
$$a_n(f):=\sum_{i=0}^n(-1)^i{n \choose i}f(n-i).$$ 
 
\begin{definit}\label{classcr} 
A function $f:\Z\rightarrow E$ is of class ${\mathcal C}^r$ if
$n^r|a_n(f)|\rightarrow 0$ in ${\mathbb R}^+$ when $n\rightarrow
+\infty$. 
\end{definit} 
 
Of course, if $f$ is of class ${\mathcal C}^r$, then it is also of
class ${\mathcal C}^s$ for any $s\leq r$. Functions of class
${\mathcal C}^0$ are functions $f$ such that $a_n(f)\rightarrow 0$ in
$E$, which can easily be seen to be equivalent to $f$ continuous (see
e.g. \cite{Co2},\S V.2.1). Hence, a function of class ${\mathcal C}^r$
for $r\in {\mathbb R}^+$ is necessarily at least continuous. Moreover,
one can write (\cite{Co2},\S V.2.1): 
\begin{eqnarray}\label{mahler} 
f(z)=\sum_{n=0}^{+\infty}a_n(f){z\choose n} 
\end{eqnarray} 
where $z\in \Z$ and ${z\choose n}:=\frac{z(z-1)\cdots
  (z-n+1)}{n!}$. The $E$-vector space of functions of class ${\mathcal
  C}^r$ is a Banach for the norm $|\!|f|\!|_r:={\rm
  sup}_n\big((n+1)^{r}|a_n(f)|\big)$ which contains locally analytic
functions on $\Z$, hence also locally polynomial functions. We denote
this Banach by ${\mathcal C}^r(\Z,E)$.

\begin{thm}\label{densitybase} 
Let $d$ be an integer such that $r-1<d$, then the $E$-vector space of
locally polynomial functions $f:\Z\rightarrow E$ of degree at most $d$
is dense in ${\mathcal C}^r(\Z,E)$. 
\end{thm} 
 
This theorem is a consequence of the Amice-V\'elu condition as will be
seen in \S\ref{Yvette} (Cor.\ref{travail}).

We now state results on the spaces ${\mathcal C}^r(\Z,E)$ when $r=n$
is in ${\mathbb Z}_{\geq 0}$ (although this is actually probably a
useless assumption for some of the results below).

For $n\in {\mathbb Z}_{\geq 0}$ and $f:\Z\rightarrow E$ any function,
define $f^{[n]}(z,h_1,\cdots,h_n):\Z\times (\Z-\{0\})^n\rightarrow E$
by induction as follows: 
\begin{eqnarray*}
f^{[0]}(z)&\!\!\!:=\!\!\!&f(z)\\
f^{[n]}(z,h_1,\cdots\!,h_n)&\!\!\!:=\!\!\!&\frac{f^{[n-1]}(z+h_n,h_1,\cdots\!,
h_{n-1})-f^{[n-1]}(z,h_1,\cdots\!,h_{n-1})}{h_n}.
\end{eqnarray*} 
 
\begin{thm}\label{ordreentier} 
A function $f:\Z\rightarrow E$ is of class ${\mathcal C}^n$ if and
only if the functions $f^{[i]}$ for $0\leq i\leq n$ extend
continuously on $\Z^{i+1}$. Moreover, the norm ${\rm sup}_{0\leq i\leq
  n}{\rm sup}_{(z,h_1,\cdots\!,h_i)\in
  \Z^{i+1}}|f^{[i]}(z,h_1,\cdots\!,h_i)|$ is equivalent to the norm
$|\!|\cdot|\!|_n$. 
\end{thm} 
 
See \cite{Co2},\S V.3.2 or \cite{Schi},\S54. Note that, in particular,
any function $f$ of class ${\mathcal C}^n$ admits a continuous
derivative $f':=f^{[1]}(z,0)$.  

\begin{thm} 
Assume $n\geq 1$. The derivative map $f\mapsto f'$ induces a
continuous topological surjection of Banach spaces ${\mathcal
C}^n(\Z,E)\twoheadrightarrow {\mathcal C}^{n-1}(\Z,E)$ which admits
a continuous section.  
\end{thm} 
 
Note that, with our definition, it is not transparent that $f'\in
{\mathcal C}^{n-1}(\Z,E)$ (and indeed, this requires a proof). We
refer to \cite{Schi},\S\S78 to 81. The topological surjection
obviously follows from the existence of a continuous section. Let us
at least explicitly give such a section $P_n:{\mathcal
  C}^{n-1}(\Z,E)\hookrightarrow {\mathcal C}^{n}(\Z,E)$: 
$$P_nf(z):=\sum_{j=0}^{+\infty}\sum_{i=0}^{n-1}
\frac{f^{(i)}(z_j)}{(i+1)!}(z_{j+1}-z_j)^{i+1}$$ 
where $z_j:=\sum_{l=0}^{j-1}a_lp^l$ if $z=\sum_{l=0}^{+\infty}a_lp^l$
with $a_l\in\{0,\cdots,p-1\}$. The fact $P^nf\in {\mathcal
  C}^{n}(\Z,E)$ if $f\in {\mathcal C}^{n-1}(\Z,E)$ is proved in
\cite{Schi},\S81. The fact $P^n$ is continuous is proved in
\cite{Schi},\S79 for $n=1$ and in \cite{Schi2},\S11 for arbitrary $n$.
 
By an obvious induction, we obtain: 
 
\begin{cor}\label{derivsurj} 
The $n$-th derivative map $f\mapsto f^{(n)}$ induces a continuous
topological surjection of Banach spaces ${\mathcal
  C}^n(\Z,E)\twoheadrightarrow {\mathcal C}^{0}(\Z,E)$ which admits a
continuous section.  
\end{cor} 
 
Let $d$ be an integer such that $d\geq n$. We have already seen that
the closure in ${\mathcal C}^n(\Z,E)$ of the vector subspace of
locally polynomial functions of degree at most $d$ is ${\mathcal
  C}^n(\Z,E)$ itself (Th.\ref{densitybase}). We we will need an
analogous result when $d<n$:

\begin{thm}\label{chikof} 
Assume $n\geq 1$ and $d<n$. The closure in ${\mathcal C}^n(\Z,E)$ of
the vector subspace of locally polynomial functions of degree at most
$d$ is the closed subspace of ${\mathcal C}^n(\Z,E)$ of functions $f$
such that $f^{(d+1)}=0$. 
\end{thm} 
 
This is proved in \cite{Schi},\S68 for $n=1$ (and $d=0$) and in
\cite{Schi2},\S8 for $n$ and $d$ arbitrary.  
 
\subsection{Tempered distributions of order $r$}\label{Yvette} 
 
We let ${\mathcal R}_E$ be the Robba ring of power series with
coefficients in $E$ converging on an open annulus and ${\mathcal
  R}_E^+\subset{\mathcal R}_E$ be the subring of power series
converging on the open (unit) disk. We thus have ${\mathcal
  R}_E^+=\{\sum_{n=0}^{+\infty}a_nX^n\mid a_n\in E,\ {\rm
  lim}_n|a_n|r^n=0\ \forall\ r\in [0,1[\}$ that we equip with the
    natural Fr\'echet topology given by the collection of norms ${\rm
      sup}_n(|a_n|r^n)$ for $0<r<1$. We start by recalling Amice's famous result:  
 
\begin{thm}\label{Yvette2} 
The map: 
\begin{eqnarray}\label{Yvette3} 
\mu\mapsto \sum_{n=0}^{+\infty}\mu\Big({z\choose n}\Big)X^n 
\end{eqnarray} 
induces a topological isomorphism between the dual of the $E$-vector
space of locally analytic functions $f:\Z\rightarrow E$ and ${\mathcal
  R}_E^+$.  
\end{thm} 

We now fix $r\in {\mathbb R}^+$. 
 
\begin{definit} 
We define the space of tempered distributions of order $r$ on $\Z$ as
the Banach dual of the Banach space ${\mathcal C}^r(\Z,E)$.  
\end{definit} 
 
People sometimes also say ``tempered distributions of order $\leq
r$''. A tempered distribution of order $r$ being also locally analytic
(as ${\mathcal C}^r(\Z,E)$ contains locally analytic functions), the
space of tempered distributions of order $r$ corresponds by
Th.\ref{Yvette2} to a subspace of ${\mathcal R}_E^+$.  
 
\begin{definit} 
An element $w=\sum_{n=0}^{+\infty}a_nX^n\in {\mathcal R}_E^+$ is of
order $r$ if $n^{-r}|a_n|$ is bounded (in ${\mathbb R}^+$) when $n$
varies.  
\end{definit} 
 
The following corollary immediately follows from Def.\ref{classcr} and
(\ref{mahler}):  
 
\begin{cor}\label{reorderr} 
The map (\ref{Yvette3}) induces a topological isomorphism between
tempered distributions of order $r$ and the subspace of ${\mathcal
  R}_E^+$ of elements $w$ of order $r$ endowed with the norm
$|\!|w|\!|_r:={\rm sup}_n((n+1)^{-r}|a_n|)$.  
\end{cor} 
 
Let $\mu$ be a locally analytic distribution, that is, an element of
the dual of the $E$-vector space of locally analytic functions
$f:\Z\rightarrow E$. If $f:\Z\rightarrow E$ is a locally analytic
function which is $0$ outside of $a+p^n\Z$ ($a\in \Z$, $n$ an
integer), we use the convenient notation
$\int_{a+p^n\Z}f(z)\mu(z):=\mu(f)$. For any $r\in {\mathbb R}^+$, we
denote $|\!|\mu|\!|_r:={\rm sup}_n\big((n+1)^{-r}|\int_{\Z}{z\choose
  n}\mu(z)|\big)\in {\mathbb R}^+\cup \infty$. We define another norm
$|\!|\mu|\!|_r'$ as follows:  
$$|\!|\mu|\!|_r':={\rm sup}_{a\in\Z}{\rm sup}_{j,n\in{\mathbb Z}_{\geq
    0}}
p^{n(j-r)}\left|\int_{a+p^n\Z}(z-a)^j\mu(z)\right|\in {\mathbb
    R}^+\cup \infty.$$  
 
\begin{lem}\label{normequiv} 
The two norms $|\!|\cdot|\!|_r$ and $|\!|\cdot|\!|_r'$ are equivalent (when defined). 
\end{lem} 
\begin{proof} 
We are going to use an intermediate norm. For
$w=\sum_{m=0}^{+\infty}a_mX^m\in {\mathcal R}_E^+$, define:  
$$|\!|w|\!|^{''}_r:={\rm sup}_n\big(p^{-nr}
{\rm sup}_m(|a_m[{\scriptstyle{\frac{m}{p^n}}}]!|)\big)\in {\mathbb R}^+\cup \infty$$ 
where $[\frac{m}{p^n}]$ is the largest integer smaller than
$\frac{m}{p^n}$. Then one can prove that there exists $c_1,c_2\in
{\mathbb R}^+$ such that $c_1|\!|w|\!|_r\leq |\!|w|\!|_r^{''}\leq
c_2|\!|w|\!|_r$ (this is purely an exercise in ${\mathcal R}_E^+$: see
\cite{Co2},Lem.V.3.19 and note that $|[\frac{m}{p^n}]!|$ is
$p^{-\frac{m}{(p-1)p^n}}$ up to a bounded scalar). If
$|\!|\mu|\!|^{''}_r:=|\!|w|\!|^{''}_r$ when $w$ is the Amice transform
of $\mu$, it is thus enough to prove that $|\!|\cdot|\!|^{''}_r$ is
equivalent to $|\!|\cdot|\!|^{'}_r$. Fix $n\in {\mathbb Z}_{\geq 0}$
and let:  
\begin{eqnarray*} 
|\!|\mu|\!|_{r,n}'&:=&p^{-nr}{\rm sup}_{a\in\Z}
{\rm sup}_{j\in{\mathbb Z}_{\geq
    0}}\left|\int_{a+p^n\Z}\Big(\frac{z-a}{p^n}\Big)^j\mu(z)\right|\\  
|\!|\mu|\!|_{r,n}^{''}&:=&p^{-nr}{\rm sup}_{m\in{\mathbb Z}_{\geq 0}}
\bigg(\left|\mu\Big({z\choose m}\left[{\frac{m}{p^n}}\right]!\Big)\right|\bigg) 
\end{eqnarray*}
(so that $|\!|\mu|\!|_{r}'={\rm sup}_n|\!|\mu|\!|_{r,n}'$ and
$|\!|\mu|\!|_{r}^{''}={\rm sup}_n|\!|\mu|\!|_{r,n}^{''}$). Note that
one can replace $a\in \Z$ by $a\in \{1\cdots,p^n\}$ in the definition
of $|\!|\mu|\!|_{r,n}'$. One can prove that each function ${z\choose
  m}[{\frac{m}{p^n}}]!$ can be written as a (finite) linear
combination of the functions ${\bf
  1}_{a+p^n\Z}\big(\frac{z-a}{p^n}\big)^j$ (for $j\in {\mathbb
  Z}_{\geq 0}$ and $a\in \{1\cdots,p^n\}$) with coefficients in
$\Z$. Conversely, one can also write the function ${\bf
  1}_{a+p^n\Z}\big(\frac{z-a}{p^n}\big)^j$ as an (infinite) linear
combination of the functions ${z\choose m}[{\frac{m}{p^n}}]!$ with
coefficients in $\Z$ converging toward $0$ (see \cite{Co2},\S
V.3.1). This implies that $|\!|\mu|\!|_{r,n}'$ is in ${\mathbb R}^+$
if and only if $|\!|\mu|\!|_{r,n}^{''}$ is in ${\mathbb R}^+$ and, if
so, $|\!|\mu|\!|_{r,n}'=|\!|\mu|\!|_{r,n}^{''}$. Likewise with
$|\!|\mu|\!|_{r}^{''}$ and $|\!|\mu|\!|_{r}^{'}$. This finishes the
proof.  
\end{proof} 
 
\begin{lem}\label{etendre} 
Let $\mu$ be a linear form on the space of locally polynomial
functions $f:\Z\rightarrow E$ of degree at most $d\in {\mathbb
  Z}_{\geq 0}$. Assume there exists a constant $C_{\mu}\in E$ such
that $\forall a\in \Z$, $\forall j\in \{0,\cdots,d\}$ and $\forall
n\in {\mathbb Z}_{\geq 0}$:  
\begin{eqnarray}\label{AMcondition}
\int_{a+p^n\Z}(z-a)^j\mu(z)\in C_{\mu}p^{n(j-r)}\O.
\end{eqnarray} 
If $d>r-1$, there is a unique way to extend $\mu$ as a locally
analytic distribution on $\Z$ such that $|\!|\mu|\!|_r'\in {\mathbb
  R}^+$ (i.e. is bounded).  
\end{lem} 
\begin{proof} 
We give the idea of the proof in the case $d=0$ and $r<1$, leaving the
technical details and the general case (which is analogous) to the
reader. We will first define $\int_{a+p^n\Z}(z-a)^j\mu(z)$ by
induction on $j$. For $j=1$ (it is known for $j=0$), we have:  
\begin{eqnarray*} 
\int_{a+p^n\Z}\!\!\!\!\!\!\!\!\!(z-a)\mu(z)=\sum_{b\equiv a(p^n)}
\int_{b+p^{n+1}\Z}\!\!\!\!\!\!\!\!\!(z-a)\mu(z)=\sum_{b\equiv
  a(p^n)}\int_{b+p^{n+1}\Z}\!\!\!\!\!\!\!\!\!(z-b)\mu(z)+\\  
\sum_{b\equiv a(p^n)}(b-a)\int_{b+p^{n+1}\Z}\!\!\!\!\!\mu(z)
=\sum_{b\equiv a(p^n)}\sum_{c\equiv b(p^{n+1})}
\int_{c+p^{n+2}\Z}\!\!\!\!\!\!\!\!\!(z-c)\mu(z)+\\ 
\sum_{b\equiv a(p^n)}\sum_{c\equiv b(p^{n+1})}(b-c)
\int_{c+p^{n+2}\Z}\!\!\!\!\!\!\mu(z)+\sum_{b\equiv a(p^n)}(b-a)
\int_{b+p^{n+1}\Z}\mu(z)=\cdots 
\end{eqnarray*} 
We have $$(b-c)\int_{c+p^{n+2}\Z}\mu(z)\in
p^{n+1}C_{\mu}p^{-r(n+2)}\O=p^{(n+1)(1-r)}p^{-r}C_{\mu}\O$$ and
$(a-b)\int_{b+p^{n+1}\Z}\!\!\mu(z)\in
p^{n(1-r)}p^{-r}C_{\mu}\O$. Decomposing again $c+p^{n+2}\Z$, we see
that we get a converging sum of terms indexed by $N$ of the type
$p^N\int_{\alpha+p^{N+1}\Z}\mu(z)\in p^{N(1-r)}p^{-r}C_{\mu}\O$ with
$N$ growing together with a remainder which is a sum of terms of the
type $\int_{\alpha+p^{N+1}\Z}(z-\alpha)\mu(z)$. Since
$\int_{\alpha+p^{N+1}\Z}(z-\alpha)\mu(z)$ must be $p$-adically very
small if $\mu$ extends thanks to the condition $|\!|\mu|\!|_r'\in
{\mathbb R}^+$ (which in particular implies
$\int_{\alpha+p^{N+1}\Z}(z-\alpha)\mu(z)\in
|\!|\mu|\!|_r'p^{(N+1)(1-r)}\O$), we immediately see that we can
canonically extend $\mu$ on locally polynomial functions of degree at
most $1$ by setting $\int_{a+p^n\Z}\!\!\!(z-a)\mu(z):=$the above
converging sum. Replacing $(z-a)$ by $(z-a)^2$, etc., an easy
induction shows that $\mu$ extends in a similar way to all locally
polynomial functions and by continuity to all locally analytic
functions and moreover that $|\!|\mu|\!|_r'\in {\mathbb R}^+$. Unicity
also follows obviously from the above decomposition and the condition
$|\!|\mu|\!|_r'\in {\mathbb R}^+$.  
\end{proof} 
 
The condition (\ref{AMcondition}) in Lem.\ref{etendre} is sometimes
called the ``Amice-V\'elu condition''.  
 
\begin{cor}\label{travail} 
1) When $d>r-1$, the subspace of locally polynomial functions
$f:\Z\rightarrow E$ of degree at most $d$ is dense in ${\mathcal
  C}^r(\Z,E)$.\\  
2) Let $\mu$ be as in Lem.\ref{etendre} with $d>r-1$, then $\mu$
uniquely extends as a tempered distribution on $\Z$ of order $r$.  
\end{cor} 
\begin{proof} 
1) If it is not dense, there exists a non zero $\mu$ which is tempered
   of order $r$ and vanishes on locally polynomial functions of degree
   at most $d$. Because of Lem.\ref{etendre}, $\mu$ vanishes on all
   locally polynomial functions, hence in particular on ${z\choose n}$
   for any $n$, hence $\mu$ has $0$ as Amice transform which is
   impossible by Th.\ref{Yvette2} since $\mu\ne 0$. 2) The unicity
   follows from 1) and the existence follows from Lem.\ref{etendre},
   Lem.\ref{normequiv} and Cor.\ref{reorderr}.  
\end{proof} 
 
Assuming Th.\ref{chikof} and using Th.\ref{ordreentier}, one can also
prove that, when $r$ is a positive integer and $d<r$, any linear form
$\mu$ as in Lem.\ref{etendre} extends uniquely to a continuous linear
form on the closed subspace of ${\mathcal C}^r(\Z,E)$ of functions $f$
such that $f^{(d+1)}=0$. We leave this as an exercise to the reader.

\section{Crystalline representations of $\G$ (C.B.)} 
We define and start studying unitary $\G$-Banach spaces $\Pi(V)$
associated to the $2$-dimen\-sional crystalline representations $V$ of
$\g$ as in the introduction. When $V$ is reducible, we prove that
$\Pi(V)$ is admissible of topological length $2$.  
 
\subsection{Preliminary Banach spaces}\label{banachdebut} 
 Let $V$ be a crystalline representation of $\g$ on a two-dimensional
$E$-vector space with distinct Hodge-Tate weights. Replacing $V$ by
$V\otimes_E\varepsilon^n$ for some $n\in {\mathbb Z}$, we can assume
that the Hodge-Tate weights of $V$ are $(0,k-1)$ with $k\in {\mathbb
  Z}$, $k\geq 2$. Assuming moreover that $V$ is $F$-semi-simple,
i.e. that the Frobenius $\varphi$ on ${\rm D}_{\rm cris}(V)$ is
semi-simple, we have seen in Lecture \ref{LB2} that there are two
elements $\alpha,\beta\in \O$ such that $\alpha\ne \beta$, ${\rm
  val}(\beta)\leq {\rm val}(\alpha)$, ${\rm val}(\alpha)+{\rm
  val}(\beta)=k-1$ and ${\rm D}_{\rm
  cris}(V)=D(\alpha,\beta)=Ee_{\alpha}\oplus Ee_{\beta}$ with
$\varphi(e_{\alpha})=\alpha^{-1}e_{\alpha}$,
$\varphi(e_{\beta})=\beta^{-1}e_{\beta}$ and:\\

1) $V$ absolutely irreducible: ${\rm val}(\alpha)>0$, ${\rm
   val}(\beta)>0$, ${\rm Fil}^iD(\alpha,\beta)=D(\alpha,\beta)$ if
   $i\leq -(k-1)$, ${\rm
   Fil}^iD(\alpha,\beta)=E(e_{\alpha}+e_{\beta})$ if $-(k-2)\leq i\leq
   0$ and ${\rm Fil}^iD(\alpha,\beta)=0$ if $i>0$,\\  
 
2) $V$ reducible and non-split: ${\rm val}(\alpha)=k-1$, ${\rm
   val}(\beta)=0$, ${\rm Fil}^iD(\alpha,\beta)=D(\alpha,\beta)$ if
   $i\leq -(k-1)$, ${\rm
   Fil}^iD(\alpha,\beta)=E(e_{\alpha}+e_{\beta})$ if $-(k-2)\leq i\leq
   0$ and ${\rm Fil}^iD(\alpha,\beta)=0$ if $i>0$,\\  
 
3) $V$ reducible and split: ${\rm val}(\alpha)=k-1$, ${\rm
   val}(\beta)=0$, ${\rm Fil}^iD(\alpha,\beta)=D(\alpha,\beta)$ if
   $i\leq -(k-1)$, ${\rm Fil}^iD(\alpha,\beta)=Ee_{\beta}$ if
   $-(k-2)\leq i\leq 0$ and ${\rm Fil}^iD(\alpha,\beta)=0$ if $i>0$.\\  
 
The genericity hypothesis alluded to in the introduction just means
that we moreover want $\alpha\ne p\beta$. Note that
$D(\alpha,\beta)\simeq D(\beta,\alpha)$ when ${\rm val}(\alpha)={\rm
  val}(\beta)$.
 
We define now preliminary $\G$-Banach spaces $B(\alpha)$, $B(\beta)$,
$L(\alpha)$, $L(\beta)$ and, when ${\rm val}(\alpha)=k-1$,
$N(\alpha)$. 
 
Let $B(\alpha)$ be the following Banach space. Its underlying
$E$-vector space is the vector space of function $f:\Q\rightarrow E$
such that $f\!\mid_{\Z}$ is of class ${\mathcal C}^{{\rm
    val}(\alpha)}$ and $(\alpha p\beta^{-1})^{{\rm
    val}(z)}z^{k-2}f(1/z)\!\mid_{\Z}$ can be extended as a function of
class ${\mathcal C}^{{\rm val}(\alpha)}$ on $\Z$. As a vector space we
thus have:  
\begin{eqnarray}\label{f1f2} 
B(\alpha)\simeq {\mathcal C}^{{\rm val}(\alpha)}(\Z,E)\oplus 
{\mathcal C}^{{\rm val}(\alpha)}(\Z,E),\ f\mapsto f_1\oplus f_2 
\end{eqnarray} 
where, for $z\in \Z$, $f_1(z):=f(pz)$ and $f_2(z):=(\alpha
p\beta^{-1})^{{\rm val}(z)}z^{k-2}f(1/z)$. Hence $B(\alpha)$ is a
Banach space for the norm:  
$$|\!|f|\!|:={\rm Max}\big(|\!|f_1|\!|_{{\mathcal C}^{{\rm val}(\alpha)}},|\!|f_2|\!|_{{\mathcal C}^{{\rm val}(\alpha)}}\big).$$ 
We endow $B(\alpha)$ with an $E$-linear action of $\G$ as follows: 
\begin{eqnarray}\label{actionjoe} 
\begin{pmatrix}a&b\\c&d\end{pmatrix}\!(f)(z)
=\alpha^{\!-{\rm val}(ad-bc)}(\alpha p\beta^{\!-1})^{\!{\rm
    val}(-cz+a)}(-cz+a)^{\!k-2}
\!f\Big(\frac{dz-b}{-cz+a}\Big) 
\end{eqnarray} 
and to check that this action induces a continuous map $\G\times
B(\alpha)\rightarrow B(\alpha)$ is a simple exercise that we leave to
the reader (recall that by the Banach-Steinhaus Theorem (see
\cite{Schn}) it is enough to check that the map $\G\rightarrow {\rm
  Hom}_E(B(\alpha),B(\alpha))$ is continuous, where the right hand
side is endowed with the weak topology of pointwise
convergence). Likewise, we define $B(\beta)$ with a continuous action
of $\G$. Note that we have a continuous $\G$-equivariant injection:  
\begin{eqnarray}\label{analytique} 
LA(\alpha):=\Big({\rm Ind}_{\B}^{\G}{\rm unr}(\alpha^{-1})
\otimes x^{k-2}{\rm unr}(p\beta^{-1})\Big)^{\rm an}\hookrightarrow B(\alpha) 
\end{eqnarray} 
given by: 
\begin{eqnarray}\label{groupeqp} 
\big(h:\G\rightarrow E\big)\mapsto \big(z\in \Q\mapsto  
h(\begin{pmatrix}0&1\\-1&z\end{pmatrix})\big) 
\end{eqnarray} 
where the left hand side is a locally analytic principal series in the
sense of Schneider and Teitelbaum (\cite{ST1}) with left action of
$\G$ by right translation on functions. We thus see that the Banach
$B(\alpha)$ is nothing else than:  
$$``\big({\rm Ind}_{\B}^{\G}{\rm unr}(\alpha^{-1})\otimes x^{k-2}{\rm
  unr}(p\beta^{-1})\big)^{{\mathcal C}^{{\rm val}(\alpha)}}.''$$  
We have analogous properties with $B(\beta)$ and $LA(\beta)$ by
interchanging $\alpha$ and $\beta$.  
 
When ${\rm val}(\alpha)=k-1$, let $N(\alpha)\subsetneq B(\alpha)$ be
the closed $E$-vector subspace of functions $f$ such that
$f_1^{(k-1)}=f_2^{(k-1)}=0$ ($k-1$-derivative) so that we have an
exact sequence of Banach spaces:  
\begin{eqnarray}\label{deriv} 
0\rightarrow N(\alpha)\rightarrow B(\alpha)\rightarrow {\mathcal C}^{0}(\Z,E)^2\rightarrow 0 
\end{eqnarray} 
where the third map sends $(f_1,f_2)$ to $(f_1^{(k-1)},f_2^{(k-1)})$. 
Note that this is well defined since $f_1$ and $f_2$ are ${\mathcal
  C}^{k-1}$ on $\Z$ (the topological surjection on the right follows
from Cor.\ref{derivsurj}). Now, let us look at the action of
$\G$. When $\chi_1$, $\chi_2$ are two characters on $\Q^{\times}$ with
values in $\O^{\times}$ (i.e. two integral characters), we denote by:  
$$\big({\rm Ind}_{\B}^{\G}\chi_1\otimes \chi_2\big)^{{\mathcal C}^0}$$ 
the Banach of continuous functions $h:\G\rightarrow E$ such that:
$$h(\begin{pmatrix}a&b\\0&d\end{pmatrix}g)=\chi_1(a)\chi_2(d)h(g)$$  
equipped with the topology of the norm $|\!|h|\!|:={\rm Sup}_{g\in
  \K}|h(g)|$ and with an action of $\G$ by right translation. It is
easily seen to be a unitary $\G$-Banach. We can also describe this
Banach as before using (\ref{groupeqp}) in terms of continuous
functions $f:\Q\rightarrow E$ satisfying a certain continuity
assumption at infinity.  
 
\begin{lem}\label{exact} 
Assume ${\rm val}(\alpha)=k-1$ (hence ${\rm val}(\beta)=0$), then the
map $f\mapsto f^{(k-1)}$ induces an exact sequence of $\G$-Banach
spaces on $E$:  
$$0\rightarrow N(\alpha)\rightarrow B(\alpha)\rightarrow \big({\rm
  Ind}_{\B}^{\G}x^{k-1}{\rm unr}(\alpha^{-1})\otimes x^{-1}{\rm
  unr}(p\beta^{-1})\big)^{{\mathcal C}^0}\rightarrow 0.$$  
\end{lem} 
\begin{proof} 
This is a standard and formal computation which is the same as in the
locally analytic context and that we leave to the reader. Note that
the characters in the right hand side parabolic induction are
$\O^{\times}$-valued.  
\end{proof} 
 
Let $L(\alpha)$ (resp. $L(\beta)$) be the following closed subspace of
$B(\alpha)$ (resp. $B(\beta)$): it is the closure of the $E$-vector
subspace generated by the functions $z^j$ and $(\alpha
p\beta^{-1})^{{\rm val}(z-a)}(z-a)^{k-2-j}$ for $a\in \Q$ and $j\in
{\mathbb Z}$, $0\leq j<{\rm val}(\alpha)$ (resp. $0\leq j<{\rm
  val}(\beta)$ with the convention $L(\beta)=0$ if ${\rm
  val}(\beta)=0$). The reader can check using (\ref{actionjoe}) that
$L(\alpha)$ (resp. $L(\beta)$) is preserved by $\G$. In order for
$L(\alpha)$ (resp. $L(\beta)$) to be a subspace of $B(\alpha)$
(resp. $B(\beta)$), it is enough to check (using (\ref{f1f2}) and up
to a translation on $z$):  
 
\begin{lem} 
For $0\leq j<{\rm val}(\alpha)$, the function $z\mapsto (\alpha
p\beta^{-1})^{{\rm val}(z)}z^{k-2-j}$ is of class ${\mathcal C}^{{\rm
    val}(\alpha)}$ on $\Z$. Likewise interchanging $\alpha$ and
$\beta$.  
\end{lem} 
\begin{proof} 
Fix $j$ as in the statement and let $f(z):=(\alpha p\beta^{-1})^{{\rm
    val}(z)}z^{k-2-j}$ which is easily seen to extend to a continuous
function on $\Z$ by setting $f(z)=0$ (look at the valuation of $f(z)$
for $z\ne 0$). Let $f_0$ be the zero function on $\Z$ and, for $n\in
{\mathbb Z}$, $n>0$, consider on $\Z$ the functions $f_n(z):=(\alpha
p\beta^{-1})^{{\rm val}(z)}z^{k-2-j}$ if ${\rm val}(z)<n$ and
$f_n(z):=0$ otherwise. It is clear that $f_n$ is of class ${\mathcal
  C}^{{\rm val}(\alpha)}$ on $\Z$ since it is locally polynomial. If
we can prove that $f_{n+1}-f_n\rightarrow 0$ in ${\mathcal C}^{{\rm
    val}(\alpha)}(\Z,E)$ when $n\rightarrow +\infty$, we deduce that
$\sum_{n=0}^{\infty}(f_{n+1}-f_n)\in {\mathcal C}^{{\rm
    val}(\alpha)}(\Z,E)$ since ${\mathcal C}^{{\rm
    val}(\alpha)}(\Z,E)$ is complete. But this function is clearly $f$
(as can be checked for any $z\in \Z$). If $B$ is any Banach space, we
have a closed embedding $B\hookrightarrow (B^*)^*$ (see
\cite{Schn},Lem.9.9), hence, to check that the sequence
$(f_{n+1}-f_n)_n$ converges toward $0$ in $B$, we can do it inside
$(B^*)^*$. It is thus enough to prove that for any tempered
distribution $\mu$ on $\Z$ of order ${\rm val}(\alpha)$, we have:  
$$\sup_{\mu}\frac{\left|\int_{\Z}(f_{n+1}(z)-f_n(z))\mu(z)\right|}{|\!|\mu|\!|_{{\rm
      val}(\alpha)}}\longrightarrow 0\ {\rm when}\ n\rightarrow
+\infty.$$  
But: 
\begin{eqnarray*} 
\int_{\Z}(f_{n+1}(z)-f_n(z))\mu(z)&=&(\alpha p\beta^{-1})^n
\Big(\int_{p^n\Z}z^{k-2-j}\mu(z) 
\ \ \ \ \ -\int_{p^{n+1}\Z}z^{k-2-j}\mu(z)\Big)\\  
&\in &{|\!|\mu|\!|'_{{\rm val}(\alpha)}}^{-1}p^{n(2{\rm
    val}(\alpha)-k+2)}p^{n(k-2-j-{\rm val}(\alpha))}\O  
\end{eqnarray*} 
using that ${\rm val}(\alpha)+{\rm val}(\beta)=k-1$ and that $\mu$ is
of order ${\rm val}(\alpha)$ (see \S\ref{Yvette}). This implies
$\int_{\Z}(f_{n+1}(z)-f_n(z))\mu(z)\in C_{\mu}{|\!|\mu|\!|_{{\rm
    val}(\alpha)}}^{-1}p^{n({\rm val}(\alpha)-j)}\O$ hence the result since
$j<{\rm val}(\alpha)$.  
\end{proof} 
 
If ${\rm val}(\alpha)=k-1$, note that the functions $f$ in
$L(\alpha)$ give rise to pairs $(f_1,f_2)$ with zero $(k-1)$-derivative, hence $L(\alpha)\subseteq N(\alpha)\subsetneq B(\alpha)$ in that case. If ${\rm
  val}(\alpha)<k-1$, the authors do not know how to
prove directly that $L(\alpha)$ is distinct from $B(\alpha)$.  

\subsection{Definition of $\Pi(V)$ and first results}\label{banachsuite} 
 
1) Assume $V$ absolutely irreducible, or equivalently $0<{\rm
   val}(\alpha)<k-1$. Then we define:  
$$\Pi(V):=B(\alpha)/L(\alpha)$$ 
with the induced action of $\G$. We will prove in the next lectures
that $\Pi(V)$ is always non zero, but this will use
$(\varphi,\Gamma)$-modules.\\  
 
2) Assume $V$ is reducible and non-split. This implies ${\rm
   val}(\alpha)=k-1$ and ${\rm val}(\beta)=0$. We define again:  
$$\Pi(V):=B(\alpha)/L(\alpha)\ne 0$$  
with the induced action of $\G$. By Lem.\ref{exact}, we have an exact sequence: 
\begin{eqnarray}\label{exact2} 
0\rightarrow \frac{N(\alpha)}{L(\alpha)}\rightarrow \Pi(V)\rightarrow
\!\big({\rm Ind}_{\B}^{\G}\!x^{k-1}{\rm unr}(\alpha^{\!-1})\otimes
x^{\!-1}{\rm unr}(p\beta^{-1})\big)^{\!{\mathcal
    C}^0}\!\!\!\!\rightarrow 0.  
\end{eqnarray} 
We will identify below the Banach $N(\alpha)/L(\alpha)$.\\ 
 
3) Assume $V$ is reducible and split. This implies ${\rm
   val}(\alpha)=k-1$ and ${\rm val}(\beta)=0$. Then we define:  
\begin{eqnarray*} 
\Pi(V)&:=&B(\beta)\oplus \big({\rm Ind}_{\B}^{\G}x^{k-1}{\rm
  unr}(\alpha^{-1})\otimes x^{-1}{\rm
  unr}(p\beta^{-1})\big)^{{\mathcal C}^0}\\  
&=&\big({\rm Ind}_{\B}^{\G}{\rm unr}(\beta^{-1})\otimes x^{k-2}{\rm
  unr}(p\alpha^{-1})\big)^{{\mathcal C}^0} \oplus \\  
&&\ \ \ \ \ \big({\rm Ind}_{\B}^{\G}x^{k-1}{\rm
  unr}(\alpha^{-1})\otimes x^{-1}{\rm
  unr}(p\beta^{-1})\big)^{{\mathcal C}^0}.  
\end{eqnarray*} 
In that case, $\Pi(V)$ can also be written: 
\begin{eqnarray*} 
\big({\rm Ind}_{\B}^{\G}{\rm unr}(\beta^{-1})\otimes
\varepsilon^{k-2}{\rm unr}(p^{k-1}\alpha^{-1})\big)^{{\mathcal C}^0}
\oplus \\  
\ \ \ \ \ \big({\rm Ind}_{\B}^{\G}\varepsilon(\varepsilon^{k-2}{\rm unr}(p^{k-1}\alpha^{-1}))\otimes \varepsilon^{-1}{\rm unr}(\beta^{-1})\big)^{{\mathcal C}^0} 
\end{eqnarray*} 
where one can see some symmetry $({\rm
  Ind}_{\B}^{\G}\chi_1\otimes\chi_2)^{{\mathcal C}^0} \oplus ({\rm
  Ind}_{\B}^{\G}\varepsilon\chi_2\otimes\varepsilon^{-1}\chi_1)^{{\mathcal C}^0}$.  
 
We will now try to ``justify'' these definitions, or at least explain
where they are coming from. We also start studying $\Pi(V)$. 
 
Let $\pi(\alpha):={\rm Sym}^{k-2}E^2\otimes_E {\rm Ind}_{\B}^{\G}{\rm
  unr}(\alpha^{-1})\otimes {\rm unr}(p\beta^{-1})$ (the parabo\-lic
induction is here the usual smooth one): $\pi(\alpha)$ is an
irreducible locally algebraic representation of $\G$ in the sense of
D. Prasad (see appendix to \cite{ST2}). Likewise, we define
$\pi(\beta):={\rm Sym}^{k-2}E^2\otimes_E {\rm Ind}_{\B}^{\G}{\rm
  unr}(\beta^{-1})\otimes {\rm unr}(p\alpha^{-1})$. There are closed
$\G$-equivariant embeddings $\pi(\alpha)\hookrightarrow LA(\alpha)$
(resp. $\pi(\beta)\hookrightarrow LA(\beta)$) identifying the source
with the subspace of functions $f:\Q\rightarrow E$ such that
$f\!\mid_{\Z}$ and $(\alpha p\beta^{-1})^{{\rm
    val}(z)}z^{k-2}f(1/z)\!\mid_{\Z}$ (resp. $(\beta
p\alpha^{-1})^{{\rm val}(z)}z^{k-2}f(1/z)\!\mid_{\Z}$) are locally
polynomial functions on $\Z$ of degree at most $k-2$ (using
(\ref{groupeqp})). Moreover, the usual intertwining operator between
smooth generic principal series induces a $\G$-equivariant
isomorphism:  
\begin{eqnarray}\label{intertw} 
I:\pi(\alpha)\simeq \pi(\beta). 
\end{eqnarray} 
See (\ref{integrale2}) and (\ref{integrale}) later in the text for an
explicit description of $I$.  
 
\begin{rem} 
{\rm The smooth principal series ${\rm Ind}_{\B}^{\G}\!{\rm
    unr}(\alpha^{\!-1})\otimes {\rm unr}(p\beta^{\!-1}\!)$ is natural
  to introduce here since it is exactly the smooth representation of
  $\G$ that corresponds under the local Langlands correspondence
  (slightly twisted) to the unramified representation of ${\rm
    W}_{\Q}$ sending an arithmetic Frobenius to
  $\begin{pmatrix}\alpha&0 \\ 0&\beta\end{pmatrix}$= matrix of
    $\varphi^{-1}$ on $D(\alpha,\beta)$.}  
\end{rem} 
 
\begin{thm}\label{compl} 
1) Assume ${\rm val}(\alpha)<k-1$ (hence ${\rm val}(\beta)>0$) then
$B(\alpha)/L(\alpha)$ (resp. $B(\beta)/L(\beta)$) is isomorphic to the
completion of $\pi(\alpha)$ (resp. $\pi(\beta)$) with respect to any
$\O$-lattice of $\pi(\alpha)$ (resp. $\pi(\beta)$) which is finitely
generated over $\G$.\\  
2) Assume ${\rm val}(\alpha)=k-1$ (hence ${\rm val}(\beta)=0$) , then
$N(\alpha)/L(\alpha)$ (resp. $B(\beta)$) is isomorphic to the
completion of $\pi(\alpha)$ (resp. $\pi(\beta)$) with respect to any
$\O$-lattice of $\pi(\alpha)$ (resp. $\pi(\beta)$) which is finitely
generated over $\G$.  
\end{thm} 
 
Note that we don't know {\it a priori} in 1) that such a lattice
exists in $\pi(\alpha)$ or $\pi(\beta)$. In case not, the statement
means $B(\alpha)/L(\alpha)=0=B(\beta)/L(\beta)$. Before giving the
proof of Theorem \ref{compl}, we mention the following four
corollaries:  
  
\begin{cor}\label{1ercas} 
Assume ${\rm val}(\alpha)<k-1$, then $B(\alpha)/L(\alpha)$ and
$B(\beta)/L(\beta)$ are unitary $\G$-Banach spaces and we have a
commutative $\G$-equivariant diagram:  
$$\begin{matrix}B(\alpha)/L(\alpha)&\simeq &B(\beta)/L(\beta)\\
  \uparrow &&\uparrow \\ \pi(\alpha)& \buildrel I\over\simeq
  &\pi(\beta)\end{matrix}$$  
where $I$ is the intertwining operator of (\ref{intertw}). 
\end{cor} 
\begin{proof} 
The statement follows from Th.\ref{compl} together with the fact that
the isomorphism $I$ obviously sends an $\O$-lattice of $\pi(\alpha)$ of
finite type over $\G$ to an $\O$-lattice of $\pi(\beta)$ of finite
type over $\G$.  
\end{proof} 
 
We will see later on that $B(\alpha)/L(\alpha)$ and
$B(\beta)/L(\beta)$ are topologically irreducible and admissible.  
 
\begin{cor}\label{bon} 
Assume ${\rm val}(\alpha)=k-1$, then we have a commutative
$\G$-equivariant diagram: 
$$\begin{matrix}N(\alpha)/L(\alpha)&\simeq &B(\beta)\\ \uparrow
  &&\uparrow \\ \pi(\alpha)& \buildrel I\over\simeq
  &\pi(\beta)\end{matrix}$$  
where $I$ is the intertwining operator of (\ref{intertw}). 
Moreover, $N(\alpha)/L(\alpha)$ is a topologically irreducible
  admissible unitary $\G$-Banach and $B(\alpha)/L(\alpha)$ is an
  admissible unitary $\G$-Banach of topological length $2$ which is a
  non-trivial extension of the Banach $\big({\rm
  Ind}_{\B}^{\G}\!\!x^{k-1}\!{\rm unr}(\alpha^{\!-1}\!)\otimes
  x^{\!-1}\!{\rm unr}(p\beta^{\!-1}\!)\big)^{{\mathcal C}^0}$ by the
  Banach $N(\alpha)/L(\alpha)\simeq B(\beta)$.  
\end{cor} 
\begin{proof} 
The first part follows from Th.\ref{compl} as previously. The unitary
$\G$-Banach $B(\beta)$ is admissible and topologically
irreducible. The admissibility is true for any Banach $\big({\rm
  Ind}_{\B}^{\G}\chi_1\otimes \chi_2\big)^{{\mathcal C}^0}$ with
integral $\chi_i$ and follows from the fact that restriction of
functions to $\K$ induces an isomorphism with $\big({\rm
  Ind}_{\BK}^{\K}\chi_1\otimes \chi_2\big)^{{\mathcal C}^0}$\!\!\!,
the dual of which is $E\otimes (\O[[\K]]\otimes_{\O[[\BK]]}\O)$ (the
map $\BK\rightarrow \O$ being
$\begin{pmatrix}a&b\\0&d\end{pmatrix}\mapsto \chi_1(a)\chi_2(d)$)
  which is obviously of finite type over $E\otimes_{\O}\O[[\K]]$ (this
  proof is extracted from \cite{ST3}). The irreducibility can be
  proved as follows: if $B$ is a closed non zero $\G$-equivariant
  $E$-vector subspace in $B(\beta)$, then $B$ is also admissible and
  $B^{\rm an}\ne 0$ by the density of locally analytic vectors (see
  Schneider and Teitelbaum's course). But the locally analytic vectors
  in $B(\beta)$ are just $LA(\beta)$ hence $\pi(\beta)\subset B$ since
  $\pi(\beta)$ is the unique irreducible subobject of $LA(\beta)$. But
  $\pi(\beta)$ is dense by Prop.\ref{compl}, hence $B=B(\beta)$. This
  proves the statement on $N(\alpha)/L(\alpha)$. The proof that
  $B(\alpha)/L(\alpha)$ is admissible of topological length $2$ is
  analogous. The exact sequence is non-split because the analogous
  exact sequence with locally analytic vectors is non-split thanks to
  Cor.\ref{vectoranal} below (and the structure of locally analytic
  principal series for $\G$). Finally, let us prove that
  $B(\alpha)/L(\alpha)$ is unitary. We will use the fact that the
  surjection $B(\alpha)/L(\alpha)\twoheadrightarrow \big({\rm
    Ind}_{\B}^{\G}\!x^{k-1}{\rm unr}(\alpha^{\!-1})\otimes
  x^{\!-1}{\rm unr}(p\beta^{-1})\big)^{\!{\mathcal C}^0}$ admits a
  continuous section (not $\G$-equivariant of course). This follows
  from (\ref{f1f2}), (\ref{deriv}) and the fact that, for any positive
  integer $n$, the $n^{\rm th}$-derivation map ${\mathcal
    C}^n(\Z,E)\twoheadrightarrow {\mathcal C}^0(\Z,E)$ admits a
  continuous section (Cor.\ref{derivsurj}). Thus, as a Banach space,
  we can write $B(\alpha)/L(\alpha)$ as:  
\begin{eqnarray}\label{ondecoupe} 
B(\beta)\oplus \big({\rm Ind}_{\B}^{\G}x^{k-1}
{\rm unr}(\alpha^{-1})\otimes x^{-1}{\rm unr}(p\beta^{-1})\big)^{{\mathcal C}^0}. 
\end{eqnarray} 
Since the matrix $\begin{pmatrix}1&0\\0&p\end{pmatrix}$ acts
  continuously on $B(\alpha)/L(\alpha)$, replacing the unit ball
  $B(\beta)^0$ of $B(\beta)$ by $p^{-n}B(\beta)^0$ for a convenient
  integer $n$, we can assume, using (\ref{ondecoupe}) and the fact
  that the two Jordan-H\"older factors of $B(\alpha)/L(\alpha)$ are
  unitary, that the unit ball of $B(\alpha)/L(\alpha)$ is preserved by
  $\begin{pmatrix}1&0\\0&p\end{pmatrix}$ and $\K$, hence by $\G$ using
    the Cartan decomposition and the integrality of the central
    character.  
\end{proof} 
 
\begin{cor}\label{vectoranal2} 
Assume ${\rm val}(\alpha)<k-1$ and $B(\alpha)/L(\alpha)\ne 0$ 
(this fact will be proved later on), then we have a continuous $\G$-equivariant injection: 
$$LA(\alpha)\oplus_{\pi(\alpha)}LA(\beta)\hookrightarrow (B(\alpha)/L(\alpha))^{\rm an}$$ 
where $\pi(\alpha)$ embeds into $LA(\beta)$ via the intertwining
(\ref{intertw}).  
\end{cor} 
\begin{proof} 
We have a continuous equivariant injection $LA(\alpha)\hookrightarrow
(B(\alpha)/L(\alpha))^{\rm an}$ by (\ref{analytique}) (it is injective
since, otherwise, $\pi(\alpha)$ would necessarily map to zero and we
know this is not true because of Th.\ref{compl}). Likewise with
$\beta$. The corollary then follows from Cor.\ref{1ercas}.  
\end{proof} 
 
Because of Cor.\ref{vectoranal} below, I am tempted to conjecture that
the injection in Cor.\ref{vectoranal2} is actually a topological
isomorphism. The statement of Cor.\ref{vectoranal2} (without the
description of the completion of $\pi(\alpha)$ as
$B(\alpha)/L(\alpha)$) was also noted in \cite{Em},Prop.2.5.  
 
\begin{cor}\label{vectoranal} 
Assume ${\rm val}(\alpha)=k-1$, then we have a topological $\G$-equivariant isomorphism: 
$$LA(\alpha)\oplus_{\pi(\alpha)}LA(\beta)\buildrel\sim\over\rightarrow
(B(\alpha)/L(\alpha))^{\rm an}$$  
where $\pi(\alpha)$ embeds into $LA(\beta)$ via the intertwining (\ref{intertw}). 
\end{cor} 
\begin{proof} 
Since all the Banach in (\ref{exact2}) are admissible by
Cor.\ref{bon}, a theorem of Schneider and Teitelbaum tells us that we
have an exact sequence of locally analytic representations from the
sequence (\ref{exact2}):  
\begin{eqnarray}\label{exactanal} 
0\rightarrow LA(\beta)\rightarrow
\Big(\frac{B(\alpha)}{L(\alpha)}\Big)^{\rm an}\!\!\!\!\rightarrow
\big({\rm Ind}_{\B}^{\G}\!x^{k-1}\!{\rm unr}(\alpha^{\!\!-1})\otimes
x^{-1}\!{\rm unr}(p\beta^{\!-1})\big)^{\rm an}
\rightarrow 0 
\end{eqnarray} 
where we have used $(N(\alpha)/L(\alpha)^{\rm an}\simeq B(\beta)^{\rm
  an}\simeq LA(\beta)$. But we also have a continuous $\G$-equivariant
injection $LA(\alpha)\hookrightarrow (B(\alpha)/L(\alpha))^{\rm an}$
as in Cor.\ref{vectoranal2}, hence a continuous $\G$-equivariant map
$LA(\alpha)\oplus_{\pi(\alpha)}LA(\beta)\rightarrow
(B(\alpha)/L(\alpha))^{\rm an}$. To prove it is a topological
bijection, note that from the exact sequence of locally analytic
representations: 
\begin{eqnarray*} 
0\rightarrow \pi(\alpha)\rightarrow LA(\alpha)\rightarrow 
\big({\rm Ind}_{\B}^{\G}x^{k-1}{\rm unr}(\alpha^{-1})\otimes
x^{-1}{\rm unr}(p\beta^{-1})\big)^{\rm an}\rightarrow 0,  
\end{eqnarray*} 
we deduce an exact sequence similar to (\ref{exactanal}) with
$(B(\alpha)/L(\alpha))^{\rm an}$ replaced by
$LA(\alpha)\oplus_{\pi(\alpha)}LA(\beta))$, together with an obvious
commutative diagram between the two exact sequences.  
\end{proof} 
 
Proof of Th.\ref{compl}: Note first that the completion of
$\pi(\alpha)$ doesn't depend on which lattice of finite type is chosen
(if any) since all these lattices are commensurable and thus give rise
to equivalent invariant norms. Moreover, this completion is also the
same as the completion with respect to any $\O$-lattice finitely
generated over $\B$. Indeed, any $\O$-lattice of finite type over $\G$
is of finite type over $\B$ since $\G=\B\K$ and the $\K$-span of any
vector in $\pi(\alpha)$ is finite dimensional. Conversely, one can
check using the compacity of $\K$ that any $\O$-lattice of finite type
over $\B$ is contained in an $\O$-lattice of finite type over $\G$ (and is
thus commensurable to it). If $\pi(\alpha)$ admits an $\O$-lattice
preserved by $\B$, then one can check that the following
$\O$-submodule:  
$$\sum_{j=0}^{k-2}\O[\B]z^j{\mathbf 1}_{\Z}+\sum_{j=0}^{k-2}
\O[\B](\alpha p\beta^{-1})^{{\rm val}(z)}z^j{\mathbf 1}_{\Q-p\Z}$$ 
(viewing $\pi(\alpha)$ as embedded into $LA(\alpha)$) is necessarily
an $\O$-lattice of $\pi(\alpha)$ which is finitely generated over
$\B$. If $\pi(\alpha)$ doesn't admit any $\O$-lattice preserved by
$\B$, then one can check that the above $\O$-submodule is the
$E$-vector space $\pi(\alpha)$. The dual of the sought after
completion of $\pi(\alpha)$ is thus isomorphic to:  
\begin{eqnarray}\label{thedual} 
\{\mu\in \pi(\alpha)^*\ \mid\ \forall g\in\B, \forall j\in
\{0,\cdots,k-2\}, |\mu(g(z^j{\mathbf 1}_{\Z}))|\leq 1\\   
\nonumber {\rm and}\ |\mu(g((\alpha p\beta^{-1})^{{\rm
    val}(z)}z^j{\mathbf 1}_{\Q-p\Z}))|\leq 1\}\otimes_{\O}E.  
\end{eqnarray} 
Granting the fact that the central character is integral, we can even
take $g\in \begin{pmatrix}1&\Q \\0&\Q^{\times}\end{pmatrix}$ in
(\ref{thedual}). For $f\in \pi(\alpha)$, seeing $f$ as a function on
$\Q$ via (\ref{groupeqp}), we write $\int_{\Q}f(z)\mu(z)$ instead of
$\mu(f)$ in the sequel. A short computation then gives that the
conditions in (\ref{thedual}) are equivalent to the existence for each
$\mu$ of a constant $C_{\mu}\in E$ such that $\forall a\in \Q$,
$\forall j\in \{0,\cdots,k-2\}$ and $\forall n\in {\mathbb Z}$,
$n>{\rm val}(a)$ if $a\ne 0$:  
\begin{eqnarray}\label{chaud1} 
\int_{a+p^n\Z}(z-a)^j\mu(z)&\in & C_{\mu}p^{n(j-{\rm
    val}(\alpha))}\O\\  
\label{chaud2}\int_{\Q-(a+p^n\Z)}\big(\frac{\alpha
  p}{\beta}\big)^{{\rm val}(z-a)}(z-a)^{k-2-j}\mu(z)&\in &
C_{\mu}p^{n({\rm val}(\alpha)-j)}\O  
\end{eqnarray} 
where $\int_U f(z)\mu(z):=\int_{\Q}f(z){\mathbf 1}_U(z)\mu(z)$ for any
open $U\subset \Q$. Let us now first assume ${\rm
  val}(\alpha)<k-1$. From (\ref{chaud1}), we deduce (modifying
$C_{\mu}$ if necessary):  
\begin{eqnarray}\label{chaud3} 
\int_{(a+p^{n-1}\Z)-(a+p^n\Z)}\!\!\!\big(\frac{\alpha
  p}{\beta}\big)^{{\rm val}(z-a)}(z-a)^{k-2-j}\mu(z)\in
C_{\mu}p^{n({\rm val}(\alpha)-j)}\O.  
\end{eqnarray} 
Writing $\Q-(a+p^n\Z)=\Q-(a+p^{n-1}\Z)\amalg (a+p^{n-1}\Z)-(a+p^n\Z)$,
using (\ref{chaud3}) for the right hand side and decomposing again
$\Q-(a+p^{n-1}\Z)=\Q-(a+p^{n-2}\Z)\amalg (a+p^{n-2}\Z)-(a+p^{n-1}\Z)$,
we see by induction that $\int_{\Q-(a+p^n\Z)}\big(\frac{\alpha
  p}{\beta}\big)^{{\rm val}(z-a)}(z-a)^{k-2-j}\mu(z)$ for $j>{\rm
  val}(\alpha)$ is uniquely determined and that (\ref{chaud2}) for
$j>{\rm val}(\alpha)$ is a consequence of (\ref{chaud1}). Furthermore,
we see by the same argument that (\ref{chaud2}) for $j={\rm
  val}(\alpha)$ (if ${\rm val}(\alpha)$ is a positive integer) and all
$a$ follows from (\ref{chaud2}) for $j={\rm val}(\alpha)$, $a=0$ and
from (\ref{chaud1}) (develop $(z-a)^{k-2-j}$ when $n\ll
0$). Moreover, we know by Lem.\ref{etendre} that (\ref{chaud1}) for
$j>{\rm val}(\alpha)$ is a consequence of (\ref{chaud1}) for $j\leq
{\rm val}(\alpha)$. Hence, (\ref{chaud1}) and (\ref{chaud2}) are
equivalent to: (\ref{chaud1}) for $j\leq {\rm val}(\alpha)$ and $a\in
\Q$ + (\ref{chaud2}) for $j<{\rm val}(\alpha)$ and $a\in \Q$ +
(\ref{chaud2}) for $j={\rm val}(\alpha)$ and $a=0$. Any $\mu\in
\pi(\alpha)^*$ is equivalent to the data of a pair $(\mu_1,\mu_2)$
where, if $f\in \pi(\alpha)$,
$\int_{\Q}f(z)\mu(z)=\int_{\Z}f(pz)\mu_1(z)+\int_{\Z}(\alpha
p\beta^{-1})^{{\rm val}(z)}f(1/z)\mu_2(z)$. Assuming $j\leq {\rm
  val}(\alpha)$, we then see by an easy computation (left to the
reader) that [(\ref{chaud1}) for $j={\rm val}(\alpha)$ +
  (\ref{chaud1}) for $j<{\rm val}(\alpha)$ and $a\in \Q^{\times}$ +
  (\ref{chaud1}) for $j<{\rm val}(\alpha)$, $a=0$ and $n\geq 0$ +
  (\ref{chaud2}) for $j={\rm val}(\alpha)$ and $a=0$ + (\ref{chaud2})
  for $j<{\rm val}(\alpha)$, $a=0$ and $n\leq 0$] is equivalent to
$\mu_1$ and $\mu_2$ satisfying condition (\ref{AMcondition}) of
Lem.\ref{etendre} with $r={\rm val}(\alpha)$ and $d=k-2$. Hence
$\mu_1$ and $\mu_2$ are in fact tempered distributions on $\Z$ of
order ${\rm val}(\alpha)$ by Cor.\ref{travail}. Thus we can integrate
{\it any} function in $B(\alpha)$ against $\mu$, in particular the
functions of $L(\alpha)$. By making $n$ tends to $-\infty$, we see
that (\ref{chaud1}) for $j<{\rm val}(\alpha)$ and $a=0$ implies $\mu$
cancels the functions $z^j\in L(\alpha)$. By making $n$ tends to
$+\infty$, we see that (\ref{chaud2}) for $j<{\rm val}(\alpha)$
implies $\mu$ cancels the functions $\big(\frac{\alpha
  p}{\beta}\big)^{{\rm val}(z-a)}(z-a)^{k-2-j}\in L(\alpha)$. Now, the
reader can check that $\mu_1$ and $\mu_2$ being tempered of order
${\rm val}(\alpha)$ together with these cancellations is in fact {\it
  equivalent} to (\ref{chaud1}) for $j\leq {\rm val}(\alpha)$ and
$a\in \Q$ + (\ref{chaud2}) for $j<{\rm val}(\alpha)$ and $a\in \Q$ +
(\ref{chaud2}) for $j={\rm val}(\alpha)$ and $a=0$. This gives a
$\G$-equivariant isomorphism between the dual of the completion of
$\pi(\alpha)$ in (\ref{thedual}) and $(B(\alpha)/L(\alpha))^*$. A
closer examination shows that this is a topological isomorphism of
Banach spaces and thus $(B(\alpha)/L(\alpha))^*$ is a unitary
$\G$-Banach space. Using the closed embedding
$B(\alpha)/L(\alpha)\hookrightarrow ((B(\alpha)/L(\alpha))^*)^*$, we
see that the $\G$-Banach $B(\alpha)/L(\alpha)$ is unitary since it is
closed in a unitary $\G$-Banach. This means that the canonical map
$\pi(\alpha)\rightarrow B(\alpha)/L(\alpha)$ induces a continuous
$\G$-equivariant morphism from the completion of $\pi(\alpha)$ (with
respect to any finite type lattice) to $B(\alpha)/L(\alpha)$. This
morphism is, as we have seen, an isomorphism on the duals for the
strong topologies, hence also for the weak topologies (using the
theory of \cite{ST3}). By biduality, we finally get 1) for
$B(\alpha)/L(\alpha)$ and also for $B(\beta)/L(\beta)$ by
interchanging $\alpha$ and $\beta$. The proof of 2) is similar (and
actually simpler since we only have to deal with $j<{\rm
  val}(\alpha)$), the only difference being that here, the conditions
(\ref{chaud1}) and (\ref{chaud2}) only allow you to integrate
functions in $N(\alpha)$ against $\mu$ (see the end of
\S\ref{Yvette}).

\section{Representations of $\G$ and $(\varphi,\Gamma)$-mod\-ules I (C.B.)} 
We now focus on the study of $\Pi(V)$ when $V$ is irreducible. Quite
surprisingly, this will crucially use the $(\varphi,\Gamma)$-module
$D(V)$ associated to $V$. We define here a topological space
$\big(\limproj_{\psi}D(V)\big)^{\rm b}$ and give a first explicit
description of it using the Wach module $N(T)$ associated to a Galois
lattice in $V$. We then give preliminary results on distributions on
$\Q$ naturally arising from $\big(\limproj_{\psi}D(V)\big)^{\rm b}$
via Amice transform.  
 
\subsection{Back to $(\varphi,\Gamma)$-modules}\label{back1} 
 
Let $V$ be as in case 1) of \S\ref{banachdebut}, $\Pi(V)$ as in \S\ref{banachsuite} and  
$D(V)$ be the $(\varphi,\Gamma)$-module associated to $V$ in \S\ref{d(v)}. Recall that we 
have defined in \S\ref{psy} a map $\psi:D(V)\twoheadrightarrow D(V)$ which is a 
left-inverse of $\varphi$ (and hence a surjection).   
Define the following $E$-vector space: 
\begin{multline*} 
\big(\limproj_{\psi}D(V)\big)^{\rm b} 
:=\{(v_n)_{n\in{\mathbb Z}_{\geq 0}}\mid v_n\in D(V),\ \psi(v_{n+1})=v_n, 
(v_n)_n{\rm\ a\ bounded\ sequence\ in\ } D(V)\}, 
\end{multline*} 
(recall that \emph{bounded} means bounded for the weak topology which 
was defined in \S\ref{defweaktopo}). 
We equip the above space with the following structures:\\ 
1) a left $\O[[X]]$-module structure, by  
$s((v_n)_n):=(\varphi^n(s)v_n)_n$, if $s \in \O[[X]]$\\ 
2) a bijection $\psi$ given by $\psi((v_n)_n):=(\psi(v_n))_n$\\ 
3) an action of $\Gamma$ given $\gamma((v_n)_n):=(\gamma v_n)_n$.\\ 
 
The aim of this lecture and the next one is to  
prove a ``canonical'' topological isomorphism: 
\begin{thm}\label{main} 
Assume that $V$ is absolutely irreducible. There is an  
isomorphism of topological $E$-vector spaces (for  
the weak topologies on both sides): 
$$\big(\limproj_{\psi}D(V)\big)^{\rm b}\simeq \Pi(V)^*$$ 
such that the action of $\begin{pmatrix}1&0\\0&p\end{pmatrix}$ on  
$\Pi(V)^*$ corresponds to $(v_n)_n\mapsto (\psi(v_n))_n$, the action  
of $\begin{pmatrix}1&0\\0&\Z^{\times}\end{pmatrix}$ to that of  
$\Gamma\simeq \Z^{\times}$ and the action of  
$\begin{pmatrix}1&\Z\\0&1\end{pmatrix}$  
to the multiplication by $(1+X)^{\Z}$. 
\end{thm} 
 
Such an isomorphism was inspired by an analogous isomorphism due to
Colmez in the case $V$ is semi-stable non crystalline (\cite{Co3}). 
 
Before proving this theorem, we will need to establish a number of 
preliminary results. Let $T$ be a $\g$-stable lattice in $V$. 
First of all, because $V$ is assumed to be 
crystalline, we can write $D(T) = \aa_{\Q} \otimes_{\Z[[X]]} N(T)$ 
where $N(T)$ is the Wach module which we have constructed in 
\S\ref{LB5}. This is a free $\O[[X]]$-module of finite rank and thus it 
is equipped with the $(p,X)$-adic topology. The first result which we 
will need in order to prove theorem \ref{main} is the following one: 
 
\begin{prop}  
Let $T$ be a $\g$-stable lattice in a crystalline,  
irreducible representation $V$ with Hodge-Tate weights  
$0$ and $r>0$ and let $\limproj_{\psi}N(T)$ denote  
the set of $\psi$-compatible sequences of elements  
of $N(T)$. Then the natural map:  
$$\big(\limproj_{\psi}N(T)\big) \otimes_{\O}E\rightarrow  
\big(\limproj_{\psi}D(V)\big)^{\rm b}$$ 
is a topological isomorphism. 
\end{prop} 
\begin{proof} 
The main content of this proposition is the  
explicit description of the 
topology of $D(V)$ in terms of $N(T)$ which we have given in 
\S\ref{LB5} and its interaction with $\psi$.  
Recall that for each $k \geq 0$, we define a semi-valuation 
$\nu_k$ on $D(T)$ as follows: if $x \in D(T)$ then $\nu_k(x)$ is 
the largest integer $j \in \mathbb{Z} \cup \{ + \infty\}$  
such that $x \in X^j N(T) + p^k D(T)$. The weak topology on $D(T)$ is 
then the topology defined by the set $\{\nu_k\}_{k \geq 0}$ of all 
those semi-valuations. 
The weak topology on $D(V)$ is the inductive 
limit topology on $D(V) = \cup_{\ell \geq 0} D(p^{-\ell}T)$.  
Concretely, if we  
have a sequence $(v_n)_n$ of elements of $D(V)$, and 
that sequence is bounded for the weak topology, then there is a 
$\g$-stable lattice $T$ of $V$ such that $v_n \in D(T)$ for every $n \geq 
0$ and furthermore, for every $k \geq 0$, there exists $f(k) \in 
\mathbb{Z}$ such that $v_n \in X^{-f(k)} N(T) + p^k D(T)$.  
Recall that we have proved (cf. Th.\ref{psionwach}) 
that if $V$ is crystalline with 
positive Hodge-Tate weights, and $\ell \geq 1$, then: 
$$\psi(X^{-\ell} N(T)) \subset p^{\ell-1} X^{-\ell} N(T) + X^{-(\ell-1)} N(T)$$ 
by iterating this $m$ times we get:  
$$\psi^m(X^{-\ell} N(T)) \subset p^{m(\ell-1)} X^{-\ell} N(T) + X^{-(\ell-1)} N(T).$$ 
Choose $k \geq 0$ and $\ell \geq 2$ such that  
$v_n \in X^{-\ell} N(T) + p^k D(T)$. Since $v_n = \psi^m  
(v_{n+m})$ for every $m \geq 1$, we see that $v_n \in p^{m(\ell-1)} X^{-\ell} N(T) +  
X^{-(\ell-1)} N(T) + p^k D(T)$ for all $m$ and by taking $m$ large 
enough we get that $v_n \in X^{-(\ell-1)} N(T) + p^k D(T)$ so that if  
$\ell \geq 2$ then we may replace it by $\ell-1$. 
Hence, for every $k \geq 0$, we have $v_n \in X^{-1} N(T) + p^k 
D(T)$. This being true for every $k$, we get $v_n \in X^{-1}N(T)$.  
We will now prove that $v_n \in N(T)$ if $V$ is irreducible (it is 
in fact enough to assume that the $\varphi$-module $\dcris(V)$  
has no part of slope $0$). Indeed, both  
$X^{-1}N(V)$ and $N(V)$ are stable by $\psi$ 
and we have seen in proposition \ref{psidcris} 
that the map $\psi : \pi^{-1} 
N(V) / N(V) \rightarrow  \pi^{-1} N(V) / N(V)$ coincides, via the 
identification $\pi^{-1} N(V) / N(V) = \dcris(V)$ given by 
$\pi^{-1}x \mapsto \overline{x}$, with the map $\varphi^{-1} : 
\dcris(V) \rightarrow \dcris(V)$. If $\dcris(V)$ has no part of 
slope $0$, then $\varphi^{-m}(x) \rightarrow 0$ as $m \rightarrow \infty$ 
and so for every $\ell$ there exists $m$ such that 
$\psi^m(X^{-1}N(T)) \subset p^{\ell}X^{-1}N(T) + N(T)$. Consequently, 
we have  $v_n \in N(T)$ for all $n \geq 0$. To finish the proof  
of the proposition, note that $\limproj_{\psi}N(T)$ is compact and that the map:  
$$\limproj_{\psi}N(T)\rightarrow \big(\limproj_{\psi}D(T)\big)^{\rm b}$$ 
is a continuous isomorphism, so that it is a topological isomorphism.  
The proposition follows by inverting $p$. 
\end{proof} 
 
Next, recall that in Prop.\ref{inccris} 
we have seen that if $V$ is a 
crystalline irreducible representation whose Hodge-Tate weights are 
$\geq 0$, then we have an injective map $N(T) \hookrightarrow 
{\mathcal R}_{\Q}^+ \otimes_{\Q} \dcris(V)$. Now let 
$V=V_{cris}(D(\alpha,\beta))$, $V$ irreducible. Given a sequence $(v_n)_n \in 
\limproj_{\psi}N(T)$ we can write $v_n = w_{\alpha,n} \otimes e_\alpha 
+ w_{\beta,n} \otimes e_\beta$. Recall that we write ${\mathcal 
R}_E^+$ for $E \otimes_{\Q} {\mathcal R}_{\Q}^+$. 
 
\begin{prop}\label{conditions} 
Given $(v_n)_n \in (\limproj_{\psi}N(T)) \otimes_{\O} E$,  
the sequences $(w_{\alpha,n})_n$ and $(w_{\beta,n})_n$  
of elements of ${\mathcal R}_E^+$ defined above satisfy the following
properties:\\  
1) $\forall\ n\geq 0$, $w_{\alpha,n}$ (resp. $w_{\beta,n}$) is of  
order ${\rm val}(\alpha)$ (resp. ${\rm val}(\beta)$) and  
$\|w_{\alpha,n}\|_{{\rm val}(\alpha)}$ (resp. $\|w_{\beta,n}\|_{{\rm val}(\beta)}$)  
is bounded independently of $n$\\ 
2) $\forall\ n\geq 0$ and $\forall\ m\geq 1$, we have 
$$\varphi^{-m}(w_{\alpha,n}\otimes e_{\alpha}+w_{\beta,n}\otimes 
e_{\beta})\in {\rm Fil}^0((\Q(\mu_{p^m})\otimes_{\Q}E)[[t]]\otimes_ED(\alpha,\beta))$$ 
3) $\forall\ n\geq 1$, $\psi(w_{\alpha,n})=\alpha^{-1}w_{\alpha,n-1}$ and  
$\psi(w_{\beta,n})=\beta^{-1}w_{\beta,n-1}$. 
\end{prop} 
\begin{proof} 
We can always change $T$ in order to have $(v_n)_n \in 
\limproj_{\psi}N(T)$ which we now assume. 
Let $c_\alpha : N(T) \rightarrow {\mathcal R}_E^+$ be the map which to 
$x=x_{\alpha} \otimes e_{\alpha} + x_{\beta} \otimes e_{\beta}$ 
assigns $x_\alpha$. 

To prove that $c_\alpha(N(T))$ is of order  
$\leq {\rm val}(\alpha)$, we use the 
following characterization of elements of order $r$ in 
${\mathcal R}_E^+$. Recall that ${\mathcal R}_E^+$ is included in a 
larger ring $\widetilde{\mathcal{R}}_E^+:=\cap_{n\geq
  0}\varphi^n(\bcris^+)$ 
on which $\varphi$ is a bijection 
and which has a topology defined by the family of semi-norms 
$\|\cdot\|_D$ where $D$ runs over all closed disks of radii $<1$. 
Fix one such closed disk $D=D(0,\rho)$ of radius $\rho < 1$.  
We say that $f \in \widetilde{\mathcal{R}}_E^+$ is of order $r$ if  
the sequence $\{ p^{-mr} \| \varphi^{-m} f \|_D \}$ is bounded 
independently of $m$. We can then define $\|f\|_{D,r}:=\sup_m(p^{-mr}  
\| \varphi^{-m} f \|_D)$. The norms $\|\cdot\|_{D,r}$ (which depend 
on the choice of $D$) are equivalent to the norms $\|\cdot\|_r$  
defined in \S\ref{c^r}. Indeed, if $D_r$ has radius $r$, then
$\| \varphi^{-m} f \|_{D_r} = \| f \|_{D_{r^{1/p^m}}}$  and one
recovers the definition of order as ``order of growth''.

Now if $x \in N(T)$ and  
$x=x_{\alpha} \otimes e_{\alpha} + x_{\beta} 
\otimes e_{\beta}$ then: 
$$\varphi^{-m}(x) = \varphi^{-m}(x_{\alpha}) 
\alpha^{-m} \otimes e_{\alpha} + \varphi^{-m} (x_{\beta}) \beta^{-m}  
\otimes e_{\beta}$$ 
and the set of $\varphi^{-m}(x)$ is bounded for any 
$\|\cdot\|_D$ norm  
(because $x \in \aplus[1/X] \otimes_{\Z} T$) 
so that the sets $\{\varphi^{-m}(x_{\alpha}) 
\alpha^{-m}\}$ and $\{\varphi^{-m}(x_{\beta}) \beta^{-m}\}$ are both 
bounded. It thus follows that $x_\alpha$ is of order ${\rm 
val}(\alpha)$ and $x_\beta$ of order ${\rm val}(\beta)$. We have  
proved that the image of $c_{\alpha}$ is made up of elements of  
order ${\rm val}(\alpha)$ and it is 
bounded in the $\|\cdot\|_{{\rm val}(\alpha)}$ norm because $N(T)$ is 
an $\O[[X]]$-module of finite type. The same holds with $\beta$ in 
place of $\alpha$ of course. For the second point, we use the fact  
that $v_n \in N(T)$ so that for all $m \geq 1$ 
we have $\varphi^{-m}(v_n) \in B_{dR}^+ \otimes_{\Q} 
V$. On the other hand by the results of \S\ref{LB41}  
(see Prop.\ref{imgiotan}) we have  
$\varphi^{-m}(v_n) \in (\Q(\mu_{p^m})\otimes_{\Q}E)((t)) 
\otimes_ED(\alpha,\beta)$ and therefore: \begin{multline*}  
\varphi^{-m}(v_n) \in  
(B_{dR}^+ \otimes_{\Q} 
V) \cap \big((\Q(\mu_{p^m})\otimes_{\Q}E)((t))  
\otimes_ED(\alpha,\beta)\big) 
\\ = {\rm Fil}^0 
\big((\Q(\mu_{p^m})\otimes_{\Q}E)((t)) 
\otimes_ED(\alpha,\beta)\big).  
\end{multline*}
Since the weights of the filtration on $D(\alpha,\beta)$ are negative, we can replace $(\Q(\mu_{p^m})\otimes_{\Q}E)((t))$ by $(\Q(\mu_{p^m})\otimes_{\Q}E)[[t]]$ in the ${\rm Fil}^0$. Finally, the third point is obvious since $v_{n-1}=\psi(v_n)$  
and $\psi$ acts as $\varphi^{-1}$ on $D_{cris}(V)$. 
\end{proof} 
 
\begin{prop}\label{conditions2} 
Conversely, if we are given two sequences $(w_{\alpha,n})_n$ and $(w_{\beta,n})_n$  
of elements of ${\mathcal R}_E^+$ satisfying properties 1)-3) of  
Prop.\ref{conditions} and if we set $v_n := w_{\alpha,n} \otimes  
e_\alpha + w_{\beta,n} \otimes e_\beta$, then $(v_n)_n \in  
(\limproj_{\psi}N(T)) \otimes_{\O} E$. 
\end{prop} 
\begin{proof} 
It is enough to show that there exists a lattice $T$ such that $v_n 
\in N(T)$ because condition 3) of Prop.\ref{conditions} ensures  
that the sequence $(v_n)_n$ is $\psi$-compatible. Choose any lattice $U$. By the  
results of \S\ref{LB4} and \S\ref{LB5} (more specifically proposition
\ref{groumpf}), we know that 
$v_n \in {\mathcal R}_E^+[1/t] \otimes_{\O[[X]]} N(U)$. 

Condition 2) of proposition \ref{conditions} then implies that 
$v_n \in {\mathcal R}_E^+[1/X] \otimes_{\O[[X]]} N(U)$, because the
filtration condition for $m \geq 1$ is equivalent to an
order-of-vanishing condition at $\epsilon^{(m)}-1$, and implies
that $v_n \in {\mathcal R}_E^+[\varphi^{m-1}(q) / t] \otimes_{\O[[X]]}
N(U)$. By using this for each $m \geq 1$, we get that 
$v_n \in {\mathcal R}_E^+[1/X] \otimes_{\O[[X]]} N(U)$.

The $\psi$-compatibility allows us to get rid of the denominators in $X$.  
Finally, condition 1) tells us that there exists $\ell$ independent of $n$  
such that $v_n \in p^{-\ell} N(U)$ so we can take $T=p^{-\ell}U$.  
\end{proof} 
 
It is important to notice that the map $\psi: N(T) \rightarrow N(T)$  
is {\it not} surjective. Recall that we have in particular proved  
in Prop.\ref{bh} the following result (we actually proved more, and we proved it for 
$\Q$-linear representations, but as usual the $E$-linear result 
follows immediately): 
 
\begin{prop}\label{berger} 
If $V$ is crystalline absolutely irreducible, 
then there exists a \emph{unique} 
non-zero $E$-vector subspace $D^0(V)$ of $D(V)$  
possessing an $\O$-lattice $D^0(T)$  
which is a compact $\O[[X]]$-submodule  
of $D(V)$ preserved by $\psi$ and $\Gamma$ with $\psi$ surjective. 
\end{prop} 
 
In fact, as we saw, the lattice $D^0(T)$ is 
of finite type over $\O[[X]]$. The following corollary will 
be used in \S\ref{results} to prove the  
irreducibility of $\Pi(V)$ (for $V$ irreducible): 
 
\begin{cor}\label{versirr} 
Let ${\mathcal M}\subseteq \limproj_{\psi}N(T)$ be a closed non zero  
$\O[[X]]$-submodule which is stable by $\psi$, $\psi^{-1}$ and $\Gamma$, then: 
$${\mathcal M}\otimes_{\O}E=(\limproj_{\psi}N(T))\otimes_{\O}E 
=\big(\limproj_{\psi}D(V)\big)^{\rm b}.$$ 
\end{cor} 
\begin{proof}
Note first that since $\limproj_{\psi}N(T)$ is compact, so is
${\mathcal M}$ (being closed in a compact set). For $m\in {\mathbb
  Z}_{\geq 0}$, denote by   
${\rm pr}_m:\limproj_{\psi}N(T)\rightarrow N(T)$,  
$(v_n)_n\mapsto v_m$ and define $M_m:={\rm pr}_m({\mathcal M})$. 
As ${\mathcal M}\ne 0$, there exists $m$ such that $M_m\ne 0$. 
Using that ${\mathcal M}$ is stable by $\psi$ and $\psi^{-1}$, 
it is straightforward to check that $M_m$ doesn't depend on $m$  
and we denote it by $M$. Via $M=M_0$, it is an $\O[[X]]$-submodule 
of $N(T)$ (necessarily of finite type since $N(T)$ is) 
stable by $\Gamma$ and $\psi$ and  
such that $\psi$ is surjective. Moreover, the canonical map  
${\mathcal M}\rightarrow \limproj_{\psi}M$ is 
easily checked to be an isomorphism  
since ${\mathcal M}$ is both dense in $\limproj_{\psi}M$ and  
compact (the density follows from $M={\rm pr}_m({\mathcal M})$  
$\forall\ m$). Doing the same with $\limproj_{\psi}N(T)$, we  
find another (finite type) $\O[[X]]$-submodule $M^0$ of $N(T)$  
containing $M$ stable by $\Gamma$ and $\psi$ with $\psi$  
surjective and such that $\limproj_{\psi}N(T)=\limproj_{\psi}M^0$.  
By Prop.\ref{berger}, we have thus $M\otimes_{\O}E=M^0\otimes_{\O}E$  
and there exists $a\in E^{\times}$ such that 
$M\subset M^0\subset aM$. Hence, we have  
$(\limproj_{\psi}M)\otimes_{\O}E\simeq 
(\limproj_{\psi}M^0)\otimes_{\O}E$  
which finishes the proof. 
\end{proof} 
 
\subsection{Back to distributions}\label{distrback} 
 
Recall from \S\ref{Yvette} that $\mu\mapsto
\sum_{n=0}^{+\infty}\int_{\Z}{z\choose n}\mu(z)X^n$ gives a bijection
between ${\mathcal R}_E^+$ and ``$E$-valued'' locally analytic
distributions on $\Z$ (some people write ``$\mu\mapsto
\int_{\Z}(1+X)^z\mu(z)$'') which induces a bijection between elements
of ${\mathcal R}_E^+$ of order $r$ and tempered distributions of order
$r$. We will need the following easy lemma (we leave its proof to the
reader):  
 
\begin{lem}\label{Psidis} 
If $w\in {\mathcal R}_E^+$ corresponds to $\mu$, then $\psi(w)\in
{\mathcal R}_E^+$ corresponds to the unique locally analytic
distribution $\psi(\mu)$ on $\Z$ such that:  
$$\int_{\Z}f(z)\psi(\mu)(z):=\int_{p\Z}f(z/p)\mu(z)$$ 
for any locally analytic function $f$ on $\Z$. 
\end{lem} 
 
We have an analogous lemma for tempered distributions of order $r$ and
functions of class ${\mathcal C}^r$ (as one easily checks that $\psi$
preserves elements in ${\mathcal R}_E^+$ of order $r$). In particular,
if $(\mu_n)_{n\in{\mathbb Z}_{\geq 0}}$ is a sequence of tempered
distributions on $\Z$ of order $r$ such that $\psi(\mu_{n+1})=\mu_n$,
we can define a linear form $\mu$ on the $E$-vector space ${\rm
  LPol}_{{\rm c},k-2}$ of locally polynomial functions
$f:\Q\rightarrow E$ with compact support of degree less than $k-2$ as
follows:  
$$\int_{\Q}f(z)\mu(z)=\int_{p^{-N}\Z}f(z)\mu(z):=\int_{\Z}f(z/p^N)\mu_N(z)$$ 
where $N\in {\mathbb Z}_{\geq 0}$ is such that supp$(f)\subset
p^{-N}\Z$. 
One readily checks using Lem.\ref{Psidis} that it doesn't depend on
such an $N$. 
 
Let $(v_n)_n=(w_{\alpha,n}\otimes e_{\alpha}+w_{\beta,n}\otimes
e_{\beta})_n\in (\limproj_{\psi}N(T))\otimes_{\O}E$. Recall this means
(Prop.\ref{conditions} and Prop.\ref{conditions2}):\\  
1) $\forall\ n\geq 0$, $w_{\alpha,n}$ (resp. $w_{\beta,n}$) is of
order ${\rm val}(\alpha)$ (resp. ${\rm val}(\beta)$) and
$|\!|w_{\alpha,n}|\!|_{{\rm val}(\alpha)}$
(resp. $|\!|w_{\beta,n}|\!|_{{\rm val}(\beta)}$) is bounded
independently of $n$\\  
2) $\forall\ n\geq 0$, $\forall\ m\geq 1$,
$\varphi^{-m}(w_{\alpha,n}\otimes e_{\alpha}+w_{\beta,n}\otimes
e_{\beta})\in {\rm
  Fil}^0((\Q(\mu_{p^m})\otimes_{\Q}E)[[t]]\otimes_ED(\alpha,\beta))$\\  
3) $\forall\ n\geq 1$, $\psi(w_{\alpha,n})=\alpha^{-1}w_{\alpha,n-1}$,
$\psi(w_{\beta,n})=\beta^{-1}w_{\beta,n-1}$.  
 
\begin{lem}\label{explicit} 
Let $w_{\alpha},\ w_{\beta}\in {\mathcal R}_E^+$, $m\in {\mathbb N}$
and $\mu_{\alpha},\ \mu_{\beta}$ the locally analytic distributions on
$\Z$ corresponding to $w_{\alpha},\ w_{\beta}$. The condition:  
$$\varphi^{-m}(w_{\alpha}\otimes e_{\alpha}+w_{\beta}\otimes
e_{\beta})\in {\rm
  Fil}^0((\Q(\mu_{p^m})\otimes_{\Q}E)[[t]]\otimes_ED(\alpha,\beta))$$  
is equivalent to the equalities:  
$$\alpha^m\int_{\Z}z^j\zeta_{p^m}^z\mu_{\alpha}(z)
=\beta^m\int_{\Z}z^j\zeta_{p^m}^z\mu_{\beta}(z)$$  
for all $j\in \{0,\cdots,k-2\}$ and all primitive $p^m$-roots $\zeta_{p^m}$ of $1$. 
\end{lem} 
\begin{proof} 
Here $\zeta_{p^m}^z$ is considered as a locally constant function
$\Z\rightarrow E$ (enlarging $E$ is necessary). Fixing $\zeta_{p^m}$ a
primitive root of $1$, we have:  
$$\varphi^{-m}(X)=(\zeta_{p^m}\otimes 1)e^{t/p^m}-1
=(\zeta_{p^m}\otimes 1)(e^{t/p^m}-1)+\zeta_{p^m}\otimes 1-1$$ 
in $(\Q(\mu_{p^m})\otimes_{\Q}E)[[t]]$ (see previous
lectures). Writing $w_{\alpha}=\sum_{i=0}^{+\infty}\alpha_iX^i$ and
$w_{\beta}=\sum_{i=0}^{+\infty}\beta_iX^i$, the condition on ${\rm
  Fil}^0$ is equivalent to:  
\begin{eqnarray*} 
\alpha^m\sum_{i=0}^{+\infty}\alpha_i\varphi^{-m}(X)^i 
e_{\alpha}+\beta^m\sum_{i=0}^{+\infty}\beta_i\varphi^{-m}(X)^ie_{\beta}\in\\
(\Q(\mu_{p^m})\otimes_{\Q}E)[[t]](e_{\alpha}+e_{\beta})\oplus   
(e^{t/p^m}-1)^{k-1}((\Q(\mu_{p^m})\otimes_{\Q}E)[[t]]e_{\alpha}) 
\end{eqnarray*} 
noting that $e^{t/p^m}-1$ generates ${\rm gr}^1(\Q[[t]])=\Q\overline
t$. Assuming $E$ contains $\Q(\mu_{p^m})$, we have: 
$$\Q(\mu_{p^m})\otimes_{\Q}E=\prod_{\Q(\mu_{p^m})\hookrightarrow E} E$$ 
and writing the binomial expansion of
$\varphi^{-m}(X)^i=(\zeta_{p^m}(e^{t/p^m}-1)+\zeta_{p^m}-1)^i$, a
simple calculation gives that the above condition is equivalent to the
equalities:  
\begin{eqnarray}\label{combin} 
\alpha^m\sum_{i=j}^{+\infty}\alpha_i{i\choose j}(\zeta_{p^m}-1)^{i-j}
=\beta^m\sum_{i=j}^{+\infty}\beta_i{i\choose j}(\zeta_{p^m}-1)^{i-j} 
\end{eqnarray} 
for $j\in \{0,\cdots,k-2\}$ and all primitive $p^m$-roots
$\zeta_{p^m}$ of $1$. Using $\alpha_i=\int_{\Z}{z\choose
  i}\mu_{\alpha}(z)$ and the elementary Mahler expansion:  
$${z\choose j}\zeta_{p^m}^{z-j}=\sum_{i=j}^{+\infty}{i\choose
  j}(\zeta_{p^m}-1)^{i-j}{z\choose i}$$  
we obtain: 
\begin{eqnarray*} 
\sum_{i=j}^{+\infty}\alpha_i{i\choose
  j}(\zeta_{p^m}-1)^{i-j}&=&\int_{\Z}
\Big(\sum_{i=j}^{+\infty}{z\choose i}{i\choose
  j}(\zeta_{p^m}-1)^{i-j}\Big)\mu_{\alpha}(z)\\  
&=&\zeta_{p^m}^{-j}\int_{\Z}{z\choose j}\zeta_{p^m}^{z}\mu_{\alpha}(z) 
\end{eqnarray*} 
and likewise with $\beta_i$ and $\mu_{\beta}$. 
Using (\ref{combin}), we easily get the result. 
\end{proof} 
 
For $n\in {\mathbb Z}_{\geq 0}$, let $\mu_{\alpha,n}$ be the tempered
distribution on $\Z$ of order ${\rm val}(\alpha)$ associated to
$\alpha^nw_{\alpha,n}$. We have
$\psi(\mu_{\alpha,n})=\mu_{\alpha,n-1}$ for $n\in {\mathbb N}$ and can
thus define ${\mu}_{\alpha}\in {\rm LPol}_{c,k-2}^*$ (algebraic dual
of ${\rm LPol}_{c,k-2}$) as before:  
\begin{eqnarray}\label{prolong} 
\int_{p^{-N}\Z}f(z)\mu_{\alpha}(z):=\int_{\Z}f(z/p^N)\mu_{\alpha,N}(z). 
\end{eqnarray} 
Likewise with $\beta$ instead of $\alpha$. By Lem.\ref{explicit}, we
have that condition 2) above on $(w_{\alpha,n},w_{\beta,n})_n$ is
equivalent to:  
\begin{eqnarray}\label{filtr} 
\alpha^{m-n}\int_{\Z}z^j\zeta_{p^m}^z\mu_{\alpha,n}(z)
=\beta^{m-n}\int_{\Z}z^j\zeta_{p^m}^z\mu_{\beta,n}(z)  
\end{eqnarray} 
for $j\in \{0,\cdots,k-2\}$, $n\geq 0$ and $m\geq 1$. 
 
Fixing embeddings $\overline{\mathbb Q}\hookrightarrow {\mathbb C}$,
$\overline{\mathbb Q}\hookrightarrow \overline\Q$ and using the
identification $\Q/\Z={\mathbb Z}[1/p]/{\mathbb Z}$, we can ``see''
$e^{2i\pi z}$ in $\overline\Q$ for any $z\in \Q$.  
 
\begin{cor}\label{equival} 
Condition 2) above is equivalent to the equalities in $\overline\Q$: 
$$\int_{p^{-N}\Z}z^je^{2i\pi
  zy}\mu_{\alpha}(z)=\Big(\frac{\alpha}{\beta}\Big)^{{\rm
    val}(y)}\int_{p^{-N}\Z}z^je^{2i\pi zy}\mu_{\beta}(z)$$  
for $j\in \{0,\cdots,k-2\}$, $y\in \Q^{\times}$ and $N>{\rm val}(y)$. 
\end{cor} 
\begin{proof} 
This follows from (\ref{prolong}) and (\ref{filtr}) noting that the
algebraic number $e^{2i\pi y/p^N}$ is a primitive $p^{N-{\rm
    val}(y)}$-root of $1$.  
\end{proof}

\section{Representations of $\G$ and $(\varphi,\Gamma)$-mod\-ules II (C.B.)} 
 
We still assume $V$ absolutely irreducible as in
\S\ref{banachdebut}. The aim of this section is to construct a
topological isomorphism $(\limproj_{\psi}D(V))^{\rm b}\simeq\Pi(V)^*$,
more precisely a topological isomorphism:  
$$(\limproj_{\psi}N(T))\otimes_{\O}E\simeq (B(\beta)/L(\beta))^*.$$ 
We then deduce from this isomorphism that $\Pi(V)$ is non zero,
topologically irreducible and admissible. These three statements were
conjectured in \cite{Br2} and \cite{Br3} (and already known in some
cases by a ``reduction modulo $p$'' method, see \cite{Br2} and Lecture
10). We also deduce a link with the Iwasawa cohomology $H^1_{\rm
  Iw}(\Q,V)$ of $V$.   
 
\subsection{The map $\Pi(V)^*\rightarrow (\lim_{\psi}D(V))^{\rm b}$}\label{unsens} 
We first construct a map from $(B(\beta)/L(\beta))^*$ to
the space $(\limproj_{\psi}N(T))\otimes_{\O}E$. Note that any $f$ in
$\pi(\alpha)$ or $\pi(\beta)$ which has compact support as a function
on $\Q$ via (\ref{groupeqp}) is an element in ${\rm LPol}_{{\rm
    c},k-2}$. Recall that $I:\pi(\alpha)\simeq \pi(\beta)$ is the
intertwining operator of (\ref{intertw}).  
 
\begin{lem}\label{crucial2} 
Let $\mu_{\alpha},\mu_{\beta}\in {\rm LPol}_{{\rm c},k-2}^*$. The
following statements are equivalent:\\  
1) For $j\in \{0,\cdots,k-2\}$, $y\in \Q^{\times}$ and $N>{\rm val}(y)$, we have: 
\begin{eqnarray}\label{rajout} 
\int_{p^{-N}\Z}z^je^{2i\pi zy}\mu_{\alpha}(z)
=\Big(\frac{\alpha}{\beta}\Big)^{{\rm val}(y)}\int_{p^{-N}\Z}z^je^{2i\pi zy}\mu_{\beta}(z). 
\end{eqnarray} 
2) For any $f\in \pi(\alpha)$ with compact support (as a function on
   $\Q$ via (\ref{groupeqp})) such that $I(f)\in \pi(\beta)$ has
   compact support (as a function on $\Q$ via (\ref{groupeqp})), we
   have:  
\begin{eqnarray}\label{debutinter} 
\int_{\Q}f(z)\mu_{\alpha}(z)=\frac{1-\frac{\beta}{\alpha}}
{1-\frac{\alpha}{p\beta}}\int_{\Q}I(f)(z)\mu_{\beta}(z).  
\end{eqnarray} 
\end{lem} 
\begin{proof} 
Note that $\alpha\ne \beta$ and $\alpha\ne p\beta$ (see
\S\ref{banachdebut}) hence the constant in (\ref{debutinter}) is well
defined and non zero. Recall the intertwining $I^{\rm sm}:{\rm
  Ind}_{\B}^{\G}{\rm unr}(\alpha^{-1})\otimes {\rm
  unr}(p\beta^{-1})\simeq {\rm Ind}_{\B}^{\G}{\rm
  unr}(\beta^{-1})\otimes {\rm unr}(p\alpha^{-1})$ is given by:  
\begin{eqnarray}\label{integrale2} 
I^{\rm sm}(h)= \bigg(g\mapsto
\int_{\Q}h\Big(\begin{pmatrix}0&-1\\1&x\end{pmatrix}g\Big)dx\bigg)  
\end{eqnarray} 
where $dx$ is Haar measure on $\Q$ (with values in ${\mathbb Q}\subset
E$). In terms of locally constant functions on $\Q$ via
(\ref{groupeqp}), we thus get:  
\begin{eqnarray}\label{integrale} 
I^{\rm sm}:f\mapsto \bigg(z\mapsto \int_{\Q}
(\alpha p\beta^{-1})^{{\rm val}(x)}f(z+x^{-1})dx=\\ 
\nonumber \int_{\Q}(\beta p\alpha^{-1})^{{\rm val}(x)}f(z+x)dx\bigg). 
\end{eqnarray} 
This is in fact an {\it algebraic} intertwining. This means that we
shouldn't bother too much about convergence problems in
$\int_{\Q}(\beta p\alpha^{-1})^{{\rm val}(x)}f(z+x)dx$ since, {\it in
  fine}, we can always replace the infinite sums around $0$ or
$-\infty$ by well defined algebraic expressions in $\beta
p\alpha^{-1}$. Let $j\in \{0,\cdots,k-2\}$ and $f_j:\Q\rightarrow E$
be a locally constant function with compact support such that $I^{\rm
  sm}(f_j)$ also has compact support. Let $\hat f_j$ be the usual
Fourier transform of $f_j$ which is again a locally constant function
on $\Q$ with compact support. Recall $\hat
f_j(x)=\int_{\Q}f_j(z)e^{-2i\pi zx}dz$ and $f_j(z)=\int_{\Q}\hat
f_j(x)e^{2i\pi zx}dx$ (see the end of \S\ref{distrback} for $e^{-2i\pi
  zx}$) where $dx$, $dz$ is Haar measure on $\Q$. The fact that
$I^{\rm sm}(f_j)$ has compact support is easily checked with
(\ref{integrale}) to be equivalent to $\hat f_j(0)=0$ (for $|z|\gg 0$,
we have $I^{\rm sm}(f_j)(z)=(\beta p\alpha^{-1})^{{\rm
    val}(z)}\int_{\Q}f_j(x)dx=(\beta p\alpha^{-1})^{{\rm val}(z)}\hat
f_j(0)$). Let $N\in {\mathbb N}$ be such that $f_j$ and $I^{\rm
  sm}(f_j)$ have support in $p^{-N}\Z$ and such that $\hat
f_j|_{p^N\Z}=0$. For $z\in p^{-N}\Z$, we have:  
\begin{eqnarray*} 
I(z^jf_j)(z)&=&z^j\int_{p^{-N}\Z}(\beta p\alpha^{-1})^{{\rm val}(x)}f_j(z+x)dx\\ 
&=&z^j\int_{p^{-N}\Z}(\beta p\alpha^{-1})^{{\rm val}(x)}
\bigg(\int_{\Q-p^N\Z}\hat f_j(y)e^{2i\pi y(z+x)}dy\bigg)dx\\ 
&=&z^j\int_{\Q-p^N\Z}\hat f_j(y)e^{2i\pi zy}\bigg(\int_{p^{-N}\Z}
(\beta p\alpha^{-1})^{{\rm val}(x)}e^{2i\pi xy}dx\bigg)dy. 
\end{eqnarray*} 
Decomposing $p^{-N}\Z=p^{-N}\Z^{\times}\amalg p^{-N+1}
\Z^{\times}\amalg\cdots$, a straightforward computation yields 
for $N>{\rm val}(y)$: 
$$\int_{p^{-N}\Z}\!\!\!\!\!(\beta p\alpha^{-1})^{{\rm val}(x)}e^{2i\pi
  xy}dx=\!\!\sum_{i=-N}^{+\infty}(\beta
p\alpha^{\!-1})^i\!\!\!\int_{p^i\Z^{\times}}\!\!\!e^{2i\pi
  xy}dx=\frac{1-\frac{\alpha}{p\beta}}{1-\frac{\beta}
{\alpha}}\Big(\frac{\alpha}{\beta}\Big)^{{\rm val}(y)}$$  
hence: 
\begin{eqnarray}\label{fourier1} 
I(z^jf_j)(z)=\frac{1-\frac{\alpha}{p\beta}}{1-\frac{\beta}
{\alpha}}\int_{\Q-p^N\Z}\hat f_j(y)z^je^{2i\pi zy}\Big(\frac{\alpha}
{\beta}\Big)^{{\rm val}(y)}dy 
\end{eqnarray} 
for $z\in p^{-N}\Z$ and $I(z^jf_j)(z)=0$ otherwise. Likewise, we have:  
\begin{eqnarray}\label{fourier2} 
z^jf_j(z)=\int_{\Q-p^N\Z}\hat f_j(y)z^je^{2i\pi zy}dy 
\end{eqnarray} 
for $z\in p^{-N}\Z$ and $z^jf_j(z)=0$ otherwise. Note that
(\ref{fourier1}) and (\ref{fourier2}) are in fact finite sums over the
same finite set of values of $y$ (as $\hat f_j$ is locally constant
with compact support). Replacing $I(z^jf_j)(z)$ and $z^jf_j(z)$ in
(\ref{variant}) below by the finite sums (\ref{fourier1}) and
(\ref{fourier2}), we then easily deduce that the following equalities
for various $f_j$ as above with support in $p^{-N}\Z$:  
\begin{eqnarray}\label{variant} 
\int_{p^{-N}\Z}I(z^jf_j(z))\mu_{\beta}(z)=\frac{1-\frac{\alpha}
{p\beta}}{1-\frac{\beta}{\alpha}}\int_{p^{-N}\Z}z^jf_j(z)\mu_{\alpha}(z) 
\end{eqnarray} 
are equivalent to the equalities in 1) (final details here are 
left to the reader). This finishes the proof. 
\end{proof} 
 
Let $\mu_{\beta}\in (B(\beta)/L(\beta))^{*}\subset \pi(\beta)^*$. 
Define a sequence of locally analytic distributions on $\Z$, 
$(\mu_{\beta,N})_{N\in {\mathbb Z}_{\geq 0}}$, as follows: 
\begin{eqnarray*} 
\int_{\Z}f(z)\mu_{\beta,N}(z):=\int_{p^{-N}\Z}f(p^Nz)\mu_{\beta}(z). 
\end{eqnarray*} 
One obviously has $\psi(\mu_{\beta,N})=\mu_{\beta,N-1}$ by Lem.\ref{Psidis}. 
 
\begin{lem}\label{bounded} 
The distributions $\mu_{\beta,N}$ are tempered of order ${\rm
  val}(\beta)$ and $|\!|\beta^{-N}\mu_{\beta,N}|\!|'_{{\rm
    val}(\beta)}$ is bounded independently of $N$ (see \S\ref{Yvette}
for the definition of this norm).  
\end{lem} 
\begin{proof}
The fact that $\mu_{\beta,N}$ is tempered of order ${\rm val}(\beta)$
is easy and left to the reader. Going back to the proof of
Th.\ref{compl}, we know that $\mu_{\beta}$ satisfies in particular
(\ref{chaud1}) and (\ref{chaud2}). Hence, for $a\in \Z$, $0\leq j\leq
k-2$ and $n\in {\mathbb Z}_{\geq 0}$:  
\begin{eqnarray*} 
\beta^{-N}\int_{a+p^n\Z}(z-a)^j\mu_{\beta,N}(z)&=&\beta^{-N}p^{Nj}
\int_{p^{-N}a+p^{n-N}\Z}(z-p^{-N}a)^j\mu_{\beta}(z)\\ 
&\in & C_{\mu_{\beta}}p^{-N{\rm val}(\beta)}p^{Nj}p^{(n-N)(j-{\rm val}(\beta))}\O\\ 
&\in &C_{\mu_{\beta}}p^{n(j-{\rm val}(\beta))}\O. 
\end{eqnarray*}
For any locally analytic distribution $\mu$ on $\Z$ and any $r\in
{\mathbb R}^+$ such that $r<k-1$, define: 
$$|\!|\mu|\!|_{r,k}':={\rm sup}_{a\in\Z}
{\rm sup}_{j\in\{0,\cdots,k-2\}}{\rm sup}_{n\in{\mathbb Z}_{\geq 0}}
p^{n(j-r)}\left|\int_{a+p^n\Z}(z-a)^j\mu(z)\right|.$$ 
An examination of the proof of Lem.\ref{etendre} shows that the norms
$|\!|\cdot |\!|_{r}'$ and $|\!|\cdot |\!|_{r,k}'$ are in fact
equivalent (see also \cite{Co2},\S V.3.6). This implies: 
$$|\!|\beta^{-N}\mu_{\beta,N}|\!|'_{{\rm val}(\beta)}\leq
c|\!|\beta^{-N}\mu_{\beta,N}|\!|'_{{\rm val}(\beta),k}\leq
c|C_{\mu_{\beta}}|$$  
for a convenient constant $c\in {\mathbb R}^+$.  
\end{proof} 
 
By Cor.\ref{1ercas}, the intertwining $I:\pi(\alpha)\simeq \pi(\beta)$
extends ``by continuity'' to an intertwining $\hat
I:B(\alpha)/L(\alpha)\simeq B(\beta)/L(\beta)$ and we define:
$$\mu_{\alpha}:=\frac{1-\frac{\beta}{\alpha}}{1-\frac{\alpha}{p\beta}}\mu_{\beta}\circ
\hat I \in (B(\alpha)/L(\alpha))^{*}.$$  
By restriction to ${\rm LPol}_{{\rm c},k-2}$, $\mu_{\alpha}$ and
$\mu_{\beta}$ define elements of ${\rm LPol}_{{\rm c},k-2}^*$
satisfying (\ref{debutinter}) by definition, hence (\ref{rajout}) by
Lem.\ref{crucial2}. As in Lem.\ref{bounded}, $\mu_{\alpha}$ gives rise
to a bounded sequence of tempered distributions $\mu_{\alpha,N}$ on
$\Z$. From Cor.\ref{equival} and Lem.\ref{bounded}, we get that Amice
transform produces from $(\mu_{\alpha,N})_N$ and $(\mu_{\beta,N})_N$ a
sequence $(v_N)_N=(w_{\alpha,N}\otimes e_{\alpha}+w_{\beta,N}\otimes
e_{\beta})_N$ satisfying conditions 1), 2), 3) before
Lem.\ref{explicit}, hence $(v_N)_N\in
(\limproj_{\psi}N(T))\otimes_{\O}E$. This defines a linear map
$(B(\beta)/L(\beta))^{*}\rightarrow
(\limproj_{\psi}N(T))\otimes_{\O}E$. The continuity follows from the
fact that the function ``Amice transform'' from the space of tempered
distributions of order $r$ on $\Z$ with norm smaller than $1$ endowed
with the weak topology of pointwise convergence to the subspace of
${\mathcal R}_E^+$ of elements $w$ of order $r$ such that
$|\!|w|\!|_r\leq 1$ endowed with its compact topology induced by that
of ${\mathcal R}_E^+$ is a topological isomorphism. Final details here
are left to the reader.  
 
\subsection{The map $(\lim_{\psi}D(V))^{\rm b}\rightarrow \Pi(V)^*$}\label{finiso} 
 
We now construct a map from $(\limproj_{\psi}N(T))\otimes_{\O}E$ to
$(B(\beta)/L(\beta))^{*}$. Let $(v_n)_n=(w_{\alpha,n}\otimes
e_{\alpha}+w_{\beta,n}\otimes e_{\beta})_n\in
(\limproj_{\psi}N(T))\otimes_{\O}E$ and define $\mu_{\alpha}$ and
$\mu_{\beta}$ as in \S\ref{distrback}. They are elements of ${\rm
  Pol}_{{\rm c},k-2}^*$ satisfying (\ref{rajout}) by
Cor.\ref{equival}, hence (\ref{debutinter}) by Lem.\ref{crucial2}.  
 
\begin{lem}\label{versfin} 
There is a unique way to extend $\mu_{\beta}$ and $\mu_{\alpha}$ as
elements of respectively $B(\beta)^*$ and $B(\alpha)^*$ such that for
any $f\in \pi(\alpha)$ (seen as a function on $\Q$ via
(\ref{groupeqp})):  
\begin{eqnarray}\label{extension} 
\int_{\Q}f(z)\mu_{\alpha}(z)=\frac{1-\frac{\beta}{\alpha}}
{1-\frac{\alpha}{p\beta}}\int_{\Q}I(f)(z)\mu_{\beta}(z). 
\end{eqnarray} 
\end{lem} 
\begin{proof} 
An easy computation using (\ref{integrale}) gives: 
\begin{eqnarray}\label{calculc} 
I(z^j{\mathbf 1}_{\Z})&=&\frac{1-\frac{1}{p}}{1-\frac{\beta}
{\alpha}}z^j{\mathbf 1}_{\Z}+z^j\Big(\!\frac{p\beta}
{\alpha}\!\Big)^{\!\!{\rm val}(z)}\!\!\!\!(1-{\mathbf 1}_{\Z})\\
\label{calculnc}
I\big(z^j\Big(\!\frac{p\alpha}{\beta}\!\Big)^{\!\!{\rm val}
(z)}\!\!\!\!(1-{\mathbf 1}_{p\Z})\big)&=&z^j{\mathbf
  1}_{p\Z}+\frac{1-\frac{1}{p}}{1-\frac{\beta}{\alpha}}z^j
\Big(\!\frac{p\beta}{\alpha}\!\Big)^{\!\!{\rm val}(z)}\!\!\!\!(1-{\mathbf 1}_{p\Z}). 
\end{eqnarray} 
Since $\int_{\Q}z^j{\mathbf 1}_{\Z}\mu_{\alpha}(z)$ and
$\int_{\Q}z^j{\mathbf 1}_{\Z}\mu_{\beta}(z)$ are well defined, we see
using (\ref{extension}) and (\ref{calculc}) that
$$\int_{\Q}z^j\big(\!\frac{p\beta}{\alpha}\!\big)^{\!\!{\rm
    val}(z)}\!(1-{\mathbf 1}_{\Z})\mu_{\beta}(z)$$ is uniquely
determined. Then, using (\ref{calculnc}) (and (\ref{extension})), we
see that $\int_{\Q}z^j\big(\!\frac{p\alpha}{\beta}\!\big)^{\!\!{\rm
    val}(z)}\!(1-{\mathbf 1}_{p\Z})\mu_{\alpha}(z)$ is also uniquely
determined. One readily checks that this defines unique extensions of
$\mu_{\beta}$ and $\mu_{\alpha}$ as linear forms on respectively
$\pi(\beta)$ and $\pi(\alpha)$. To check that $\mu_{\beta}$ is in fact
in $B(\beta)^*$, it is enough by Cor.\ref{travail} to check that the
relevant Amice-V\'elu condition (as in Lem.\ref{etendre}) is satisfied
by $\mu_{\beta}$ on each copy of $\Z$ via (\ref{f1f2}). We already
know it for the first one since $w_{\beta,0}$ is of order ${\rm
  val}(\beta)$. For the second one, a straightforward computation via
the definition of $f_2$ in (\ref{f1f2}) (exchanging $\alpha$ and
$\beta$) shows that we have to find a constant $C_{\mu_{\beta}}\in E$
such that:  
\begin{eqnarray}\label{ch1} 
\int_{a^{-1}+p^{n-2{\rm
      val}(a)}\Z}\Big(\!\frac{p\beta}{\alpha}\!\Big)^{\!\!{\rm val}
(z)}z^{k-2-j}(1-az)^j\mu_{\beta}(z)\in C_{\mu_{\beta}}
p^{n(j-{\rm val}(\beta))}\O 
\end{eqnarray} 
for $a\in \Z-\{0\}$, $j\in\{0,\cdots,k-2\}$ and $n>{\rm val}(a)$, and: 
\begin{eqnarray}\label{ch2} 
\int_{\Q-p^{-n}\Z}\Big(\!\frac{p\beta}{\alpha}\!\Big)^{\!\!{\rm
    val}(z)}z^{k-2-j}\mu_{\beta}(z)\in C_{\mu_{\beta}}p^{n(j-{\rm
    val}(\beta))}\O  
\end{eqnarray} 
for $j\in\{0,\cdots,k-2\}$ and $n\in {\mathbb Z}_{\geq 0}$. From
condition 1) before Lem.\ref{explicit}, we know there exists
$C_{\mu_{\beta}}\in E$ such that
$\beta^{-N}\int_{a+p^n\Z}(z-a)^j\mu_{\beta,N}(z)\in
C_{\mu_{\beta}}p^{n(j-{\rm val}(\beta))}\O$ for $a\in \Z$ and $N,n\in
{\mathbb Z}_{\geq 0}$. For $a\in \Q^{\times}$, let $u_a:=p^{{\rm
    val}(a)}a^{-1}\in \Z^{\times}$. By (\ref{prolong}), the left hand
side of (\ref{ch1}) is equal to (up to a unit in $\O$):  
$$\Big(\!\frac{p\beta}{\alpha}\!\Big)^{\!\!-{\rm
    val}(a)}p^{(j-k+2){\rm val}(a)}\int_{u_a+p^{n-{\rm
      val}(a)}\Z}\!\!\!\!z^{k-2-j}(z-u_a)^j\mu_{\beta,{\rm
    val}(a)}(z)$$  
and we see, writing $z^{k-2-j}=(z-u_a+u_a)^{k-2-j}=\cdots$ and using
    the above bound with $N={\rm val}(a)$, that it belongs to:  
$$C_{\mu_{\beta}}p^{({\rm val}(\alpha)-{\rm val}(\beta)-1){\rm
	val}(a)}p^{(j-k+2){\rm val}(a)}p^{{\rm val}(\beta){\rm
	val}(a)}p^{(n-{\rm val}(a))(j-{\rm val}(\beta))}\O$$  
which is exactly $C_{\mu_{\beta}}p^{n(j-{\rm val}(\beta))}\O$ using
${\rm val}(\alpha)+{\rm val}(\beta)=k-1$. This proves (\ref{ch1}). To
prove (\ref{ch2}), one uses the equality (similar to (\ref{calculc})):  
\begin{eqnarray*} 
\Big(\frac{p\beta}{\alpha}\!\Big)^{{\rm val}(z)}z^{k-2-j}(1-{\mathbf
  1}_{p^{-n}\Z})&=&p^{-n}I(z^{k-2-j}{\mathbf 1}_{p^{-n}\Z})-\\  
&&\frac{1-\frac{1}{p}}{1-\frac{\beta}{\alpha}}\Big(\frac{\alpha}{p\beta}\Big)^nz^{k-2-j}{\mathbf
  1}_{p^{-n}\Z}  
\end{eqnarray*} 
together with (\ref{extension}), the equalities:  
\begin{eqnarray*} 
p^{-n}\int_{p^{-n}\Z}z^{k-2-j}\mu_{\alpha}(z)&=&p^{(j-k+1)n}
\int_{\Z}z^{k-2-j}\mu_{\alpha,n}(z)\\    
\Big(\!\frac{\alpha}{p\beta}\!\Big)^n\int_{p^{-n}\Z}z^{k-2-j}
\mu_{\beta}(z)&=&\Big(\!\frac{\alpha}{p\beta}\!\Big)^np^{(j-k+2)n}
\int_{\Z}z^{k-2-j}\mu_{\beta,n}(z) 
\end{eqnarray*} 
and the fact $\alpha^{-n}\int_{\Z}z^{k-2-j}\mu_{\alpha,n}(z)$ and
$\beta^{-n}\int_{\Z}z^{k-2-j}\mu_{\beta,n}(z)$ are bounded
independently of $n$. The proof that $\mu_{\alpha}$ extends as a
distribution on $B(\alpha)$ is similar.  
\end{proof} 
 
By Lem.\ref{versfin}, we have a continuous injective map
$(\limproj_{\psi}N(T))\otimes_{\O}E\hookrightarrow B(\beta)^{*}$,
$(w_{\alpha,n}\otimes e_{\alpha}+w_{\beta,n}\otimes
e_{\beta})_n\mapsto \mu_{\beta}$ (the continuity follows as in
\S\ref{unsens} from the continuity of Amice transform and the
injectivity is straightforward using Lem.\ref{versfin}). With the
action of $\G$ deduced from (\ref{actionjoe}) on the dual
$B(\beta)^*$, it is an easy and pleasant exercise on Amice transform
that we leave to the reader (using e.g. the formula
$\int_{\Z}(1+X)^z\mu(z)$) to check that the induced action of
$\begin{pmatrix}1&\Z\\0&\Q^{\times}\end{pmatrix}$ on
$(\limproj_{\psi}N(T))\otimes_{\O}E$ is indeed as in Th.\ref{main}. In
particular, it preserves $\limproj_{\psi}N(T)$.  
 
\begin{lem}\label{versfin2} 
The distribution $\mu_{\beta}\in B(\beta)^*$ defined in
Lem.\ref{versfin} sits in $(B(\beta)/L(\beta))^*$.  
\end{lem} 
\begin{proof} 
The non zero compact $\O$-submodule of $B(\beta)^*$ which is the image
of $\limproj_{\psi}N(T)$ is preserved by the action of $\B$ since
$\begin{pmatrix}1&\Z\\0&\Q^{\times}\end{pmatrix}$ generates
$\begin{pmatrix}1&\Q\\0&\Q^{\times}\end{pmatrix}$ and the central
character is integral. Let:  
$$\pi(\beta)^0:=\left\{f\in \pi(\beta),\ \int_{\Q}f(z)\mu(z)\in \O\
\forall\ \mu\in\limproj_{\psi}N(T)\right\}.$$  
We have that $\pi(\beta)^0\!\!\subset \pi(\beta)$ is thus preserved by
$\B$. By \cite{Schn},Lem.13.1(iii), we also have that
$\pi(\beta)^0\otimes_{\O}E=\pi(\beta)$ (using that $B(\beta)^*\hookrightarrow \pi(\beta)^*$ and thus that $\limproj_{\psi}N(T)$ is bounded inside $\pi(\beta)^*$). By the same proof as for Th.\ref{compl} (using that any
$\mu\in \limproj_{\psi}N(T)$ is such that $\forall\ g\in\B$, $\forall\
f\in \pi(\beta)^0$, $|\mu(g(f))|\leq 1$), we deduce that any $\mu\in
B(\beta)^*$ which is in the image of $\limproj_{\psi}N(T)$ cancels
$L(\beta)$. Whence the result.  
\end{proof} 
 
By Lem.\ref{versfin2}, the above continuous injection
$(\limproj_{\psi}N(T))\otimes_{\O}E\hookrightarrow B(\beta)^*$ takes
values in the subspace $(B(\beta)/L(\beta))^*$. It is easily checked
that this map is an inverse to the previous one (\S\ref{unsens}), thus
finishing the proof of Th.\ref{main}. In particular, we have an
isomorphism:  
$$D(V)^{\psi=1}=\big((\limproj_{\psi}D(V))^{\rm b}\big)^{\psi=1}\simeq
{\Pi(V)^*}^{\begin{pmatrix}1&0\\0&p^{\mathbb Z}\end{pmatrix}}.$$  
Using Fontaine's canonical isomorphism (see e.g. \cite{Co2},\S4.3): 
\begin{eqnarray}\label{fontaine} 
H^1_{\rm Iw}(\Q,V):=\big(\limproj_n H^1({\rm
  Gal}(\overline\Q/\Q(\mu_{p^n})),
T)\big)\otimes_{\O}E\simeq D(V)^{\psi=1} 
\end{eqnarray} 
for any Galois lattice $T$ in $V$, we deduce the quite surprising corollary: 
 
\begin{cor}\label{brfoiw}
We have an isomorphism of $\O[[\Z^{\times}]]$-modules: 
$$H^1_{\rm Iw}(\Q,V)\simeq {\Pi(V)^*}^{\begin{pmatrix}
1&0\\0&p^{\mathbb Z}\end{pmatrix}}$$ 
where the $\O[[\Z^{\times}]]$-module structure on the left hand side
is coming from the action of $\Gamma$ and on the right hand side from
the action of $\begin{pmatrix}1&0\\0&\Z^{\times}\end{pmatrix}$.  
\end{cor} 
 
When $V$ is reducible (as in 2) and 3) of \S\ref{banachdebut}), it is
not true that we have an isomorphism $(\limproj_{\psi}D(V))^{\rm
  b}\simeq \Pi(V)^*$ (although it is ``almost true''). But if
$\alpha\ne p^{k-1}$ and $\beta\ne 1$ (which always holds if $V$ is
coming from a modular form by the Weil conjectures), one can prove
that there is still an isomorphism of $\O[[\Z^{\times}]]$-modules
$H^1_{\rm Iw}(\Q,V)\simeq {\Pi(V)^*}^{\begin{pmatrix}1&0\\0&p^{\mathbb
      Z}\end{pmatrix}}$.

\subsection{Irreducibility and admissibility}\label{results} 
 
We deduce from Th.\ref{main} that, for $V$ irreducible as before,
$\Pi(V)$ is non zero, topologically irreducible and admissible.  
 
\begin{cor} 
The Banach $B(\alpha)/L(\alpha)$ and $B(\beta)/L(\beta)$ are non zero. 
\end{cor} 
\begin{proof} 
This follows from $(\limproj_{\psi}D(V))^{\rm b}\ne 0$ which follows
from $D(V)^{\psi=1}\ne 0$ which follows e.g. from (\ref{fontaine}).  
\end{proof} 
 
\begin{cor}\label{irreductible} 
The Banach $B(\alpha)/L(\alpha)$ and $B(\beta)/L(\beta)$ are topologically irreducible. 
\end{cor} 
\begin{proof} 
Assume there exists a (necessarily unitary) $\G$-Banach $V$ together
with a $\G$-equivariant topological surjection
$B(\alpha)/L(\alpha)\twoheadrightarrow V$. We can find norms on the
two Banach such that we have a surjection of unit balls
$(B(\alpha)/L(\alpha))^0\twoheadrightarrow V^0$. By duality, we get a
$\G$-equivariant continuous injection of compact $\O$-modules:  
$${\rm Hom}_{\O}(V^0,\O)\hookrightarrow {\rm
  Hom}_{\O}((B(\alpha)/L(\alpha))^0,\O).$$  
We have seen that $\limproj_{\psi}N(T)$ is a compact $\O$-lattice of
  $(B(\alpha)/L(\alpha))^*$ which is preserved by $\B$. Hence,
  multiplying by a scalar in $E^{\times}$ if necessary, we get a
  continuous $\B$-equivariant injection of compact $\O$-modules ${\rm
  Hom}_{\O}(V^0,\O)\hookrightarrow \limproj_{\psi}N(T)$. If $V\ne 0$,
  Cor.\ref{versirr} together with the translation of the action of
  $\B$ in terms of $\psi$ and $\Gamma$ imply ${\rm
  Hom}_{\O}(V^0,\O)\otimes E\buildrel\sim\over\rightarrow
  (B(\alpha)/L(\alpha))^*$, which implies
  $B(\alpha)/L(\alpha)\buildrel\sim\over\rightarrow V$. Hence
  $B(\alpha)/L(\alpha)$ (and $B(\beta)/L(\beta)$ by Cor.\ref{1ercas})
  is topologically irreducible.  
\end{proof} 
 
We see from the previous proof that $B(\alpha)/L(\alpha)$ and
$B(\beta)/L(\beta)$ are even topologically irreducible as
$\B$-representations.  
 
\begin{cor}\label{admissible} 
The Banach $B(\alpha)/L(\alpha)$ and $B(\beta)/L(\beta)$ are admissible. 
\end{cor} 
\begin{proof} 
We don't know if the compact $\O$-module $\limproj_{\psi}N(T)$ is
preserved by $\K$ inside $$(\limproj_{\psi}D(V))^{\rm b}\simeq
(B(\alpha)/L(\alpha))^*,$$ but we can replace it by ${\mathcal
  M}:=\cap_{g\in \K}g(\limproj_{\psi}N(T))\subset \limproj_{\psi}N(T)$
which is a non zero compact $\O[[X]]$-submodule as in
Cor.\ref{versirr} (since it is easily checked to be preserved by
$B(\Q)$ using the Iwasawa decomposition of $\G$). Note that there are
now two natural structures of $\O[[X]]$-module on ${\mathcal M}$ : the
first is the one already defined and the second is:  
$$(\lambda,v)\in \O[[X]]\times {\mathcal M}\mapsto
\begin{pmatrix}0&1\\1&0\end{pmatrix}\lambda\begin{pmatrix}0&1\\1&0\end{pmatrix}v.$$
  Equivalently, the first is such that multiplication by $(1+X)^{\Z}$
  corresponds to the action of $\begin{pmatrix}1&\Z\\0&1\end{pmatrix}$
    and the second is such that multiplication by $(1+X)^{\Z}$
    corresponds to the action of
    $\begin{pmatrix}1&0\\\Z&1\end{pmatrix}$. We have seen that
      ${\mathcal M}=\limproj_{\psi}M$ (see proof of Cor.\ref{versirr})
      with $M$ of finite type over $\O[[X]]$, $M\subset N(T)$ and
      $\psi$ surjective on $M$. Let ${\rm pr}_0:{\mathcal
	M}\twoheadrightarrow M$ be the projection on the first
      component and define ${\mathcal N}:={\rm Ker}({\rm
	pr}_0)\subsetneq {\mathcal M}$. The map ${\mathcal N}\mapsto
      M$, $v\mapsto {\rm
	pr}_0\big(\begin{pmatrix}0&1\\1&0\end{pmatrix}v\big)$ is
	injective: if $v$ maps to $0$, its associated distribution
	$\mu\in B(\alpha)^*\simeq {\mathcal C}^{{\rm
	    val}(\alpha)}(\Z,E)^*\oplus {\mathcal C}^{{\rm
	    val}(\alpha)}(\Z,E)^*$ (see (\ref{f1f2})) is zero on both
	copies of ${\mathcal C}^{{\rm val}(\alpha)}(\Z,E)$ hence is
	zero in $(B(\alpha)/L(\alpha))^*$. Thinking in terms of
	distributions again, we also see that ${\mathcal N}$ is an
	$\O[[X]]$-module for the first structure but only a
	$\varphi(\O[[X]])$-module for the second structure (recall
	from (\ref{actionjoe}) and Th.\ref{main} that multiplication
	by $\varphi(1+X)=(1+X)^p$ on distributions correspond to the
	change of variables $z\mapsto z+p$ on functions). Moreover,
	for this second action, the above injection ${\mathcal
	  N}\hookrightarrow M$ is $\varphi(\O[[X]])$-linear. Since $M$
	is of finite type over $\O[[X]]$, hence over
	$\varphi(\O[[X]])$, we thus get that the
	$\varphi(\O[[X]])$-module ${\mathcal N}$ for the second action
	of $\varphi(\O[[X]])$ is of finite type. Now fix elements
	$(e_1,\cdots,e_m)\in {\mathcal M}$ (resp. $(f_1,\cdots,f_n)\in
	{\mathcal N}$) such that ${\rm pr}_0(e_i)$ (resp. $f_i$)
	generate $M$ over $\O[[X]]$ (resp. $\mathcal N$ over
	$\varphi(\O[[X]])$). Let $v\in {\mathcal M}$. There exist
	$\lambda_1,\cdots,\lambda_m$ in $\O[[X]]$ such that $v-\sum
	\lambda_ie_i\in {\mathcal N}$. Then there exist
	$\mu_1,\cdots,\mu_n$ in $\varphi(\O[[X]])$ such that $v-\sum
	\lambda_ie_i=\sum
	\begin{pmatrix}0&1\\1&0\end{pmatrix}\mu_i\begin{pmatrix}0&1\\1&0\end{pmatrix}f_i$. 
Since the $\lambda_i$ correspond to the action of elements in the
group algebra of $\begin{pmatrix}1&\Z\\0&1\end{pmatrix}$ and the
  $\begin{pmatrix}0&1\\1&0\end{pmatrix}
\mu_i\begin{pmatrix}0&1\\1&0\end{pmatrix}$ 
correspond to the action of elements in the group algebra 
of $\begin{pmatrix}1&0\\p\Z&1\end{pmatrix}$, we see in 
particular that $\mathcal M$ is of finite type over the 
group algebra of $\K$, whence the admissibility.      
\end{proof}

\section{Reduction modulo $p$ (L.B.)} 
In the $\ell$-adic case ($\ell\ne p$), Vign\'eras has proved that one
can reduce the classical local Langlands correspondence over
$\overline{{\mathbb Q}_{\ell}}$ modulo the maximal ideal of
$\overline{{\mathbb Z}_{\ell}}$ and thus obtain a new correspondence
between semi-simple representations over $\overline{{\mathbb
    F}_{\ell}}$ (\cite{Vi}).  
 
It is tempting to do the same here with $V$ and $\Pi(V)$. That is, for
any $V$ as in \S\ref{banachdebut}, we have associated a non zero
unitary admissible $\G$-Banach $\Pi(V)$. Let $\Pi(V)^0\subset \Pi(V)$
(resp. $T\subset V$) be any unit ball (resp. any $\O$-lattice) which
is preserved by $\G$ (resp. $\g$). In this last lecture, we want to
state and prove some cases of a conjecture relating the
semi-simplification of $\Pi(V)^0\otimes_{\O}{\mathbb F}_E$ and the
semi-simplification of $T\otimes_{\O}{\mathbb F}_E$.   
 
\subsection{Statement of the conjecture} 
In order to state the conjecture, we need to make  
lists of some $2$-dimensional $\overline\F$-representations  
of $\g$ and of some $\overline\F$-representations of $\G$.  
 
Let us start with those of $\g$. Let $\omega$ be the $\mod{p}$  
cyclotomic character seen as a character of $\Q^{\times}$ via class field theory, 
and let $\mu_{-1}$ be the unramified quadratic character of $\Q^{\times}$. 
We write $\i$ for the inertia subgroup of $\g$. Let $\omega_2  
: \i \ra \overline\F^{\times}$ be Serre's fundamental character of level $2$  
and for $s \in \{ 0,\cdots,p\}$, let ${\rm ind}(\omega_2^s)$  
be the unique (irreducible) representation of $\g$ whose determinant  
is $\omega^s$ and whose restriction to the inertia subgroup of  
${\rm Gal}(\overline{\Q}/\mathbb{Q}_{p^2})$ is $\omega_2^s \oplus \omega_2^{ps}$. 
We know that as $r$ runs over $0,\dots,p-1$ and as $\eta$ runs over all  
characters $\eta: \Q^{\times} \ra \overline\F^{\times}$, the representations 
$\rho(r,\eta):= ({\rm ind}(\omega_2^r)) \otimes \eta$ run over all  
irreducible $2$-dimensional $\overline\F$-representations of $\g$  
and that the only isomorphisms are: 
\begin{align*} 
\rho(r,\eta) & \simeq \rho(r,\eta \mu_{-1}) \\ 
\rho(r,\eta) & \simeq \rho(p-1-r,\eta \omega^r) \\ 
\rho(r,\eta) & \simeq \rho(p-1-r,\eta \omega^r \mu_{-1}) 
\end{align*} 
 
Now, we turn to representations of $\G$.  
For $r \in \{0,\cdots,p-1\}$, let  
${\rm Sym}^r \overline\F^2$ denote the natural representation of $\K$ 
acting through ${\rm GL}_2(\F)$ which we extend to $\K \Q^{\times}$ by sending $p$ to $1$. 
Let \[ {\rm ind}_{\K \Q^{\times}}^{\G}  
{\rm Sym}^r \overline\F^2 \] be the $\overline\F$-vector space of functions $f:\G 
\ra {\rm Sym}^r \overline\F^2$ which are compactly supported modulo $\Q^{\times}$ 
and such that $f(kg)={\rm Sym}^r(k)(f(g))$ with 
the left action of $\G$.  
 
For $r \in \{0,\cdots,p-1\}$,   
and $\lambda \in \overline\F$, let  
\[ \pi(r,\lambda,\eta) := \left[ \left( {\rm ind}_{\K \Q^{\times}}^{\G}  
{\rm Sym}^r \overline\F^2 \right) / (T-\lambda) \right]  
\otimes (\eta \circ \det), \] 
where $T$ is a ``Hecke operator'' which corresponds to the double 
coset \[ \K\Q^{\times} \begin{pmatrix} 1 & 0 \\ 0 & p \end{pmatrix} \K. \] 
 
Those representation are irreducible, 
and all smooth irreducible $\overline\F$-representations of $\G$  having 
a central character are given by (see \cite{BL1},\cite{BL},\cite{Br1}): 
\begin{enumerate} 
\item $\eta \circ \det$ where $\eta:\Q^{\times} 
\ra \overline\F^{\times}$ is a smooth character. 
\item the ${\rm Sp} \otimes (\eta \circ \det)$ where ${\rm Sp}$ is  
the ``Special'' representation. 
\item the $\pi(r,\lambda,\eta)$  
where $\eta:\Q^{\times} \ra \overline\F^{\times}$ is a smooth character 
and $r \in \{0,\cdots,p-1\}$ and  
$\lambda \in \overline\F \setminus \{-1,1\}$. 
\end{enumerate} 
 
When $\lambda \neq 0$, the $\pi(r,\lambda,\eta)$ have no non  
trivial intertwinings (unless $\lambda = \pm 1$, in which case 
there are some ``easy'' ones, see \cite{BL}) and when $\lambda=0$, 
the only isomorphisms are: 
\begin{align*} 
\pi(r,0,\eta) & \simeq \pi(r,0,\eta \mu_{-1}) \\ 
\pi(r,0,\eta) & \simeq \pi(p-1-r,0,\eta \omega^r) \\ 
\pi(r,0,\eta) & \simeq \pi(p-1-r,0,\eta \omega^r \mu_{-1}). 
\end{align*} 
 
This certainly suggests that there is a correspondence between the  
representations $\rho$ of $\g$ and the representations $\pi$ of $\G$. 
More precisely, we define the following ``correspondence'' (\cite{Br1}): 
\begin{definit}\label{corrmodp} 
If $\lambda=0$, then  
\[ \pi(r,0,\eta) \leftrightarrow \rho(r,\eta) \]  
and if $\lambda \neq 0$, then 
\[ \pi(r,\lambda,\eta)^{\rm ss} \oplus  
\pi([p-3-r],\lambda^{-1},\omega^{r+1} \eta)^{\rm ss}  
\leftrightarrow 
\begin{pmatrix} 
{\rm unr}(\lambda^{-1}) \omega^{r+1} & 0 \\ 0 & {\rm unr}(\lambda)
\end{pmatrix} \otimes \eta \] 
where ${\rm unr}(\lambda)$ is the unramified character sending Frobenius  
to $\lambda$, ${\rm ss}$ denotes the semisimplification and $[x]$ is  
the integer  in $\{0,\cdots,p-2\}$ which is congruent to $x \mod{p-1}$. 
\end{definit} 
 
On the other hand, we have associated in the previous  
chapters a unitary representation $\Pi(V)$ of  
$\G$ to a crystalline representation 
$V$ of $\g$. Let $V \mod{p}$ denotes the semisimplification of  
$T/pT$ where $T$ is some Galois invariant lattice of $V$ and  
let $\overline{\Pi}(V)$ be the semisimplification of $\Pi^0(V)/p$ 
where $\Pi^0(V)$ is some $\G$-invariant unit ball of $\Pi(V)$. 
 
\begin{conj}\label{chrisconj} 
The representation $\Pi^0(V)/p$ is of finite length  
(so that $\overline{\Pi}(V)$ is well-defined) and 
$\overline{\Pi}(V)$ corresponds to $V \mod{p}$ 
under the correspondence of definition \ref{corrmodp}. 
\end{conj} 

One application of this conjecture is that it is in principle  
easier to compute $\overline{\Pi}(V)$ than to compute $V \mod{p}$  
given $\dcris(V)$. Note that using \cite{BL} and the definition of $\Pi(V)$, this conjecture is easily checked to be true when $V$ is reducible. We thus only consider $V$ irreducible in the sequel.

For $a_p\in E$ with positive valuation, let $V_{k,a_p}$ be the crystalline representation with Hodge-Tate weights $(0,k-1)$ such that 
$\dcris(V_{k,a_p}^*) = D_{k, a_p} = E e_1 
\oplus E e_2$ where: 
\[ \begin{cases} \varphi(e_1) = p^{k-1} e_2 \\ 
\varphi(e_2) = -e_1 + a_p e_2  
\end{cases} 
\text{and}\quad 
{\rm Fil}^i D_{k, a_p} = \begin{cases} 
D_{k, a_p} & \text{if $i \leq 0$,} \\ 
E e_1 & \text{if $1 \leq i \leq k-1$,} \\ 
0 & \text{if $i \geq k$.} 
\end{cases} \] 
 
In the notation of the previous lectures, we have 
$V_{k,a_p}=V(\alpha,\beta)$ where  
$\alpha$ and $\beta$ are the roots 
of the polynomial $X^2- a_p X + p^{k-1} = 0$. 
The description with $k$ and $a_p$ is a bit more convenient for 
computations associated to Wach modules. 
 
It is possible to compute $\overline{\Pi}(V_{k,a_p})$  
``by hand'' for small values of $k$ and one deduces  
from \ref{chrisconj} some predictions for what 
$V_{k,a_p} \mod{p}$ should be. Indeed, we have the following theorem 
(see \cite[th\'eor\`eme 1.4]{Br2}): 
 
\begin{thm}\label{br1p4} 
Suppose that ${\rm val}(a_p)>0$.
\begin{enumerate} 
\item If $2 \leq k \leq p+1$, 
then $\overline{\Pi}(V_{k,a_p}) = \pi(k-2,0,1)$. 
\item If $k=p+2$ and $p \neq 2$, then  
\begin{enumerate} 
\item if ${\rm val}(a_p)<1$, then $\overline{\Pi}(V_{k,a_p}) = 
  \pi(p-2,0,\omega) \simeq \pi(1,0,1)$. 
\item if ${\rm val}(a_p) \geq 1$, then $\overline{\Pi}(V_{k,a_p}) = 
\pi(p-2,\lambda,\omega)^{\rm ss} \oplus 
\pi(p-2,\lambda^{-1},\omega)^{\rm ss}$ where 
$\lambda^2-\overline{(a_p/p)}\lambda+1=0$.  
\end{enumerate} 
If $p+3 \leq k \leq 2p$ and $p \neq 2$, then  
\begin{enumerate} 
\item if ${\rm val}(a_p)<1$, then $\overline{\Pi}(V_{k,a_p}) = 
  \pi(2p-k,0,\omega^{k-1-p}) \simeq \pi(k-1-p,0,1)$.  
\item  if ${\rm val}(a_p) = 1$, then $\overline{\Pi}(V_{k,a_p}) = 
\pi(k-3-p,\lambda,\omega)^{\rm ss} \oplus 
\pi(2p-k,\lambda^{-1},\omega^{k-1-p})^{\rm ss}$ where $\lambda = 
\overline{(k-1)a_p/p}$.  
\item if ${\rm val}(a_p) > 1$, then $\overline{\Pi}(V_{k,a_p}) = 
\pi(k-3-p,0,\omega)$. 
\end{enumerate} 
\end{enumerate} 
\end{thm} 
 
In the remainder of this chapter, we will prove some of the 
formulas for $V_{k,a_p} \mod{p}$ which are predicted by the above 
theorem (via Conj.\ref{chrisconj}). In particular, we will  
explain the proof of the following  
theorem: 
 
\begin{thm}\label{blz} 
\begin{enumerate} 
\item If $k \leq p+1$, then $V_{k,a_p} \mod{p} =  
{\rm ind} (\omega_2^{k-1})$.  
\item If $k=p+2$ and ${\rm val}(a_p) > 1$ then  
\[ V_{k,a_p} \mod{p} = {\rm unr}(\sqrt{-1}) \omega  
\oplus {\rm unr}(-\sqrt{-1}) \omega. \] 
\item If $p+3 \leq k \leq 2p-1$  
and ${\rm val}(a_p) > 1$ then 
$V_{k,a_p} \mod{p} =  
{\rm ind} (\omega_2^{k-1})$.  
\end{enumerate} 
\end{thm} 
 
When $a_p=0$, one can explicitly describe $V_{k,0}$ in terms of 
Lubin-Tate characters and so we can also describe $V_{k,0} \mod{p}$ 
for all $k$.  
When $k \leq p$, then the theory of Fontaine-Laffaille gives us 
an explicit description of $V_{k,a_p} \mod{p}$ regardless of the 
valuation of $a_p$ and the theorem is then straightforward. 
Our method of proof for theorem \ref{blz} is to show that  
once we fix $k$, then if we vary $a_p$ a little bit, $V_{k,a_p} \mod{p}$ 
does not change. In order to do this, we compute the Wach module  
associated to $V_{k,a_p}$ as explicitly as possible.  
 
\subsection{Lifting $\dcris(V)$} 
 
In this paragraph, we'll explain the general strategy for constructing 
the Wach module $N(V)$ attached to a crystalline representation $V$ if 
we are given the filtered $\varphi$-module $\dcris(V)$. Recall from theorem 
\ref{sdlkfj} in \S \ref{LB54} that $N(V) / X \cdot N(V)$ is a filtered 
$\varphi$-module isomorphic to $\dcris(V)$. Furthermore, the functor $V \mapsto 
\dcris(V)$ is fully faithful, so if we can construct \emph{some} Wach module  
$N$ such that $N / X \cdot N = \dcris(V)$, then we necessarily have $N=N(V)$.  
 
Note that both the $\varphi$-module structure and the filtration on  
$N/XN$ depend solely on the map $\varphi : N \ra N$ and not at all on 
the action of $\Gamma_{\Q}$. Furthermore, given an admissible filtered 
$\varphi$-module $D$, it is easy to construct  
a $\varphi$-module $N$ over $\O[[X]]$ such that $N/XN=D$.  
The problem then is that it is very hard in general to find an action 
of $\Gamma_{\Q}$ which commutes with $\varphi$.  
 
Therefore, in order to construct the Wach module $N(V)$  
starting from $\dcris(V)$,  
one has to find a ``lift'' of $\varphi$ which is such that we can then define 
an action of $\Gamma_{\Q}$ which will commute with that lift. I do not 
know of any systematic way of doing this in general. 
 
If $a_p=0$, then we can carry out this program easily enough  
and the resulting Wach module was given in \S \ref{LB51} when $k=2$. 
Let us generalize it to arbitrary $k \geq 2$.  
The Wach module $N(V_{k,0}^*)$ is of rank $2$ 
over $\bplus_{\Q}$ generated by $e_1$ and $e_2$ with 
$\varphi(e_1)=q^{k-1} e_2$ and $\varphi(e_2)= - e_1$. 
The action of $\gamma \in \Gamma_{\Q}$ 
is given by: 
\begin{align*}  
\gamma(e_1) & = \left(\frac{\log^+(1+X)}{\gamma(\log^+(1+X))} 
\right)^{k-1} e_1 \\ 
\gamma(e_2) & = \left(\frac{\log^-(1+X)}{\gamma(\log^-(1+X))}  
\right)^{k-1} e_2  
\end{align*}  
where \[ \log^+(1+X) =\prod_{n \geq 0} \frac{\varphi^{2n+1}(q)}{p}  
\quad\text{and}\quad 
\log^-(1+X) =\prod_{n \geq 0} \frac{\varphi^{2n}(q)}{p},  \] 
as previously.
 
\subsection{Continuity of the Wach module} 
One can also think about the above constructions 
in terms of matrices. Given a  
Wach module $N$, we can fix a basis and let $P \in   
{\rm M}(d,\bplus_{\Q})$ be the matrix of Frobenius 
in that basis. The problem is then to construct an  
action of $\Gamma_{\Q}$, that  
is for every $\gamma \in \Gamma_{\Q}$,  
we need to construct a matrix $G \in  
{\rm GL}(d,\bplus_{\Q})$ such that $P \varphi(G) = G \gamma(P)$. 
In practice we need do this only for a topological  
generator $\gamma$ of $\Gamma_{\Q}$ (which is procyclic for $p \neq 
2$).   
 
In this paragraph, we use this approach to show  
that if we know a Wach module $N$ and we ``change'' $P$ slightly, 
then this defines another Wach module. More to the point, we  
have the following proposition: 
 
\begin{prop}\label{wachperturb} 
Let $N$ be a Wach module and let $P$ and $G$ be the matrices of $\varphi$ 
and $\gamma \in \Gamma_{\Q}$. Define $\alpha(r)=v_p (1-\eps(\gamma)) \cdot  
(1-\eps(\gamma)^2) \cdots (1-\eps(\gamma)^r)$. If $H_0 \in {\rm M} 
(d, p^{\alpha(k-1)+C} \aplus_{\Q})$ with $C \geq 0$, then there exists  
$H \in {\rm M} (d, p^C \aplus_{\Q})$ such that if we set $Q=({\rm 
Id}+H)P$, then the $\varphi$-module defined by the matrix $Q$ admits a 
commuting action of $\Gamma_{\Q}$. 
\end{prop} 
 
The proof of this rests on the following technical lemma: 
 
\begin{lem}\label{wachclose} 
If $N$ is a Wach module and if $G \in {\rm GL}(d,\aplus_{\Q})$  
is the matrix of $\gamma \in \Gamma_{\Q}$, and if $H_0 \in {\rm M} 
(d, p^{\alpha(k-1)+C} \aplus_{\Q})$ with $C \geq 0$,  
then there exists $H \in {\rm M} 
(d, p^C \aplus_{\Q})$ such that: 
\begin{enumerate} 
\item $H \equiv H_0 \mod{X}$. 
\item $({\rm Id}+H)^{-1} G \gamma({\rm Id}+H) \equiv G \mod{X^{k-1}}$.  
\end{enumerate} 
\end{lem} 
 
The idea of the proof of proposition \ref{wachperturb} 
is then that if $P,G$ are the matrices of $\varphi$ 
and $\gamma$ on $N$, then the lemma above and 
the fact that $P \varphi(G) = G \gamma(P)$ imply that  
$G - ({\rm Id}+H) P \varphi(G) \gamma(P({\rm Id}+H))^{-1} \equiv 0 
\mod{X^{k-1}}$ and one can then construct a matrix for $\gamma$ by 
successive approximations starting from $G$ (see \cite[\S 3.1]{BLZ} 
for more details). 
 
The proof of lemma \ref{wachclose} is also by successive 
approximations and we sketch it below. 
\begin{proof} 
If we write explicitly the condition (2) of the proposition,  
then we get: 
\begin{multline*} 
(\on{Id} + X G_1 + X^2 G_2 + \cdots )\gamma(H_0 + X H_1 + \cdots)  
=  (H_0 + X H_1 + \cdots) (\on{Id} + X G_1 + X^2 G_2 + \cdots )  
\end{multline*} 
and bearing in mind that $\gamma(X^r) \equiv \eps(\gamma)^r X^r 
\mod{X^{r+1}}$, we can solve the above equation $\mod{X^r}$ in $r$ 
successive steps. The first is to find $H_1$ such that 
$(1-\eps(\gamma))H_1 = G_1 H_0 - H_0 G_1$ so that \[ H_1 \in {\rm M} 
(d, p^{\alpha(k-1)+C-v_p(1-\eps(\gamma))} \aplus_{\Q}). \]  
The second is 
to find $H_2$ such that $(1-\eps(\gamma)^2)H_2 =$ (some $\Z$-linear 
combination of products of $G_0, G_1, G_2, H_0, H_1$) so that  \[ H_2   
\in {\rm M} (d, p^{\alpha(k-1)+C-v_p(1-\eps(\gamma))  
- v_p(1-\eps(\gamma)^2)} \aplus_{\Q}). \] Eventually, we get 
$H_0,\cdots, H_{k-1} \in {\rm M} (d, p^C \aplus_{\Q})$  
and we can then set 
$H=H_0+XH_1+\cdots+X^{k-1}H_{k-1}$. 
\end{proof} 
 
We can then apply proposition \ref{wachperturb} to the matrices 
\[ P = \begin{pmatrix} 
0 & -1 \\ q^{k-1} & 0  
\end{pmatrix}  
\qquad\text{and}\qquad 
H_0 = \begin{pmatrix} 
0 & 0 \\ -a_p & 0  
\end{pmatrix} \] 
to get a Wach module $N_{k,a_p}$ as soon as the valuation of $a_p$ is 
large enough, and a straightforward computation shows that 
$N_{k,a_p}/X N_{k,a_p}$ is then isomorphic to $\dcris(V_{k,a_p}^*)$. 
 
An explicit computation of $\alpha(k-1)$ and some refinements of the 
computations above specific to this case then give us the 
proof of theorem \ref{blz}. 
 
\begin{rem} 
To conclude, let us point out that using the description of $\Pi(V)^*$ 
with $(\varphi,\Gamma)$-modules given in the previous lectures, 
we can prove that $\Pi^0(V)/p$ is of 
finite length in general and we can also prove conjecture 
\ref{chrisconj} in many cases, including all the cases when $T/pT$ is 
irreducible (forthcoming work).  
\end{rem}

\end{document}